\documentclass[a4paper]{article}
\usepackage{mathtools}
\usepackage{amsfonts}
\usepackage{amssymb}
\usepackage{amsthm}
\usepackage{algorithm}
\usepackage{algpseudocode}
\usepackage{graphicx}
\usepackage{float}
\usepackage[pdftitle={Finite-state transducers for substitution tilings},
            pdfauthor={Simon Tatham},hidelinks]{hyperref}
\usepackage[format=plain,justification=centering]{caption}

\newcommand{\N}{\mathbb{N}}
\newcommand{\R}{\mathbb{R}}

\newcommand{\Int}[2]{\mathbf{I}_{#1,#2}}
\newcommand{\Ext}[2]{\mathbf{E}_{#1,#2}}
\newcommand{\Tile}[1]{\bigl(#1\bigr)}
\newcommand{\Edge}[2]{\bigl(#1,#2\bigr)}
\newcommand{\Sup}[2]{\bigl(\xleftarrow{#1}#2\bigr)}
\newcommand{\HTPF}{\textit{HTPF}}
\newcommand{\HHTPFFF}{\textit{HHTPFFF}}
\renewcommand{\leq}{\leqslant}
\renewcommand{\geq}{\geqslant}

\newtheorem{theorem}{Theorem}

\title{Finite-state transducers for substitution tilings}
\author{Simon Tatham}
\date{2025-12-18}

\begin{document}

\maketitle

\begin{abstract}
We present a suite of algorithmic techniques for handling substitution
tilings by treating a tile's hierarchy of supertiles in a purely
combinatorial fashion using finite state automata. The resulting
techniques are very convenient for practical generation of patches of
aperiodic tilings such as hats, Spectres and Penrose tiles, both
random and deliberately selected. They also permit some analyses of
the represented tiling. A particular product of this process is two
substitution systems for the hat tiling which are `unambiguous' in
that a single tile address uniquely determines the rest of the plane.
\end{abstract}

\section{Introduction}

Well-known aperiodic tilings such as the Penrose P2 and P3 tilings
\cite{Penrose1979}, the Ammann-Beenker tiling (\cite{10.5555/19304},
Figures~10.4.14 and~10.4.15), the hat tiling \cite{Smith_Jul2024} and
the Spectre tiling \cite{Smith_Sep2024}, are generated by
\emph{substitution systems}, in which one tiling of the plane is
converted (`\emph{deflated}') into another, by consistently
substituting a set of smaller tiles for each original tile. In some
substitution systems, every level of this system looks the same; in
others, a set of \emph{metatiles} is repeatedly deflated, and then a
different final substitution step turns the metatiles into the tiles
of the desired end-product tiling.

In some substitution systems, the tilings before and after deflation
correspond geometrically, with the \emph{subtiles} deflated from each
\emph{supertile} occupying exactly the same region of the plane (under
an appropriate scaling), or failing that, occupying a region that
overlaps the edges of the original supertile in a consistent fashion.
In others, such as the {\HTPF} system of metatiles for the hat tiling,
the geometry is much more distorted by the deflation process, and even
the metatiles themselves change their shape slightly in every
deflation. It is therefore desirable to find ways to deal with these
tilings which do not need to consider the geometry of any tiles except
those in the final end-product tiling.

The substitution system itself provides an alternative way to identify
each tile in a tiling of the plane, via an infinitely long
\emph{address} \cite{Goodman-Strauss1999} describing the relationship
of a particular tile $t_0$ to its supertile $t_1$ (both what type of
tile $t_1$ is, and which subtile of $t_1$ is $t_0$), and the
relationship of that in turn to the second-order supertile $t_2$, and
so on. These addresses can be regarded as strings of symbols from a
finite alphabet, and \cite{Goodman-Strauss1999} suggests that finite
state automata can be used to manipulate those strings. In particular
an automaton is suggested which accepts the addresses of \emph{two}
tiles, one symbol at a time from each (as if they were combined by the
\texttt{zip} operation in programming languages such as Python), and
determines whether the two addresses describe tiles that are adjacent
in the tiling.

This article develops that idea into a full system of algorithms for
working with tilings in this way. The automaton described in
\cite{Goodman-Strauss1999} is modified slightly so that its input
includes a choice of \emph{edge} of each of the two input tiles, and
reports not only whether the tiles could be adjacent at all, but
whether they are adjacent along that specific pair of edges. A second
automaton is constructed from this, which receives only a single edge
address as input, and generates the address of its neighbour as
output. This second automaton cannot always be constructed, but if it
cannot, the substitution system can be modified to remove the
obstacle, with the modified system always being a `refinement' of the
old one -- identical except that some tile types have been duplicated
into multiple `subtypes' carrying additional information about their
role in the tiling. Algorithms for all of these constructions are
presented.

Some results of the refinement process are also described, including
two substitution systems for the hat tiling not previously published
to the best of my knowledge.

\section{Specifying a substitution system}\label{combsystem}

In order to handle a substitution tiling algorithmically, it must
first be specified, in a form that provides all the necessary data.

The algorithms described in this article work entirely on a
\emph{combinatorial} description of the tiling. The only information
required is about what tiles and edges exist, which of them are
adjacent to which, and which map to which under deflation. Lengths,
angles and position in the plane are not required. Therefore,
\emph{calculation} in real or complex numbers is also not
required:\ all these algorithms are discrete in character,
manipulating strings or sequences of elements from small finite sets.

We shall present small examples first, and work up to the full
description of a combinatorial substitution system at the end of the
section.

Section~\ref{source} describes a geometric representation of the
tiling from which the combinatorial details can conveniently be
generated. This geometric representation is not used directly by any
of the algorithms, although the proofs that the algorithms work will
depend on the combinatorial data corresponding to some geometric
tiling.

\subsection{Constraints}\label{constraints}

The term `substitution system', or `substitution tiling', has both
general and specific usages.

In a general sense, a `substitution system' might be any set of rules
that transforms one family of tiles (or perhaps other objects)
systematically into another. An example in this general sense might be
a dissection of a single polygon into multiple polygons similar to
itself, but of different sizes. Iterating this substitution leads to
the number of differently sized polygons growing without bound. One
might also choose to iterate the substitution a different number of
times for each sub-polygon, e.g.\ continuing to substitute polygons as
long as they exceed some threshold size.


Here, we are concerned with the more specific usage of `substitution
system' that applies to aperiodic tilings. For the purposes of this
article, we expect a substitution system to meet some more stringent
constraints:

\begin{description}
  \item[Every tile is substituted the same number of times.] We only
    consider systems in which an entire tiling of the plane is
    transformed into an entire new tiling, by applying the
    substitution rules to every tile exactly once. That transformation
    can be iterated, so that each tile goes through substitution $n$
    times, but we still expect that every tile has been substituted
    the same number of times.

  \item[Finitely many tile types up to congruence.] We consider
    systems with a finite set of `prototiles' (reference subsets of
    the plane) fixed in advance, so that no matter how many times the
    substitution rules are applied, every tile remains congruent (not
    merely similar) to one of those prototiles. The example above
    violates this assumption because of the unbounded set of different
    sizes of tile.

  \item[Deflation makes more tiles.] We expect that after applying the
    substitution system to a whole tiling to create a deflated tiling,
    the total number of tiles is increased. This need not be true for
    a single tile in a single deflation, but it must be true that as
    $n\to\infty$, the total number of tiles in the $n$-times deflation
    of any starting tile increases without bound.

    For example, there exists a substitution system that interchanges
    between the Penrose P2 and P3 tilings in each deflation, with the
    property that each P2 kite turns into just one P3 thick rhomb. But
    in the next deflation that P3 rhomb turns into two P2 tiles (a
    kite and a dart), and that number grows further as you iterate
    more times.
\end{description}

\subsection{Sub-tiles, sub-edges, and adjacency map}\label{adjmap}

\begin{figure}
\centering
\includegraphics*[alt={The kite and dart from the Penrose P2 tiling, each cut into two separate triangles, labelled A and B for the kite and U and V for the dart. A rule is shown for turning each triangle into two or three smaller triangles of the same types.}]{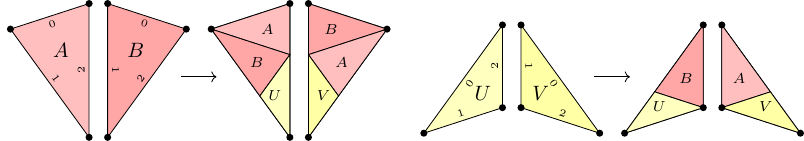}
\caption{Substitution system for Penrose P2 tiling using Robinson triangles} \label{robinson1}
\end{figure}

As an initial example, we shall use a substitution system for the
Penrose P2 tiling, in which the usual `kite' and `dart' tiles are each
subdivided into two isosceles \emph{Robinson triangles}
\cite{Robinson}. During deflation, each triangle is subdivided into
either 2 or 3 smaller triangles, according to the rules shown in
figure~\ref{robinson1}.

To describe this substitution system combinatorially, one says that
there are four \emph{types} of tile, here referred to as $A,B,U,V$.
For each tile we index its edges, in an arbitrary manner:\ here, we
adopt the convention (also shown in figure~\ref{robinson1}) that edges
are indexed 0,1,2 in anticlockwise order around the triangle, with the
base of each isosceles triangle having index 0.

During deflation, each \emph{supertile} generates a number of
\emph{subtiles}; again we must assign each subtile an arbitrary index,
and specify which tile type each subtile is. Indices can be anything
you like as long as they are distinguishable. Here we shall keep
things simple for a computer implementation, by using integers
$0,1,2,\ldots$ for the subtiles of each tile.\footnote{This Robinson
triangle system, and the corresponding one for the Penrose P3 tiling,
both have the convenient property that no two subtiles of the same
tile have the same type. As a result, you don't really need a subtile
index at all: you could reuse the subtile's type as its index. But
this property is rare and we shall not depend on it.}

Next, we must specify how the tiles arising from deflation are
arranged relative to each other, by providing an \emph{adjacency map}
for each supertile, listing pairs of subtile edges that meet in its
deflation. For example, in the deflation of the $A$ triangle, the
adjacency map says that edge 0 of the subtile of type $B$ and edge 1
of the subtile of type $U$ are adjacent to each other.

Finally, we must describe how adjacency of supertiles is related to
adjacency of subtiles. If that $A$ supertile is deflated \emph{twice},
so that its $B$ and $U$ subtiles each turn into two or three
sub-subtiles, then the adjacency map for the $B$ subtile shows how its
three sub-subtiles are arranged relative to each other, and similarly
for the $U$, but to lay out the full set of sub-subtiles of the
original $A$, we must also know which sub-subtiles of the $B$ are
adjacent to which of the $U$.

To answer this, we also specify what happens to the edges of tiles
during deflation. Each edge of a supertile is divided into a number of
\emph{sub-edges}, and we specify how many sub-edges arise from each
edge. For example, when deflating the $B$ tile, edges 1 and 2 are each
divided into two sub-edges, and edge 0 is not divided at all,
generating only a single (trivial) sub-edge.

When the edges of two supertiles meet, their sub-edges must match up,
in reversed order. For this reason it is convenient to index sub-edges
in a way that is symmetric about 0, so that sub-edge $i$ of one
supertile meets sub-edge $-i$ of the other. So when an edge is divided
into $n$ sub-edges, we assign them indices
$-n+1,-n+3,\ldots,n-3,n-1$.\footnote{Indexing the sub-edges in the
more obvious fashion, from 0 at one end, is also possible, but less
convenient. In this representation, sub-edge $i$ of one supertile
meets sub-edge $n-1-i$ of the other, so every time an algorithm needs
to compute one of those indices from the other, it must look up the
value of $n$ based on the tile type and edge index in question. The
symmetric indexing approach allows this common query to be performed
by simply negating $i$, with no auxiliary table lookup.} In this
example, every edge either has a single sub-edge with index $0$, or
two sub-edges with indices $-1,+1$.

Finally, we extend the adjacency map for each supertile so that
it also lists edge pairs consisting of one edge of a sub-tile, and one
sub-edge of a supertile edge. Formally, we use a symbol of the form
$\Int{i}{j}$ to denote edge $j$ of sub-tile $i$ of the supertile, and
$\Ext{u}{v}$ to denote sub-edge $v$ of edge $u$ of the supertile
(connoting `interior' and `exterior' edges respectively). Then the
adjacency map specifies an involution on the set of all of these
symbols, with no fixed point.

\begin{figure}
\centering
\includegraphics*[alt={The deflation rule for the A triangle from the previous figure, but this time every interior and exterior edge is marked with an identifying symbol. A table on the right shows the pairs of symbols that correspond to opposite sites of the same edge.}]{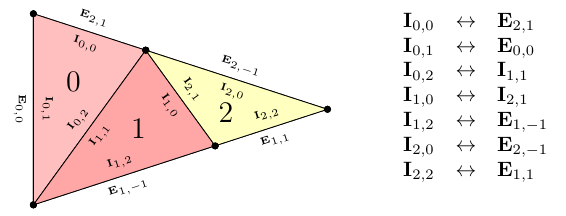}
\caption{Adjacency map for the type-$A$ Robinson triangle} \label{robinson2}
\end{figure}

Figure~\ref{robinson2} shows an example:\ the full adjacency map for
the deflation of a triangle of type $A$, and the geometric diagram it
is derived from. The other adjacency maps are calculated similarly.
The full set can be found in Figure~\ref{p2-triangles-adjmap}.

\subsection{Spurs}\label{spurs}

As a second example, we show an alternative substitution system for
the Penrose P2 tiling, this one based on whole kites and darts rather
than Robinson triangles. This system can be derived from the previous
one by first merging the deflation diagrams of the $A$ and $B$
triangles to obtain the deflation of a whole kite, and similarly
merging the $U$ and $V$ diagrams to deflate a whole dart. After this,
every $A$ and $B$ sub-tile in either deflation diagram appears next to
its counterpart, so we merge those into whole kites, obtaining a
deflation of a kite and a dart into kites, $U$ triangles, and $V$
triangles. Finally, we replace each $U$ triangle with a whole dart
sticking out past the edge of the supertile's outline, and to
compensate, we remove the $V$ triangles completely. (As if $U$ had
bulged along the line connecting it to $V$, and $V$ had shrunk
correspondingly, until $U$ had occupied all of the dart and $V$ had
reached zero size.)

\begin{figure}
\centering
\includegraphics*[alt={The kite and dart from the Penrose P2 tiling, this time whole rather than cut into triangles. For each tile type a rule is shown for decomposing it into two or three smaller tiles, not exactly matching the original tile outline. In particular, a zero-thickness spur appears in the deflation diagram for the dart, where one of the concave edges of the original dart was.}]{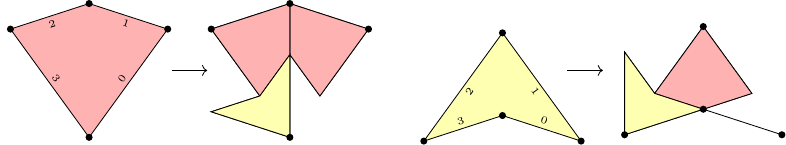}
\caption{Substitution system for Penrose P2 whole tiles}\label{p2whole}
\end{figure}

The resulting substitution system is shown in figure~\ref{p2whole}. In
this system, the deflation of a tile no longer occupies precisely the
same geometric space as the original tile. Some supertile edges, under
deflation, become crooked paths of three sub-edges with angles between
them. However, the system is still consistent:\ any two edges which can
appear next to each other in a full Penrose tiling of the plane have
matching shapes in their deflations. For example, when two kites meet
along a long edge, sharing their apex, their deflations fit too,
because the protrusion from one of the long edges fits exactly into
the hole in the other.

The new feature in this system is a zero-thickness \emph{spur},
sticking out from the right of the dart's deflation, running along a
former edge of the vanished $V$ triangle. This spur is invisible in
the final deflated tiling:\ it protrudes from one supertile's deflation
between the deflation of two others, and contributes no subtiles to
the output. But for the purposes of the algorithms presented here, a
spur such as this cannot be disregarded completely, because the
deflation of edge 1 of a dart must still have the same number of
sub-edges as that of any edge it can meet (such as edge 2 of a dart,
or edge 3 of a kite), so that every sub-edge has its counterpart in
the other supertile.

We represent this in the combinatorial description of the tiling by
allowing the adjacency map to contain entries which pair up two
\emph{exterior} edges (symbols of type $\Ext{u}{v}$), along with the
entries already seen which pair up two interior edges, or one interior
to one exterior edge. For example, in this diagram, the `upper' side
of the spur is the first of three sub-edges of the dart's edge 1.
Since those three sub-edges are indexed symmetrically as $-2,0,+2$,
the spur sub-edge will be $\Ext{1}{-2}$ in our notation. The `lower'
side of the spur is the single sub-edge of the dart's edge 0. So the
adjacency map for the dart will contain an entry pairing up
$\Ext{1}{-2}\leftrightarrow\Ext{0}{0}$. The full adjacency maps for
this system, including this pairing of two external edges, can be seen
in Figure~\ref{p2-whole-adjmap}.

\subsection{Edges are not always polygon edges}\label{edgesubdivision}

A background assumption made by this representation of a tiling is
that, whenever two tiles $t,u$ are adjacent in the tiling, they meet
`edge to edge':\ along the shared boundary, vertices of $t$ correspond
exactly to vertices of $u$, and each edge of tile $t$ on the boundary
meets the whole of a single edge of tile $u$. If this were not true,
the adjacency map described in previous sections could not be written
down in the form of a one-to-one mapping, because some edge
$\Int{t}{k}$ would not meet another single edge, but might instead meet
two or more other edges, or \emph{parts} of two other edges.

\begin{figure}
\centering
\includegraphics*[alt={Two diagrams of hat monotiles. In one, the longest edge of one hat meets the same edge of a second hat. In the other, the long edge of one hat meets two other hats, with all three hats sharing a vertex at the centre of the long edge.}]{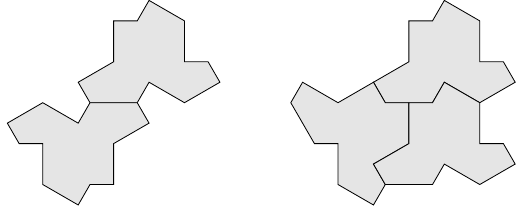}
\caption{The long edge of a hat monotile can meet one hat or two}\label{hats}
\end{figure}

If `edge' is taken to mean an edge of a polygon in the usual sense --
a maximal straight segment of the shape's boundary -- then this
property is not always satisfied. The hat monotile
\cite{Smith_Jul2024}, regarded as a polygon, has 13 edges, consisting
of six of length $1$, six of length $\sqrt3$, and one of length $2$.
When it tiles the plane, the length-$2$ edge of the hat sometimes
meets the length-$2$ edge of a neighbouring hat, but more often meets
length-$1$ edges of two different hats, so that the midpoint of the
length-$2$ edge is a point where three tiles meet. Figure~\ref{hats}
shows both of these possibilities.

Fortunately, the solution is simple. For the purposes of describing a
tiling in this combinatorial form, we regard the length-$2$ edge of
the hat as two separate length-$1$ edges with a vertex in between
them, so that the hat has 14 vertices and 14 edges, and it so happens
that two consecutive edges are collinear.

In general, for the purposes of a combinatorial description of this
kind, one chooses how to divide the boundary of each tile into `edges'
by placing a vertex at any point on the boundary where it is possible
for two other tiles to meet, and then propagating those vertices
transitively to other tiles by placing a vertex anywhere that a
neighbouring tile's vertex can appear. Once this procedure has
iterated to closure, the edges can be taken to be maximal segments of
each tile's boundary not containing a vertex. (Additional vertices can
be placed if desired for any other reason.)

\subsection{Layers}\label{layers}

Some substitution tilings use an auxiliary system of `metatiles',
which do not appear in the final tiling of the plane. One way to
generate the hat tiling \cite{Smith_Jul2024} uses a system of
metatiles labelled $H,T,P,F$, with a set of substitution rules to
deflate metatiles into other metatiles, and another set to deflate
metatiles into hats. So one can generate a large patch of hat tiles by
starting with a single metatile, say $H$; deflating it to a larger and
larger patch of smaller sub-metatiles, for as many iterations as
desired; and finally applying the other set of substitution rules just
once, to turn each of the smallest sub-sub-$\cdots$-sub-metatiles into
a small number of actual hats.

To express a tiling of this type, we introduce the concept of a
\emph{layer}, defining a subset of the tiles involved in the entire
substitution system. Each complete tiling of the plane, between
deflations, contains tiles from only one layer. For example, in the
{\HTPF} hat system, there are two layers:\ a `hats' layer containing
the hat and its reflection, and a `metatiles' layer containing the
$H,T,P,F$ metatiles.

To define the relationship between layers, we say that the final
tiling of the plane we are trying to generate is the \emph{base layer}
of the system. Then, for each layer $L$, we specify a \emph{parent
layer} $P$, containing whatever set of supertiles deflate to the tiles
in $L$. For example, in the {\HTPF} system, the parent of the `hats'
layer is the `metatiles' layer, and the parent of the `metatiles'
layer is \emph{also} the `metatiles' layer.

These parenthood relationships inductively define a sequence of
higher- and higher-order layers, indexed by $\N_0$:\ we define the 0th
layer to be the base layer, and the $(n+1)$th layer to be the parent
of the $n$th. In general this sequence will be eventually
periodic:\ with a finite number of layers, a layer must be repeated,
and then the same sequence of parents appears for ever.

However, the period can be greater than 1:\ the layer structure of a
substitution system need not terminate in a layer which is its own
parent. An example, already mentioned in section~\ref{constraints}, is
that the Penrose P2 tiling (kites and darts) can be deflated to the P3
tiling (thin and thick rhombs), and \emph{vice versa}, creating a
two-layer substitution system in which each layer's parent is the
other layer, and each deflation swaps between a P2 and P3 tiling of
the plane. (Some details of this are given in \cite{DEBRUIJN198139}.)

\subsection{Complete description of a combinatorial system}\label{combcomplete}

To summarise, here is the full set of information required to describe
a tiling substitution system in the combinatorial manner described
here.

\begin{itemize}
  \setlength\itemsep{0em}\setlength\parskip{0em}
\item \textbf{Tiles}. List the tile types that exist. For each one,
  state how many edges it has.
\item \textbf{Layers}. List the layers that exist. For each one,
  identify its parent layer, and list the tile types that appear in
  it.
\item \textbf{Deflations}. For each pair $(L,P)$ of a layer $L$ and
  its parent layer $P$, give information about how each tile type $t$
  in $P$ deflates to tiles in $L$. Specifically:
  \begin{itemize}
    \setlength\itemsep{0em}\setlength\parskip{0em}
  \item \textbf{Number of sub-edges}. For each edge of $t$, state how
    many sub-edges it turns into during deflation from $P\to L$.
  \item \textbf{Set of sub-tiles}. List the sub-tiles of $t$ when
    deflated from $P\to L$. Assign each one a distinct index, and give
    its tile type.
  \item \textbf{Adjacency map}. As described in section~\ref{adjmap},
    give the adjacency map for the deflation of $t$ from $P\to L$,
    mapping symbols of the forms $\Int{t}{e}$ and $\Ext{e}{s}$ to each
    other in pairs.
  \end{itemize}
\end{itemize}

\subsection{Geometric `source code'}\label{source}

All the algorithms described in this article can work using only the
combinatorial representation of a tiling described above. However,
that representation is not convenient for data entry. For more
complicated substitution systems, when tiles have many edges and
deflation of a tile creates many sub-tiles, the adjacency maps can be
large and unwieldy (see Figures~\ref{hath8} and~\ref{hath7} for
examples!). Also, since so much information about tile shape and
position is missing from the combinatorial data, there is very little
chance to cross-check it automatically to catch human error.

In practice the author has generally found it more convenient to
describe the tiling by its geometric details, and have a computer
program calculate the combinatorial data from that. However, this is
not \emph{required}:\ if the data listed in section~\ref{combcomplete}
can be obtained by some other means, that is just as acceptable.

To describe a tiling geometrically for conversion into this
combinatorial form, we introduce the extra concept of classifying the
tile edges into \emph{edge types}. Any two tile edges that meet in a
tiling must have the same type; any two edges that behave differently
under deflation must have different types.

Edge types are often \emph{directed}, by arbitrarily nominating one of
their end points to be a `source' and the other a `target'. This is
required if any deflation of the edge is asymmetric, and introduces
the extra constraint that when two edges of a directed type meet in
the tiling, the direction must match as well as the type, with the
target vertices of the two tiles' edges coinciding. Alternatively, an
edge type can be undirected, in which case all its deflations must be
symmetric under $180^\circ$ rotation, including types and directions
of the sub-edges.

For each edge type, the geometric tiling description specifies a
vector in $\R^2$, or (equivalently and more conveniently) a complex
number. The description of each tile states the type and direction of
every edge, and a rotation of that edge relative to its reference
orientation. Finally, for each deflation from layer $P\to L$, the
description specifies not just \emph{how many} sub-edges each edge
type deflates to, but what all their types and directions are, and by
what angle each sub-edge is rotated relative to the original edge.

This is enough information to determine the precise geometric outline
of the deflation of a tile $t$ from layer $P\to L$:\ for each edge of
$t$, look up its sequence of sub-edges, and concatenate those
sequences, each transformed by the same rotation that was specified
for that edge of the original tile.\footnote{In some substitution
systems, such as the one for the Spectre tiling given in
section~\ref{spectre}, a particular deflation $P\to L$ can be
\emph{reflecting}, meaning that after constructing each deflated tile
outline in this manner it is reflected before filling it with
subtiles.} In particular, this makes it possible to automatically
identify pairs of edges in the deflated tile outline which retrace
each other to form the zero-thickness spurs described in
section~\ref{spurs}.

After discounting those spurs, the remaining outline of the tile
deflation encloses a region of the plane which must be filled with
sub-tiles from the layer $L$. Each tile can be placed by specifying
just one entry from the adjacency map. For example, one might place
the first sub-tile, with index $t$, by saying that its edge
$\Int{t}{i}$ matches some $\Ext{u}{v}$ identifying a boundary edge. That
is enough information to locate the tile in the plane, and any other
pairs of matching edges can then be detected automatically. Each
further sub-tile is placed similarly, either adjacent to another
boundary edge, or adjacent to an edge of an already placed sub-tile.

A full example of this process is shown in a later section.
Figure~\ref{htpf-prototiles} shows the six tiles involved in the
{\HTPF} system for the hat tiling (the four metatiles, and the two
handednesses of hat), with the metatiles using five distinct edge
types and the hats two. Figures~\ref{htpf-edge-meta-meta}
and~\ref{htpf-edge-meta-hats} show the edge deflation rules for the
processes of deflating the metatiles into more metatiles, or to hats,
respectively. Finally, Figures~\ref{htpf-adjmap-meta-meta}
and~\ref{htpf-adjmap-meta-hats} show the result of deflating each
tile's outline according to each set of those rules, and the correct
tiling of the interior of the patch with metatiles or hats
respectively. (The edge types and directions are not shown, but if
they were, the edge types in the patch outlines would correctly match
those of the subtiles filling them.) Several of these diagrams show
zero-thickness spurs, such as the one in the top left of the $H\to H$
deflation which matches $\Ext{2}{4}\leftrightarrow\Ext{3}{-2}$.

This geometric representation is more convenient for data entry than
specifying the full combinatorial adjacency data directly. The user
need only specify one entry of the adjacency map per sub-tile, and all
the others are calculated automatically, by finding pairs of edges
sharing both endpoints\footnote{\label{adjmultiple}It is possible to
find that a deflation diagram includes a pair of points $u,v$ such
that there are \emph{multiple} edges $u\to v$, making the adjacency
map non-unique. Examples of this are shown in Figures~\ref{spectreh7}
and~\ref{hath7}. In this situation one expects the same number
of edges from $u\to v$ as there are from $v\to u$, and then pairs them
up arbitrarily. The different choices will only affect questions that
do not change the output tiles, such as which side of one
zero-thickness spur another one is considered to lie.}, and matching
up each edge from $u\to v$ with one from $v\to u$.

Moreover, the system of edge types permits a large amount of automatic
cross-checking:\ \emph{every} entry in the adjacency map, whether
manually specified or automatically calculated, can be checked to
ensure it really does match two edges of the same type and direction.
For this reason, it is usually convenient to treat the geometric
representation as `source code' -- a concise high-level description
entered manually by a human -- and generate the combinatorial
representation by a process strongly analogous to software
compilation, in that it first performs a thorough check for type
mismatches and other inconsistencies, and then generates all the
low-level details that were too tedious to write by hand.

Of course, one other advantage of specifying the geometric shapes of
the tiles is that \emph{after} the algorithms in this article have
been used to determine the structure of a patch of tiling, the
geometric information can be reused to draw the tiling as output!

\subsection{Addresses of tiles and edges}\label{address}

In a plane tiling generated by a substitution system, we can assign to
any tile an \emph{address} \cite{Goodman-Strauss1999} which gives its
position within the combinatorial tile hierarchy.

If the original tiling $\mathcal{T}_0$ is deflated from some tiling of
supertiles $\mathcal{T}_1$, then one can imagine identifying a tile
$t_0\in\mathcal{T}_0$ by (inductively) identifying its supertile
$t_1\in\mathcal{T}_1$, and stating which subtile of $t_1$ is $t_0$.
Our combinatorial data has assigned each of $t_1$'s subtiles a
distinct index. So this suggests it might only be necessary to specify
that index, so that the full address of $t_0$ consists of a sequence
of indices:\ the index of $t_0$ as a subtile of $t_1$, the index of
$t_1$ as a subtile of its supertile $t_2$, and so on.

However, this is not necessarily enough information to even identify
the \emph{type} of the original tile $t_0$. For example, in the system
of Penrose tiles shown in figure~\ref{p2whole}, suppose the subtiles
are indexed so that subtile 1 of a kite is one of the sub-kites, and
subtile 1 of a dart is the sub-dart. Then the index sequence
$1,1,1,\ldots$ could represent a kite all of whose supertiles are
kites, \emph{or} a dart all of whose supertiles are darts. An
ambiguous sequence of indices of this kind could potentially begin at
any point in the sequence, even if the lower-order indices in the
sequence were unambiguous.

Therefore, we must also specify the type of the starting tile and each
of its supertiles. So the full address of a tile is a sequence of the
form
$$t_0 \xleftarrow{i_1} t_1 \xleftarrow{i_2} t_2 \xleftarrow{i_3} t_3 \cdots$$
in which the $(t_n)$ are tile types, and the $(i_n)$ are subtile
indices. The notation $t_{n-1} \xleftarrow{i_n} t_n$ means that the
$(n-1)$th- and $n$th-order supertiles have types $t_{n-1}$ and $t_n$
respectively, and that the former is the $(i_n)$th subtile of the
latter.

For a full tiling of the plane, the address is an infinite sequence,
describing larger and larger supertiles without bound. In other
situations one can also consider a \emph{finitely} long address, which
stops after specifying some particular supertile type $t_n$, and says
nothing about what larger tiling (if any) that supertile might be part
of.

We shall present algorithms to handle these addresses using finite
state automata. We must therefore represent each address in the form
of a sequence of symbols from a finite alphabet. The alphabet we
choose consists mainly of symbols of the form $\Sup{i}{t}$, where $t$
is a tile type and $i$ is one of its subtile indices in some
deflation\footnote{\emph{Which} deflation? In a tiling with multiple
layers, the meaning of a symbol $\Sup{i}{t}$ can vary depending on
where in the address it appears. For example, in the {\HTPF} hat
system, if it appears as the first supertile symbol, it will describe
one of the hats in the deflation from a metatile to the base hat
layer, but if it appears any later, it will describe one of the
\emph{metatiles} in the deflation from metatiles to other metatiles.
This does not cause a problem in the algorithms, because they can
always keep track of what layer they are currently in, and
disambiguate the meaning of each symbol.}. These symbols can represent
every part of an address except for the initial tile type, for which
we require a set of extra symbols specifying \emph{only} a tile type.
So every tile address starts with a bare tile type, followed by an
infinite sequence of $\Sup{i}{t}$ symbols:
$$\Tile{t_0}, \Sup{i_1}{t_1}, \Sup{i_2}{t_2}, \Sup{i_3}{t_3}, \ldots$$

The main aim of our algorithms is to determine the address of a tile's
neighbour from that of the tile itself. Therefore, we must have a way
to specify \emph{which} neighbour is desired. To place the two tiles
relative to each other in the plane, we must also know how the
neighbour is oriented relative to the initial tile. In other words, we
want to specify a particular \emph{edge} of the original tile $t_0$,
and receive an answer telling us not only what type of tile is
adjacent to $t_0$ along that edge, but also which of its edges
coincides with the edge we specified.

Therefore, we need another kind of address that specifies a particular
edge of a tile. We modify the initial symbol so that instead of being
a plain tile type, it is a pair $\Edge{t}{e}$, indicating edge $e$ of
tile $t$ according to whatever arbitrary indexing scheme we have
assigned to our tile edges. So the address of a tile \emph{edge}
consists of one of these modified initial symbols, but thereafter,
uses the same set of supertile symbols $\Sup{i}{t}$ as the address of
a tile:
$$\Edge{t_0}{e}, \Sup{i_1}{t_1}, \Sup{i_2}{t_2}, \Sup{i_3}{t_3}, \ldots$$

\section{Relation to regular languages}\label{reglang}

We have described representations of the address of a tile, and the
address of a tile edge, in terms of sequences of symbols from a
finitely large alphabet. The aim was to calculate further addresses
within the same tiling by string-processing techniques, rather than
geometric ones.

The simplest type of string processing known to computer science is a
finite state automaton. These are reliably fast, and easy to analyse
and reason about. So if we can possibly use finite state automata to
process tile addresses, it's worth a try.

Let $\mathcal{E}$ denote the set of finite strings of symbols of the
correct form to be the edge address of a tile within an $n$th-order
supertile for some $n$. That is, strings of the form
$$\Edge{t_0}{e}, \Sup{i_1}{t_1}, \Sup{i_2}{t_2}, \ldots, \Sup{i_n}{t_n}$$
satisfying the necessary consistency criteria: $e$ must be the index
of an edge of tile type $t_0$, and each $t_{n}$ must be the correct
tile type to be subtile $i_{n+1}$ of a tile of type $t_{n+1}$.

It is easy to see that $\mathcal{E}$ is a regular language, because a
DFA can be constructed for it trivially, consisting of a special
\textsc{start} state and a state for every $(\textit{layer},
\textit{tile type})$ pair, with a transition between two states of the
latter type whenever the supertile type and subtile index in a symbol
are consistent with the source state of the transition.

Similarly, let $\mathcal{T}$ denote the language of strings giving the
address of a tile rather than an edge, identical to $\mathcal{E}$
except that the first symbol is a bare tile type $t_0$ instead of a
pair $(t_0,e)$. This too is regular, by the same reasoning.

Let $\mathcal{T}\times\mathcal{T}$ denote the set of ordered pairs of
two tile addresses of the same length, represented as a sequence of
ordered pairs containing one symbol from each address, as if they had
been combined by Python \texttt{zip}. That is, if
$(\sigma_0,\sigma_1,\ldots,\sigma_n),(\sigma'_0,\sigma'_1,\ldots,\sigma'_n)\in\mathcal{T}$,
then
$((\sigma_0,\sigma'_0),(\sigma_1,\sigma'_1),\ldots,(\sigma_n,\sigma'_n))\in\mathcal{T}\times\mathcal{T}$.
We define $\mathcal{E}\times\mathcal{E}$ similarly.

Let the \emph{neighbour language}
$\mathcal{N}\subset\mathcal{E}\times\mathcal{E}$ denote the set of
pairs of edge addresses ($A,A'$) which are neighbours within the
$n$-times deflation of some common supertile. That is, the tiles
identified by the two addresses are adjacent in the $n$th-order
deflation of the supertile, and moreover, they meet along the edges
specified by $A$ and $A'$.

\begin{theorem}\label{adjreg}
$\mathcal{N}$ is a regular language.
\end{theorem}

\begin{proof}

We begin by showing this for a restricted class of substitution
systems:\ those in which the deflation of each tile type can be scaled
to occupy the same region of the plane as the original tile, with no
zero-thickness spurs or geometric distortion. An example is the
Robinson-triangle system for the P2 tiling, shown in
Figure~\ref{robinson1}.

If a substitution system has this property under a single deflation,
then it has it under any finite number of deflations. So if two tiles
$t,u$ are adjacent within the deflation of an order-$n$ supertile,
then for every $0<i<n$, either their order-$i$ supertiles coincide, or
are themselves adjacent along a whole edge, with the edge separating
$t,u$ being part of that longer supertile edge.

Therefore, as we ascend from the tiles $(t,u)$ to their supertiles
$(t_1,u_1),\ldots,(t_n,u_n)$, each pair of supertiles is adjacent,
until the first pair that coincide. After that, all further supertile
pairs coincide.

This allows us to directly construct a state machine which matches
$\mathcal{N}$. For the initial segment of a string before two
supertiles coincide, each state $S(L,t,u,e,f)$ describes a layer $L$ of
the system; a pair of tile types $t,u$ within that layer; and a pair
of edge indices $e,f$ of those tiles, indicating their adjacency: edge
$e$ of tile $t$ meets edge $f$ of tile $u$.

Transitions can be determined using the adjacency maps, such as the
one shown in Figure~\ref{robinson2}. In a state $S(L,t,u,e,f)$, if
there exist tile types $t',u'$ in $L$'s parent layer $P$, and indices
$i,j,k,e',f'$ such that
\begin{itemize}
\setlength\itemsep{0em}\setlength\parskip{0em}
\item subtile $i$ of tile type $t'$ (when deflated from $P\to L$) is type $t$
\item subtile $j$ of tile type $u'$, similarly, is type $u$
\item the adjacency map for deflating tile $t'$ maps
$\Int{i}{e}\leftrightarrow\Ext{e'}{k}$
\item the adjacency map for deflating tile $u'$ maps
$\Int{j}{f}\leftrightarrow\Ext{f'}{-k}$
\end{itemize}
then this is saying that $t,u$ occupy compatible positions along the
boundaries of two edges of $t',u'$, such that if edge $e'$ of $t'$
meets edge $f'$ of $u'$, then their subtiles would meet in the way
described by the original state $S(L,t,u,e,f)$. Hence, a transition
exists from $S(L,t,u,e,f)\to S(P,t',u',e',f')$, on a symbol pair
$\Bigl(\Sup{i}{t'},\Sup{j}{u'}\Bigr)$.

At the point where the supertiles of $t,u$ coincide, we enter a second
family of states, describing a single tile in a particular layer, say
$A(L,t)$. These are the accepting states of our automaton. The first
transition into an accepting state occurs when subtiles $i,j$ of a
supertile type $s$ match the types $t,u$ from an existing state, and
those subtiles meet internally to the supertile's deflation, along the
correct edges. That is, if $s$'s adjacency map includes the entry
$\Int{i}{e}\leftrightarrow\Int{j}{f}$, then a transition exists from
$S(L,t,u,e,f)\to A(P,s)$ on the symbol pair
$\Bigl(\Sup{i}{s},\Sup{j}{s}\Bigr)$. Thereafter, transitions between
accepting states require both symbols in the pair to be identical, and
to continue to spell out a valid tile address, exactly as in the
natural DFA matching $\mathcal{T}$.

Now we must extend this construction to more general substitution
systems: those in which tile edges deflate to non-straight outlines,
and especially those with zero-thickness spurs. In these systems, it
need not be true that two tiles $t,u$ sharing an edge have supertiles
$t',u'$ also sharing an edge.

\begin{figure}
\centering
\includegraphics*[alt={Increasingly complex outlines derived from the P2 kite and dart shapes, in each step transforming the edges via the deflation map. After many iterations, each shape has converged to an irregular fractal.}]{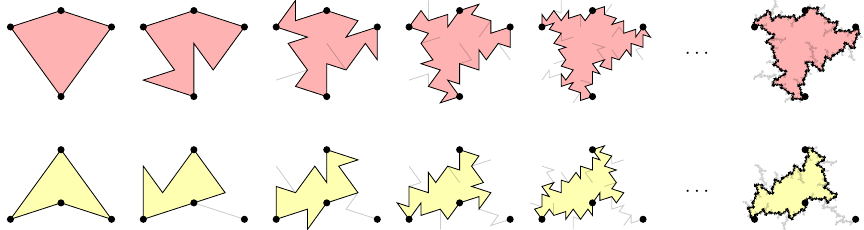}
\caption{Converging to fractal tile outlines for the P2 whole-tiles
  system}\label{p2fractal}
\end{figure}

\begin{figure}
\centering
\includegraphics*[alt={On the left, a group of six tiles is shown, consisting of the deflations of two larger kite tiles. The dart subtile of the top kite shares an edge with a kite subtile of the bottom kite, but the supertiles themselves only touch at a vertex. On the right, the same arrangement is shown with the fractal versions of the tile boundaries, and both the subtiles and supertiles share an identical segment of fractal boundary.}]{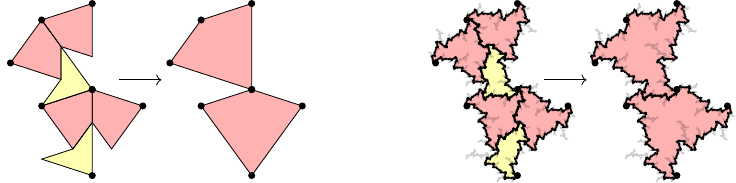}
\caption{Adjacent tiles may not have adjacent supertiles, but the fractal versions do}\label{p2fractalexample}
\end{figure}

In systems meeting the constraints described in
section~\ref{constraints}, iterating the geometric deflation process
and re-scaling each deflated tile shape to the same overall area as
the original tile causes the tile outlines to converge to a fractal.
An example, based on the P2 whole-tile substitution system in
Figure~\ref{p2whole}, is shown in Figure~\ref{p2fractal}. Transforming
the substitution system into this form, with every tile having a
fractal outline and every tile edge being a fragment of that outline,
restores the property that the deflation of a tile occupies the same
plane region as the tile itself. Hence, when two tiles are adjacent
along a segment of their outline, their two supertiles must also be:
see Figure~\ref{p2fractalexample} for an example case. So we can still
construct a set of automaton states based on two tiles sharing a
segment of boundary.

We must also show that there are \emph{finitely} many such states,
which is no longer obvious when all we know is that two tiles share
\emph{some segment of boundary} instead of one of a discrete set of
edges.

To see this, we observe that the fractal tile shapes are bounded.
In a tiling of the plane using tiles from some particular layer, there
are only finitely many ways to tile a disc of radius $r$ with a given
tile $t$ at a fixed location in the centre, and the fractal forms
of those patches exhibit all the ways that the fractal outline of
tile $t$ can be adjacent to that of any other tile. No other forms of
adjacency need be considered, because even if two fractal tile
outlines $t,u$ physically fit together in some other way, that
adjacency cannot occur between a pair of supertiles within the context
of a legal tile address, since any such pair of supertiles form part
of a tiling of the whole plane.
\end{proof}

Calculating the fractal tile shapes and their adjacencies is
potentially awkward. We shall show a further useful property of
$\mathcal{N}$, which we shall use in the next section to derive an
easier, purely combinatorial algorithm for finding a recogniser for
$\mathcal{N}$.

The extra property relates to \emph{reversibility}. A DFA that
recognises a regular language $L$ can always be turned into a
recogniser for the language $L^R$ consisting of the reversal of every
string of $L$, by simply reversing all the arrows in the state
transition diagram and interchanging the starting and accepting
states. Usually the reversed recogniser is no longer deterministic,
and must be determinised again if a DFA for $L^R$ is desired.
\cite{10.1145/322326.322334} defines the concept of a
\emph{zero-reversible acceptor} to be a DFA for a language $L$ which
remains deterministic under this reversal\footnote{Why
`zero-reversible' and not just `reversible'?
\cite{10.1145/322326.322334} defines the weaker condition of
$k$-reversibility to mean that the reversed automaton need not be
deterministic already, but if multiple transitions on the same symbol
are available then $k$ extra symbols of lookahead suffice to rule out
all but one of them. `Reversible' by itself is used to mean
`$k$-reversible for some $k\geq0$'.}, because there is a unique
accepting state, and no two distinct states have transitions to the
same state on the same symbol. A zero-reversible \emph{language} is
one that admits a zero-reversible acceptor.

Our neighbour language $\mathcal{N}$ cannot in general be
zero-reversible, because it has an accepting state for every supertile
type, i.e.\ usually more than one. However, it can be written as the
disjoint union of zero-reversible languages in a natural way.

We define the sublanguage $\mathcal{N}(L,t)\subset\mathcal{N}$ to be
the set of edge address pairs which represent adjacent tiles whose
highest-order supertiles \emph{not} shared are in layer $L$ of the
substitution system, and their lowest-order shared supertile is a tile
of type $t$ in layer $P$, the parent of $L$. Sublanguages of the form
$\mathcal{N}(L,t)$ are pairwise disjoint, and partition $\mathcal{N}$.

\begin{theorem}\label{zerorev}
$\mathcal{N}(L,t)$ is zero-reversible.
\end{theorem}
\begin{proof}
From the proof of Theorem~\ref{adjreg}, we can directly construct a
DFA for $\mathcal{N}$ in which each non-accepting state describes two
tiles $s,s'$ in layer $L$ adjacent along subsets $b,b'$ of their
boundary. (In a system with spurs, these boundary subsets may not
correspond to single edges of the original tiles, but nonetheless,
there will be a finite set of possible boundary subsets that $s,s'$
can legally share.) Accepting states describe a single tile $s$ in
layer $L$. We denote these types of state as $S(L,s,b,s',b')$ and
$A(L,s)$ respectively.

We consider $b,b'$ to be the maximal subsets of each tile's boundary
shared with the other's. In some substitution systems, the geometry of
two tiles means that if they are adjacent along some pair of edges
$e,e'$ then they are also adjacent along another edge pair $f,f'$. For
these purposes we regard $b,b'$ as being the union of both of those
edges, and any further edges the two tiles might share. Thus, we do
not consider $S(L,s,e,s',e')$ and $S(L,s,f,s',f')$ to be distinct DFA
states.

A DFA for $\mathcal{N}(L,t)$ can be constructed from the DFA for
$\mathcal{N}$ itself by removing all the accepting states except for
$A(P,t)$, and giving that state no outgoing transitions, so that the
two input tile addresses must specify just enough symbols to reach the
first shared supertile.

This DFA has a unique accepting state, by construction:\ we made it by
deleting all the other accepting states. So we only need to show that
no state has two distinct predecessors with transitions to it on the
same symbol.

A transition to the accepting state is of the form $S(L,s,b,s',b')\to
A(P,t)$, on the symbol $\Bigl(\Sup{i}{t},\Sup{i'}{t}\Bigl)$. Such a
transition is only valid if subtiles $i$ and $i'$ of $t$ have types
$s,s'$ respectively. Then subtiles $i,i'$ of $t$ have a fixed
geometric relationship, and there is only one maximal pair of boundary
segments $b,b'$ along which they can meet. So $A(P,t)$ has at most one
predecessor for any symbol.

Similarly, consider a transition between two non-accepting states
$S(M,s,b,s',b')\to S(Q,t,c,t',c')$ (where $Q$ is the parent layer of
$M$) on the symbol $\Bigl(\Sup{i}{t},\Sup{i'}{t'}\Bigr)$. Again, $s$
is constrained to be the tile type of subtile $i$ of $u$, and $s'$
similarly; and placing the fractal tile shapes of $t,t'$ next to
each other along boundary segments $c,c'$ fixes their relative
position in the plane, and leaves only one possibility for which
segments $b,b'$ of those subtiles can meet.
\end{proof}

\section{Algorithms for handling addresses}\label{algorithms}

In the previous section we showed
that the combinatorics of a tiling can be handled by finite state
automata. In this section we present concrete algorithms for
constructing two different forms of automaton. The ultimate aim is to
create a deterministic finite-state transducer which takes a tile edge
address as input, and produces its neighbour edge address as output.

\subsection{Computing neighbours recursively}\label{rec}

We first present a recursive algorithm to calculate the neighbour of a
tile. For some applications, this recursive algorithm is adequate by
itself to generate patches of tiling. It is also used to generate the
input to the algorithm in section~\ref{adjmatcher}.

As input, this algorithm requires a combinatorial tiling description,
and the address of a specific tile edge, in the form shown in
section~\ref{address}:
$$\Edge{t_0}{e}, \Sup{i_1}{t_1}, \Sup{i_2}{t_2}, \Sup{i_3}{t_3}, \ldots$$
It computes and returns an address in the same form, describing the
neighbour of the input tile along that edge, and which edge of it
meets $e$.

The algorithm begins by inspecting the deflation diagram of tile
$t_1$, one of whose tiles is $t_0$. The tiling description provides an
adjacency map describing this deflation diagram. In that map, the
symbol $\Int{i_1}{e}$ describes edge $e$ of the subtile with index
$i_1$ -- that is, exactly the edge we have been asked about. So the
first step is to look up that symbol in the map, to see what it is
paired with.

In the simple case, the result of that lookup is another interior
symbol, i.e.\ one of the form $\Int{i'_1}{e'}$. This means that the
sub-tile adjacent to the starting one is another sub-tile of the same
supertile $t_1$. So the algorithm need only replace the two low-order
symbols of the input address, by looking up the tile type $t'_0$ of
the new sub-tile, and returning the modified address
$$\Edge{t'_0}{e'}, \Sup{i'_1}{t_1}, \Sup{i_2}{t_2}, \Sup{i_3}{t_3}, \ldots$$

In the more complicated case, $\Int{i_1}{e}$ turns out to be paired
with an exterior edge $\Ext{u}{v}$, so that edge $e$ of the input tile
$t_0$ is part of the boundary of its supertile $t_1$. In that case, we
must find the neighbouring \emph{supertile} of $t_1$ -- we'll call
that $t'_1$ -- and decide which subtile of that is the answer to the
original query.

To find the neighbour of $t_1$, we invoke this entire algorithm a
second time, recursively. The input address is formed by removing the
two low-order symbols from our original input, and replacing them with
a symbol describing the supertile $t_1$, and in particular the edge of
it that contains our original query edge. The adjacency map lookup
result $\Ext{u}{v}$ says that the original query edge is part of edge
$u$ of $t_1$, in particular sub-edge $v$ of it. So we must find out
what supertile borders on edge $u$ of $t_1$, and therefore we
recursively query the address
$$\Edge{t_1}{u}, \Sup{i_2}{t_2}, \Sup{i_3}{t_3}, \Sup{i_4}{t_4}, \ldots$$
We also slightly modify the tiling description provided to the
recursive call, by changing its base layer to be the one containing
$t_1$, i.e.\ the parent of the original base layer.

The recursive invocation of the algorithm returns the full address of
an edge of another tile in the parent layer:
$$\Edge{t'_1}{u'}, \Sup{i'_2}{t'_2}, \Sup{i'_3}{t'_3}, \Sup{i'_4}{t'_4}, \ldots$$
So we know that edge $u'$ of this supertile borders on edge $u$ of the
original supertile. The smaller tile edge we really wanted to know
about was not the whole of edge $u$, but sub-edge $v$ of it. This
meets sub-edge $-v$ of the new supertile's edge (as discussed in
section~\ref{adjmap}). So to determine what sub-tile of $t'_1$ we
need, we look up the symbol $\Ext{u'}{-v}$ in the adjacency map for
$t'_1$.

If this lookup in the adjacency map returns an interior symbol
$\Int{i'_1}{e'}$, then we can terminate the algorithm by adding those
details to the address returned from the recursive call. Find out the tile type $t'_0$ of subtile $i'_1$ of $t'_1$, and return
$$\Edge{t'_0}{e'}, \Sup{i'_1}{t'_1}, \Sup{i'_2}{t'_2}, \Sup{i'_3}{t'_3}, \ldots$$

There is one remaining possibility. After a recursive call, the lookup
of $\Ext{u'}{-v}$ in the new supertile's adjacency map might return
\emph{another} exterior symbol. This occurs when the algorithm
encounters one of the zero-thickness spurs described in
section~\ref{spurs}. We have tried to step in to the deflation of tile
$t'_1$, and found that we've landed on a zero-thickness spur, not a
real tile -- so we must step straight off the other side of the spur
to find out what real tile lies beyond it.

In this situation, the algorithm does the same thing it would have
done if the \emph{original} lookup had returned an exterior symbol. It
makes another recursive call to compute the neighbour $t''_1$ of
$t'_1$, and loops round again, until it reaches a supertile
$t^{(k)}_1$ in which the adjacency map lookup finally returns an
interior symbol.

A summary of this algorithm, in pseudocode, is shown as
Algorithm~\ref{algrec}.

\begin{algorithm}[t]
\caption{Recursively compute the neighbour of an edge address}\label{algrec}
\begin{algorithmic}
\State let $\Edge{t_0}{e}, \Sup{i_1}{t_1}, \textit{address-tail} = \textit{input-address}$;
\State $\textit{lookup-result} \leftarrow \hbox{look up }\Int{i_1}{e}\hbox{ in adjacency map for }t_1$;
\While {\textit{lookup-result} is of the form $\Ext{u}{v}$}
    \State let $\textit{supertile-address} = (t_1,u),\textit{address-tail}$;
    \State recurse to find $(t'_1, u'),\textit{new-address-tail} = \hbox{neighbour of }\textit{supertile-address}$;
    \State $\textit{lookup-result} \leftarrow \hbox{look up }\Ext{u'}{-v}\hbox{ in adjacency map for }t'_1$;
    \State replace $(t_1,\textit{address-tail})\leftarrow(t'_1,\textit{new-address-tail})$;
\EndWhile
\State now expect that \textit{lookup-result} is of the form $\Int{i'_1}{e'}$;
\State let $t'_0 = \hbox{type of subtile }i'_1\hbox{ of }t_1$;
\State return $(t'_0,e'), (i'_1,t_1), \textit{address-tail}$;
\end{algorithmic}
\end{algorithm}

\subsection{Computing neighbours using a transducer}

The recursive algorithm described above is adequate for some purposes,
but in complicated substitution systems it has to recurse repeatedly
up and down its input string rewriting symbols, which makes it slower,
and also difficult to reason about. Moreover, Section~\ref{infsup}
will describe a class of tiling instances it cannot handle at all.

A finite-state transducer is a better way to calculate neighbour
addresses, if one can be constructed. In this section we present
algorithms to do so.

\subsubsection{Recogniser for neighbour edge pairs}\label{adjmatcher}

In section~\ref{reglang} we showed that the neighbour language
$\mathcal{N}$, consisting of all zipped-together pairs of tile edge
addresses of the same finite length describing two tiles sharing an
edge, is a regular language. Our first step is to construct a DFA to
recognise $\mathcal{N}$.

One way to construct this DFA would be to follow the structure of the
proof of Theorem~\ref{adjreg}:\ make a state for every possible
adjacent pair of tiles or their fractalisations, together with
\textsc{start} and the required accepting states, and fill in the
transitions by considering the geometry of the tiles.

However, using the geometry of the tiling complicates implementation
(the coordinates may well involve algebraic numbers, requiring either
elaborate exact computation or risky floating-point approximation), as
well as requiring more information as input. In this section we
instead present a purely combinatorial algorithm to construct the same
DFA from only the data described in section~\ref{combcomplete}, based
on Theorem~\ref{zerorev} which tells us that $\mathcal{N}$ is the
union of a set of zero-reversible languages.

In \cite{10.1145/322326.322334} an algorithm is given to infer a DFA
for a zero-reversible language from a finite corpus of strings of the
language. Paraphrased, one first constructs an NFA consisting of a
single \textsc{start} state, an \textsc{accept} state, and a separate
path from \textsc{start} to \textsc{accept} for each string in the
corpus. Then merge\footnote{The description in
\cite{10.1145/322326.322334} suggests constructing the initial DFA in
the form of a trie, so that states with the same predecessor start off
already merged. This is equivalent.} any two states $A,B$ which have
transitions $U\xrightarrow{s}A,U\xrightarrow{s}B$ from the same
predecessor $U$ on the same symbol $s$, and conversely, merge any two
states with transitions to the same \emph{successor} on the same
symbol, i.e.\ $A\xrightarrow{s}V,B\xrightarrow{s}V$ for some $s,V$.
Iterate both of those merging operations to closure. The first policy
guarantees to make the NFA deterministic; the second makes its
reversal deterministic too, so that the resulting DFA is still a DFA
if the arrows are reversed and \textsc{start} swapped with
\textsc{accept}.

\begin{algorithm}[t]
\caption{Build a DFA for neighbouring edge address pairs}\label{algadj}
\begin{algorithmic}
\State make an empty automaton with just a \textsc{start} state;
\State add all accepting states $A(L,t)$ and their transitions;
\State $\mathit{done}\leftarrow\emptyset$;\Comment{set of edge address prefixes already processed}
\For{$k$ in $1,2,\ldots$}\Comment{length of strings in input corpus}
  \For{$A_1$ in edge addresses of length $k$ with no prefix in \textit{done}}
    \State let $A_2$ = neighbour of $A_1$ according to Algorithm~\ref{algrec};
    \If{Algorithm~\ref{algrec} did not fail by overrunning the end of $A_1$}
      \State call $\textsc{AddString}(\mathit{automaton},\mathrm{zip}(A_1,A_2))$;\Comment{see Algorithm~\ref{algadjaddstring}}
      \State $\mathit{done}\leftarrow\mathit{done}\cup\{A_1\}$;
    \EndIf
  \EndFor
  \State call $\textsc{MergeStates}(\mathit{automaton})$;\Comment{see Algorithm~\ref{algadjmerge}}
  \State return if $\textsc{IsComplete}(\mathit{automaton})$;\Comment{see Algorithm~\ref{algadjcomplete}}
\EndFor
\end{algorithmic}
\end{algorithm}

In the process, the DFA is transformed from one that matched
\emph{only} the strings in the finite input corpus into one that
matches their `zero-reversible closure', the smallest zero-reversible
language containing the input corpus. For example, if the input corpus
was $\{abc,abdbc\}$, then its zero-reversible closure would be the
language $ab(db)^*c$, because the state reached via $ab$ and the one
reached via $abdb$ have been merged, permitting any number of
repetitions of $db$ to cycle back to that state before finishing the
string with a $c$.

The recursive algorithm in the previous section can be used to make a
corpus of strings from $\mathcal{N}$, by feeding it every possible
input string $A$ up to some chosen length $k$, and if it successfully
returns an output string $B$ at all, adding $\mathrm{zip}(A,B)$ to the
corpus. This algorithm can fail, by recursing so deeply that it finds
itself needing to read a symbol from beyond the end of $A$. For this
application, we must detect that error, and handle it by abandoning
the attempt to determine the neighbour address of $A$ and moving on to
the next possibility.

Given a corpus of this kind, it is straightforward to classify each
string into one of the sublanguages $\mathcal{N}(L,t)$ described in
section~\ref{reglang}, partitioning the corpus into sub-corpuses. Each
of these sub-corpuses is a sample from a zero-reversible language, so
the state-merging algorithm described above can be used to make a DFA
for it. Then those DFAs can be recombined into one with multiple
accepting states.

\begin{algorithm}[t]
\caption{Subroutine to add new states to an NFA to make it recognise
  one extra string, labelling the states with layers of the
  substitution system}\label{algadjaddstring}
\begin{algorithmic}
\Procedure{AddString}{\textit{automaton},\textit{string}}
\State let $(p_1,p_2,\ldots,p_k)=\mathit{string}$;\Comment{each $p_i$ is a symbol pair}
\State add states $s_1,\ldots,s_{k-1}$ to automaton;
\State set $\mathrm{layer}(s_1)=\mathrm{parent}(B)$, $\mathrm{layer}(s_{i+1})=\mathrm{parent}(\mathrm{layer}(s_i))$;
\State let $L=\mathrm{layer}(s_{k-1})$, $t=\hbox{tile type specified in both symbols of }p_k$;
\State add transitions $\textsc{start}\xrightarrow{p_1}s_1\xrightarrow{p_2}\cdots\xrightarrow{p_{k-1}}s_{k-1}\xrightarrow{p_k}A(L,t)$;
\EndProcedure
\end{algorithmic}
\end{algorithm}

\begin{algorithm}[t]
\caption{Subroutine to merge states in an NFA, turning it into a DFA
  matching a union of zero-reversible languages}\label{algadjmerge}
\begin{algorithmic}
\Procedure{MergeStates}{\textit{automaton}}
  \Repeat
    \If{$\exists$ transitions $A\xrightarrow{p}C,B\xrightarrow{p}C$ and $\mathrm{layer}(A)=\mathrm{layer}(B)$}
      \State merge states $A,B$ and all their transitions;
    \EndIf
    \If{$\exists$ transitions $A\xrightarrow{p}B,A\xrightarrow{p}C$ and $\mathrm{layer}(B)=\mathrm{layer}(C)$}
      \State merge states $B,C$ and all their transitions;
    \EndIf
  \Until{neither type of merge is possible}
\EndProcedure
\end{algorithmic}
\end{algorithm}

Rather than literally separating the input corpus into subsets, it is
simpler to build the whole combined DFA in one go, by making a single
NFA with one \textsc{start} state and a full set of $A(L,t)$ accepting
states. Then, as before, create a separate path of states from
\textsc{start} to an appropriate accepting state for each string of
the input corpus. Inductively annotate every state with a layer of the
substitution system:\ \textsc{start} is associated with the base layer,
and any successor of a state $s$ is associated with the parent of $s$'s
layer. Now repeatedly merge states which have the same predecessor on
the same symbol, or which have the same successor on the same symbol,
\emph{except} that two states may only be merged if they are annotated
with the same layer. This precaution avoids the mistake of merging
states from different layers $L,L'$ with transitions to the same
accepting state $A(P,t)$, where $P$ is the common parent of both $L$
and $L'$. (Recall from Theorem~\ref{zerorev} that the sublanguages
$\mathcal{N}(L,t)$ are distinguished by the last layer \emph{before}
entering an accepting state, so such a merge would mix up the states
from two of those sublanguages' DFAs, and create nonsense.)

This procedure terminates in a DFA recognising some sublanguage of the
full language $\mathcal{N}$. (Specifically, the union of the
zero-reversible closure of each sub-corpus; each of those closures is
a sublanguage of one $\mathcal{N}(L,t)$.) However, it may not be all
of $\mathcal{N}$.

\begin{algorithm}[t]
\caption{Subroutine to check an automaton has a neighbour for every legal input}\label{algadjcomplete}
\begin{algorithmic}
\Function{IsComplete}{$A: \mathit{automaton}$}
  \State let $B$ be a DFA matching any valid edge address;
  \State $\mathit{discovered}\leftarrow\{(\textsc{start}_A,\textsc{start}_B)\}$;
  \State $\mathit{visited}\leftarrow\emptyset$;
  \While{$\mathit{discovered}\setminus\mathit{visited}$ contains a state pair $(s_A,s_B)$}
    \State add $(s_A,s_B)$ to $\mathit{visited}$;
    \For{\textbf{all} $(b,t_B)$ such that $B$ has a transition $s_B\xrightarrow{b}t_B$}
      \If{$\not\exists (a,t_A):A\hbox{ has a transition }s_A\xrightarrow{(a,b)}t_A$}
        \State return \textbf{false};\Comment{automaton is incomplete!}
      \Else
        \State add ($t_A,t_B$) to $\mathit{discovered}$ for all such pairs $(a,t_A)$;
      \EndIf
    \EndFor
  \EndWhile
  \State return \textbf{true};
\EndFunction
\end{algorithmic}
\end{algorithm}

By Theorem~\ref{adjreg}, the whole neighbour language $\mathcal{N}$ is
regular. Therefore, there is some finite length $k$ such that every
state of a minimal DFA for $\mathcal{N}$ is reached by some edge
address pair of length $\leq k$. However, this $k$ is not known in
advance, and may be greater than the length of strings we generated
for our input corpus. So the final part of our algorithm is to detect
if the output DFA is complete, and if not, retry with a corpus of
longer address pairs. Since $\mathcal{N}$ is regular, sooner or later,
an attempt will succeed.

To detect completeness, we reinterpret our DFA for a language of
edge-address pairs as an NFA for a language of \emph{single} edge
addresses, by replacing each symbol pair $(a,b)$ on a transition with
just $a$. (Choosing $b$ instead should give the same results, because
the construction process is symmetric.) Then we search that NFA to
find any legal edge address which it does not accept. This is done by
constructing the DFA for the language $\mathcal{E}$ of \emph{all}
legal edge addresses, and performing a graph search over the cross
product of the two automata to find pairs of states reachable from
\textsc{start} via the same sequence of symbols. If this search
encounters any state pair and symbol from which the DFA for
$\mathcal{E}$ has a transition and the newly built NFA does not, then
we have not constructed a recogniser for all of $\mathcal{N}$, and
must try again.

A pseudocode sketch of the full algorithm for constructing a DFA for
$\mathcal{N}$ is shown in Algorithm~\ref{algadj}, with subroutines
given in more detail by Algorithms~\ref{algadjaddstring},
\ref{algadjmerge} and~\ref{algadjcomplete}.

\begin{figure}
\centering
\includegraphics*[alt={A finite state machine derived from the Penrose P2 tiling, in its Robinson-triangle form. There is a single start state, four accepting states labelled with the triangle types, and 18 intermediate states labelled with pairs of triangle types and edge indices. Each transition is labelled with a pair of edge address symbols.}, scale=0.9]{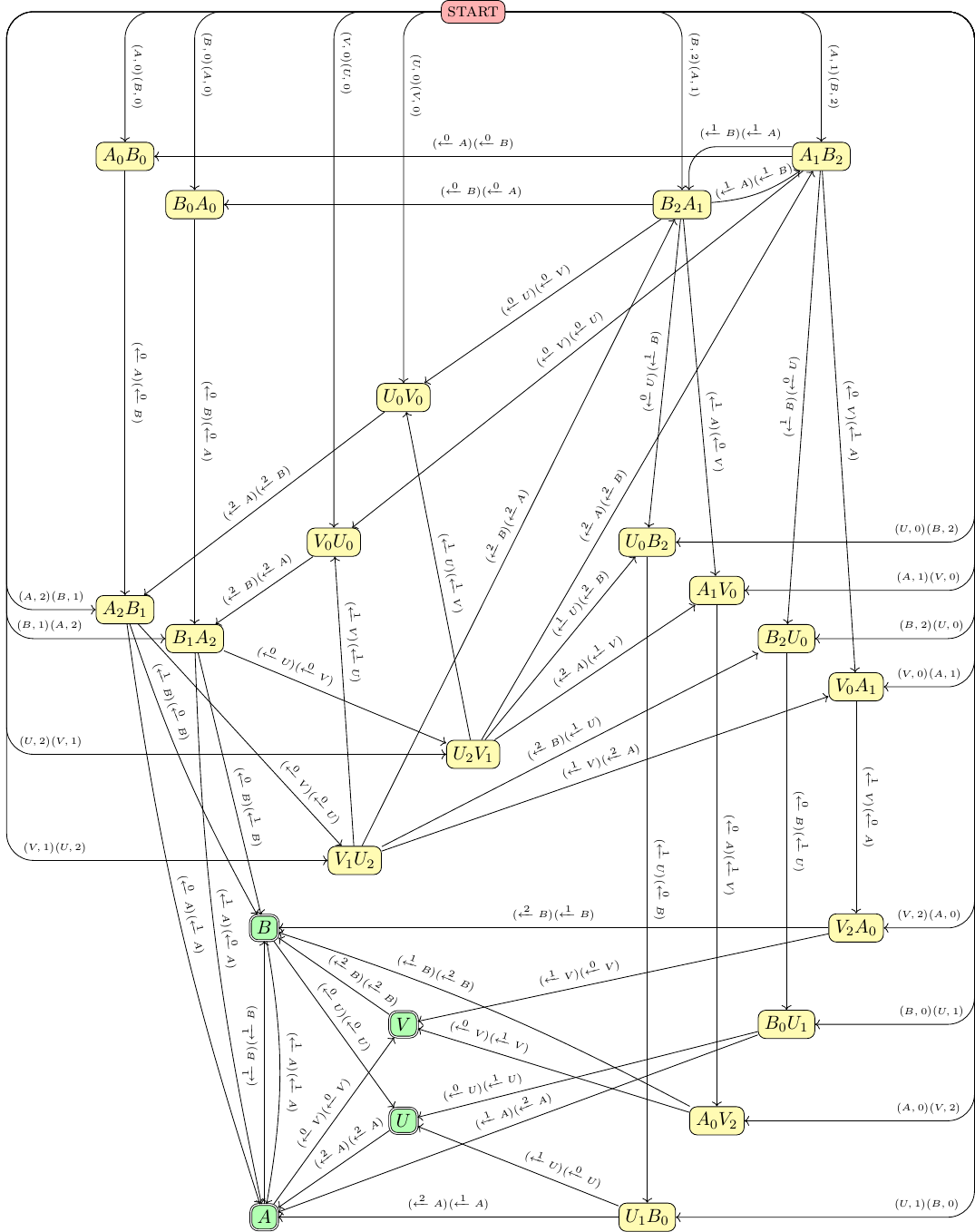}
\caption{Example recogniser for neighbour edge pairs, in the P2 Robinson
  triangle system from Figure~\ref{p2-triangles-adjmap}. \\
Non-accepting states are labelled with two tile types and the edges of those tiles that meet. \\
Accepting states (double outline) are labelled with a single tile type. \\
} \label{p2recogniser}
\end{figure}

A full example of a recogniser generated by this algorithm is shown in
Figure~\ref{p2recogniser}. In this case, since the substitution system
is simple enough that adjacent tiles always have adjacent supertiles,
it is possible to label each state of the recogniser with a specific
pair of tile edges. In more complicated recognisers this is not always
possible.

\subsubsection{Deterministic transducer}\label{transducer}

The DFA built in the previous section solves the problem of deciding
whether two tile edge addresses are neighbours of each other. However,
this is not the problem we typically want to solve! Usually we have
\emph{one} tile edge address, and want to calculate the address of its
neighbour.

By changing our point of view, we can reinterpret our existing DFA as
a finite-state \emph{transducer}, with each transition on a symbol
pair $(a,b)$ reinterpreted as a transition on the input symbol $a$,
producing $b$ as output.

This change of viewpoint renders the automaton nondeterministic.
Interpreted as a recogniser for neighbouring pairs, it could have two
state transitions $S\xrightarrow{a,b}T,S\xrightarrow{a,c}U$, with the
left-hand symbol $a$ in common, and going to different states
depending on whether the right-hand symbol is $b$ or $c$. When $b$ or
$c$ is provided as input, the choice of which transition to take is
forced, and the automaton can be deterministic. But interpreted as a
transducer accepting $a$ as input and emitting one of $b$ or $c$ as
output, this is less helpful -- the automaton now has multiple choices
of transition to take from that state on the input $a$.

The next algorithm attempts to convert this nondeterministic
transducer into a deterministic one. This transformation loses the
invariant that every transition outputs exactly one symbol:\ instead,
each transition on an input symbol outputs a \emph{string} of output
symbols, which may be empty, or contain more than one symbol. So the
output of the automaton will lag behind its input, in situations such
as the one just described, where the input symbol $a$ does not give
enough information to know whether the corresponding output symbol
should be $b$ or $c$. However, if the transducer has finitely many
states, the lag is bounded, so with a sufficiently long input string
every output symbol will be generated eventually.

\begin{algorithm}[t]
\caption{Try to make a deterministic transducer from a recogniser for $\mathcal{N}$}\label{determinise}
\begin{algorithmic}
\State $\textit{discovered}\leftarrow\bigl\{\{(\textsc{start},\epsilon)\}\bigr\}$;
\State $\textit{visited}\leftarrow\emptyset$;
\State $\textit{transitions}\leftarrow\emptyset$;
\While{$\mathit{discovered}\setminus\mathit{visited}$ contains an element \textit{src-dfa-state}}
  \State add \textit{src-dfa-state} to \textit{visited};
  \For{\textbf{all} symbols $a$}
  \State $\textit{dest-dfa-state}\leftarrow\emptyset$;
    \For{$(\textit{src-rec-state},\textit{output})$ in \textit{src-dfa-state}}
      \For{each $\textit{src-rec-state}\xrightarrow{(a,b)}\textit{dest-rec-state}$ in the recogniser for $\mathcal{N}$}
        \State add $(\textit{dest-rec-state},\textit{output}\circ b)$ to \textit{dest-dfa-state};
      \EndFor
    \EndFor
    \State let $\textit{prefix}=\hbox{longest common prefix of all output strings in }\textit{dest-dfa-state}$;
    \State remove \textit{prefix} from start of every output string in \textit{dest-dfa-state};
    \State add \textit{dest-dfa-state} to \textit{discovered};
    \State mark \textit{dest-dfa-state} as accepting if the NFA states in it are accepting;
    \State add $\textit{src-dfa-state}\xrightarrow{a}\textit{dest-dfa-state}$ to \textit{transitions} with output $\textit{prefix}$;
  \EndFor
\EndWhile
\State return $(\textit{discovered},\textit{transitions})$;
\end{algorithmic}
\end{algorithm}

The procedure for doing this is very similar to the procedure for
converting an ordinary (recogniser-style) NFA into a DFA, in which
each DFA state is identified with a subset of the states of the NFA.
To apply this procedure to transducers rather than recognisers,
permitting uncertain output to be delayed until all but one of the
possibilities has been eliminated, each DFA state must be identified
not with a plain set of NFA states, but with a set of pairs
$(\textit{NFA~state},\textit{string})$, where the \textit{string} is
`pending output':\ output generated by the path through the NFA to that
state, which has not yet been emitted as output from the DFA.

To calculate a transition from one DFA state to the next on a symbol
$a$, we enumerate all the transitions on $a$ from each of the
corresponding NFA states, and in each case, construct an output
$(\textit{state},\textit{string})$ pair with the NFA's output symbol
appended to the output string. Then, before converting this set of
pairs into a DFA state, we examine it to find out whether all the
state pairs agree on a prefix of their output strings. If so, then we
remove that shared prefix from all of the output strings, and record
it as the \emph{deterministic} output of the DFA transition we are
calculating.

For each state of the DFA, we expect its set of associated NFA states
to either all be accepting, or all non-accepting. In the former case
the DFA state is also marked as accepting.

(Proof:\ if the DFA construction made a state containing an accepting
and a non-accepting state, it would mean that there was an $m$-symbol
input string $A$ with two possible output strings $B,C$ such that the
pair $(A,B)$ reached an accepting state of the NFA while $(A,C)$
reached a non-accepting state. There must be some extension of $(A,C)$
which does reach an accepting state, say by appending extra symbols
$a,c$ to the two strings to make an $n$-symbol pair $(Aa,Cc)$. Then
$(Aa,Ba)$ must also be acceptable to the NFA, because after $(A,B)$
reaches an accepting state, further symbol pairs with both elements
equal are legal inputs and remain in an accepting state. So $Aa$ is an
$n$-symbol input with $Ba,Cc$ as distinct possible outputs, and since
both end in an accepting state, the $n$th-order supertile type is the
same in both cases. But this means that $Ba$ and $Cc$ identify
different lowest-level tiles within that $n$th-order supertile -- a
contradiction, since two different subtiles of that supertile cannot
both be adjacent to the same input tile along the same edge.)

\begin{figure}
\centering
\includegraphics*[alt={A second finite state machine derived from the Penrose P2 Robinson-triangle tiling. There is a single start state, four accepting states labelled with the triangle types, and 26 intermediate states. Unlike the previous state machine, each transition is labelled with one input edge address symbol and a string of output symbols with lengths ranging from 0 to 3. Each intermediate state is labelled with a set of states of the previous machine, with each of those states accompanied by another string of unemitted output symbols.}, scale=0.369]{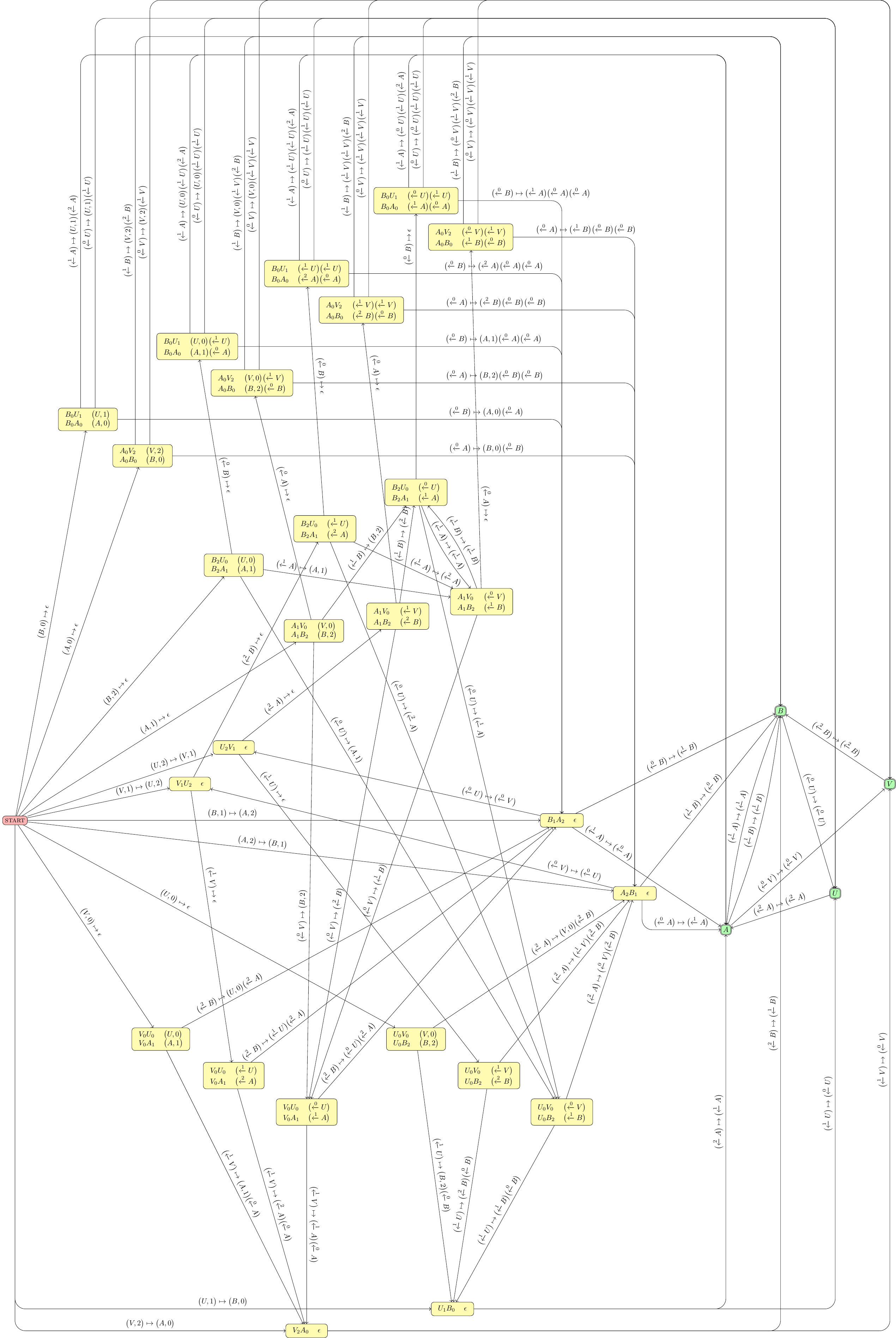}
\caption{Example transducer for the P2 Robinson triangle system
  (Figure~\ref{p2-triangles-adjmap}). \\ States other than
  \textsc{start} and accepting states are annotated with a set of
  (recogniser state, pending output) pairs, one on each line. The
  recogniser states are the same as the ones shown in
  Figure~\ref{p2recogniser}.} \label{p2transducer}
\end{figure}

Algorithm~\ref{determinise} shows this procedure in pseudocode. An example transducer built by the algorithm is shown in Figure~\ref{p2transducer}.

Including these pending output strings as part of the description of a
DFA state makes the set of potential states infinite, because strings
can be arbitrarily long. So, unlike the technique for determinising an
NFA \emph{recogniser} for a language, this algorithm can fail to
terminate. Indeed, for some important substitution systems, it does.
We discuss this in the next section.

\subsection{Ambiguous substitution systems}\label{ambiguous}

The transducer built by Algorithm~\ref{determinise} can emit no
symbols, or more than one symbol, on any given transition. But it is
the determinisation of a transducer that emitted exactly one output
symbol per input symbol, so it retains the invariant that the number
of output symbols it has emitted, plus the length of any pending
output string in its current state, must equal the number of input
symbols it has so far received.

Hence, if Algorithm~\ref{determinise} terminates, then there is an
integer $k$ such that the length of any pending output string in any
DFA state is $<k$. That is, the transducer is able to determine the
$n$th output symbol of any edge address after examining at most $n+k$
input symbols.

Conversely, if there is a $k$ such that $n+k$ input symbols are always
enough to rule out all but one possibility for the $n$th output
symbol, then Algorithm~\ref{determinise} must terminate, because the
set of possible $(\textit{NFA~state},\textit{string})$ pairs becomes
finite if the strings all have length $\leq k$.

We shall define a substitution system as \textbf{unambiguous} if
Algorithm~\ref{determinise} terminates successfully, and otherwise,
\textbf{ambiguous}.

\subsubsection{Detecting ambiguous systems}\label{detectfailure}

Conveniently, the failure of Algorithm~\ref{determinise} to terminate
can be detected reliably, enabling it to be converted into a
terminating algorithm which either determinises its input transducer,
or reports that it cannot.

Given a nondeterministic transducer which is not determinisable, there
is no upper bound on the length of `pending output' strings that can
appear in DFA states constructed by Algorithm~\ref{determinise}. In
particular, this means that for any $k$, there exists some input
string $I$ of length $m+k$ for which only the first $m$ output symbols
can be uniquely determined. In the nondeterministic transducer, this
means that there exist two legal paths $P_1,P_2$ from \textsc{start}
for that input string which differ in their $(m+1)$st output symbol,
say one path visits states $\textsc{start},s_1,s_2,\ldots,s_{m+k}$ and
the other
$\textsc{start},s_1,s_2,\ldots,s_m,s'_{m+1},s'_{m+2},\ldots,s'_{m+k}$.
By choosing $k$ sufficiently large we can find two integers $m<u<v<k$
after that symbol for which the state pairs $(s_u,s'_u)=(s_v,s'_v)$.
The segment of input and output between those positions can be pumped:
if the transducer is given the first $v$ symbols of $I$, followed by
infinitely many repetitions of the segment between $u$ and $v$, there
is a valid output sequence corresponding to each of these, by
similarly repeating each path's output between $u,v$. So the $(m+1)$st
output symbol never becomes unambiguous, because both of these paths
with different values for that symbol can be legally continued
indefinitely.

To detect this condition in practice, given a set of states
constructed by an unfinished execution of Algorithm~\ref{determinise},
it is enough to form a directed graph $G$ as follows, and check it for
cycles. Each vertex of $G$ corresponds to a set of ordered pairs of
the form $(\textit{NFA~state},\textit{symbol})$, and is made from a
state of the half-constructed DFA (which is a set of pairs
$(\textit{NFA~state},\textit{string})$) by discarding all but the
first symbol in each pending output string. (States in which the
pending output strings are empty need not be considered.) Edges of $G$
correspond to only those state transitions in the DFA which generate
no deterministic output.

Then a cycle in $G$ corresponds to a path between two DFA states which
generated no definite output and only appended more symbols to each
possible pending output, with the property that repeating the same
sequence of input symbols will have the same effect indefinitely.
Conversely, if the transducer is not determinisable, then a situation
of this kind must eventually arise from the construction given above,
because after each repetition of the section of $I$ between positions
$u,v$, the transducer must have returned to a state where \emph{two}
particular $(m+1)$th output symbols are still possible, and even if
other possible $(m+1)$th symbols are ruled out by going round the
cycle, eventually the set of possible symbols must stop decreasing
while still being greater than 1.

Hence, we enhance Algorithm~\ref{determinise} so that, interleaved
with the main algorithm, it periodically pauses the discovery of new
transducer states, constructs the graph $G$ described above, and
checks it for a cycle. If a cycle is ever found, the algorithm reports
failure; sooner or later, either that must happen, or the
determinisation must succeed.

\subsubsection{Neighbourhood-based refinement}\label{refinement}

The existence of ambiguous substitution systems is theoretically
interesting, but practically awkward, because a deterministic
transducer is a very convenient thing to have for a tiling, and it
would be nice if that convenience were always available.

The essential problem, in an ambiguous system, is that there is some
sequence of input symbols for which it is difficult to determine even
the first output symbol, i.e.\ the lowest-order tile adjacent to the
tile $t$ described by the first input symbol. (The start of the
difficult sequence need not coincide with the beginning of the entire
input edge address, although it often does in practice.) One way to
work around a problem like this is to provide extra information in
advance. In this case, that means \emph{refining} the tiling, by
cloning the tile type $t$ into multiple subtypes, and arranging that
each of the possible neighbour tiles can only legally appear next to
one of those subtypes. This requires more information to be provided
in the input tile address, but for some practical purposes this is no
difficulty (if you are generating a \emph{random} patch of tiling,
then you can invent more detailed tile subtypes at random as easily as
less detailed ones), and we shall also see in section~\ref{infsup}
that for some purposes it is unavoidable.

\begin{figure}
\centering
\includegraphics*[alt={A simple substitution system showing a `chair', or L-triomino, subdivided into four smaller chairs. The sub-chairs are labelled 0,1,2,3 and the edges of the original chair are numbered from 0 to 7.}]{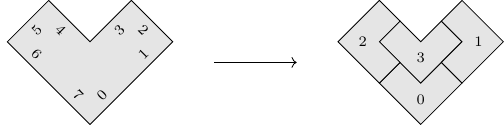}
\caption{Substitution system for the simple `chair' tiling} \label{chair}
\end{figure}

Of course, it is one thing to say `clone tile $t$', but another to do
it in practice. How many clones of $t$ do you need? When $t$ appears
in the deflation diagram of another tile type (or even of one of the
clones of $t$ itself), which kind of $t$ should it be replaced with?
How many other tile types need to be cloned so that different versions
of them can deflate to different kinds of $t$? Where does the cloning
end?

In this section we present an algorithm which automates this process,
by taking an input substitution system and outputting a refinement of
it:\ a finite number of distinct versions of each original tile type,
and a description of how to derive the new system's deflation diagrams
from the old system. This description only needs to provide new
information not present in the original system's deflation diagrams,
relating to which clone of each subtile appears where, and whether it
depends on which clone of the supertile is being deflated. All the
other details, like edge subdivisions and neighbour maps, are
unchanged.

Since the purpose of creating clones of a tile is to restrict what
neighbours can appear next to them, the essential idea is to make a
clone of each tile type corresponding to each of its possible
neighbourhoods.

\begin{figure}
\centering
\includegraphics*[alt={Seven diagrams showing a chair surrounded by other chairs, leaving no edge of the central chair exposed. The pattern of outer chairs is different in every case. The central chairs are numbered 1 to 7, each with a distinctive colour.}]{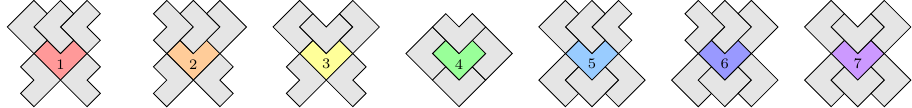}
\caption{All the neighbourhoods of a single chair} \label{chair7}
\end{figure}

A natural definition of the `neighbourhood' of a tile in a tiling is
to list, for each of the tile's edges, what type of tile is adjacent
to that edge, and which edge of the neighbour tile. We show an example
of classifying tiles by this type of neighbourhood, based on the
hierarchical substitution system of L-triominoes (or `chairs') from
\cite{Goodman-Strauss1999} with a single tile type. Figure~\ref{chair}
shows the substitution rule for the chair tiling. It also numbers all
the edges of the chair tile, and assigns identifiers for the four
subtiles in a deflation ($B,L,R,C$ for `bottom', `left', `right' and
`centre'). Figure~\ref{chair7} shows the seven possible neighbourhoods
of a single chair that can arise in a tiling generated by this system.

\begin{figure}
\centering
\includegraphics*[alt={The type-1 chair diagram from the previous figure, shown transformed by deflating the central chair and all the outer chairs into four subchairs. The leftmost subchair of the central chair is highlighted, along with all its immediate neighbours, showing that it matches the pattern of the type-2 chair from the previous figure.}]{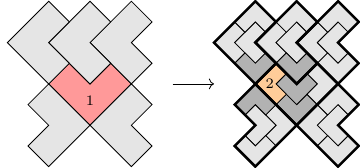}
\caption{Deriving the subtype of a tile in a deflation diagram} \label{refine-subtile}
\end{figure}

If we are to make a refined substitution system using these seven
possibilities as the subtypes of the chair tile, then we must
construct the deflation diagram for each subtype. These deflation
diagrams will be refinements of the diagram for the chair itself
(Figure~\ref{chair}); the only new decision to make, for each subtype
of chair being deflated, is which subtype of chair appears in each
position in its deflation diagram.

\begin{figure}
\centering
\includegraphics*[alt={Seven copies of the chair deflation diagram, one for each of the seven chair types from a previous figure. In each diagram, the subchairs are also marked with their types and colours, derived by the process shown in the previous figure.}]{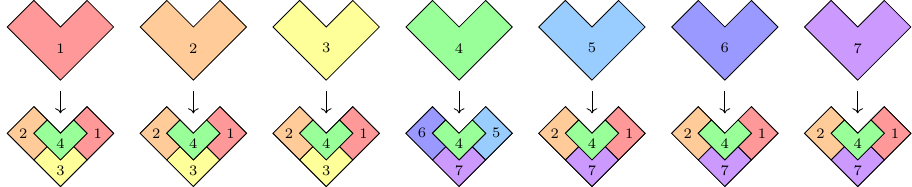}
\caption{Full deflation diagrams for the 7-chair refinement} \label{chair7-deflations}
\end{figure}

To determine this, we start from the neighbourhood diagram for the
subtype of the larger chair; deflate that entire diagram into smaller
chairs; and from there, determine the neighbourhood (and hence the
subtype) of each chair deflated from the original tile.
Figure~\ref{refine-subtile} shows an example of this, demonstrating
that when deflating a chair of subtype 1 (according to the numbering
in Figure~\ref{chair7}), the leftmost subtile has subtype 2. Repeating
this procedure for all the subtiles of all seven chair types, we
obtain a full set of deflation diagrams for the refined chair tiling.
Figure~\ref{chair7-deflations} shows these diagrams in full.

The chair substitution system is ambiguous. With just one chair type,
Algorithm~\ref{algadj} successfully constructs a recogniser for
neighbouring edge address pairs, but Algorithm~\ref{determinise} fails
to convert it into a deterministic transducer. However, the refined
system with seven chair types is unambiguous:\ it \emph{does} admit a
deterministic transducer.

To see why this helps, here is a concrete reason why the original
chair tiling is ambiguous. Consider a tile which is subtile 0 of
its supertile, which is subtile 0 of its supertile in turn, and
so on $n$ times. That is, its address is of the form
$$\textrm{chair}\underbrace{\xleftarrow{0}\textrm{chair}\xleftarrow{0}\textrm{chair}\xleftarrow{0}\cdots\xleftarrow{0}\textrm{chair}\xleftarrow{0}}_{\textstyle n\hbox{ steps}}\textrm{chair}$$
so that the lowest-order tile occupies the very bottom corner of the
$n$th-order supertile, in the orientation shown in figure~\ref{chair}.
What borders that tile on one of its lower edges -- say, edge 0 from
the same figure? It depends on what borders the $n$th-order supertile,
and the various possibilities do not all lead to the same orientation
of the \emph{smallest} chair:\ this $n$th-order chair might be either
back-to-back with the next one or nestled inside it, and the two
touching lowest-order chairs will have the same relationship as the
high-order ones in each case. So a transducer for this substitution
system cannot generate even the first output symbol without looking at
the $(n+1)$th input symbol. But $n$ was arbitrary. So the transducer
cannot have finitely many states.

The refined chair tiling resolves this problem by requiring the input
tile address to be specified in more detail. The deflation diagrams in
Figure~\ref{chair7-deflations} show that subtile 0 of any chair has
either type 3 or type 7:\ for chair types 1,2,3, subtile 0 has type 3,
and for types 4,5,6,7, it has type 7. In particular, each of types 3
and 7 has a subtile 0 of the same type as itself. So the input address
must now be specified in one of the following forms:
$$\begin{array}{ll}
\textrm{chair}_3\xleftarrow{0}\textrm{chair}_3\xleftarrow{0}\textrm{chair}_3\xleftarrow{0}\cdots\xleftarrow{0}\textrm{chair}_3\xleftarrow{0}\textrm{chair}_t &\quad t\in\{1,2,3\} \\
\textrm{chair}_7\xleftarrow{0}\textrm{chair}_7\xleftarrow{0}\textrm{chair}_7\xleftarrow{0}\cdots\xleftarrow{0}\textrm{chair}_7\xleftarrow{0}\textrm{chair}_t &\quad t\in\{4,5,6,7\} \\
\end{array}$$
This allows the transducer to know, when it sees the \emph{first}
symbol of the address, whether the $n$th symbol is going to be in the
set $\{1,2,3\}$ or the set $\{4,5,6,7\}$. More immediately, because
the lowest-order tile has either type 3 or 7, the neighbourhood
diagrams in Figure~\ref{chair7} show what orientation of chair tile
borders it along the lower edges. So the transducer has information it did not have in the original substitution system.

(This does not derive from any magical foresight involved on the part
of the transducer! Rather, changing the address notation in this way
has required the \emph{user providing the input tile address} to
compensate for the transducer's inability to look arbitrarily far into
the future, by giving the necessary information early in the address
rather than leaving it until later.)

This example with the chair tiling illustrates the basic procedure for
refining a substitution system. However, it isn't as simple as that in
every case, because if the substitution system has zero-thickness
spurs, then knowing the type of every formal neighbour of a tile $t$
may not be enough information to know the type of every neighbour of
$t$'s subtiles in a particular deflation:\ even if the neighbours of
$t$ fully surround $t$ leaving no exposed edge, the same might not be
true after deflating all the tiles.

To account for this, the refinement algorithm must be prepared to use
a more general definition of `neighbourhood'. The strategy is to make
a list of \emph{queries} about tiles in the tiling:
\begin{itemize}
\setlength\itemsep{0em}\setlength\parskip{0em}
\item the initial query $Q_0$ is simply `What type is the central tile?'
\item immediate followups are of the form `If the answer to query
  $Q_0$ was a tile of type $t$, what type is that tile's neighbour
  along edge $e$, and what edge of that neighbour meets edge $e$?'
\item further followups name the full answer to the previous query,
  including a tile type and an edge: `If the answer to query $Q_n$ was
  edge $e$ of a tile of type $u$, what type is that tile's neighbour
  along edge $f$, and what edge of that neighbour meets edge $f$?'
\end{itemize}
A query need only be answered if its precondition is satisfied. If a
query $Q_n$'s predecessor query $Q_m$ had an answer other than the one
assumed by $Q_n$, then the answer to $Q_n$ is $\bot$, or `not
applicable'. The same is true if the predecessor query's answer was
itself $\bot$.

\begin{algorithm}[t]
\caption{Answer the list of refinement queries for a tile}\label{refine-query}
\begin{algorithmic}
\Procedure{AnswerQueries}{information $T_0$ identifying the starting tile}
\State $A_0\leftarrow\textit{tile-type}(T_0)$;\Comment{$Q_0$ always asks for the starting tile type}
\For{queries $Q_n$ for $n=1,2,\ldots$}
  \State let $(m,u,e)$ = $Q_n$;\Comment{`if $Q_m$ was a tile of type $u$, try its edge $e$'}
  \If{$A_m$ has tile type $u$}
    \State $(T_n,e')\leftarrow\hbox{neighbour of tile }T_m\hbox{ across edge }e$;
    \State $A_n\leftarrow(\textit{tile-type}(T_n), e')$;
  \Else\Comment{including the case $A_m=\bot$}
    \State $A_n\leftarrow\bot$;
  \EndIf
\EndFor
\State return the list of $(A_n)$;
\EndProcedure
\end{algorithmic}
\end{algorithm}

\begin{algorithm}
\caption{Compute the neighbour of a subtile in the deflation of a
  supertile subtype}\label{refine-subtypes}
\begin{algorithmic}
\Procedure{SupertileNeighbour}{$\textit{answers},\textit{query},e$}
\State let $\textit{answer}=\hbox{answer to \textit{query} in \textit{answers}}$;
\If{\textit{answer} named edge $e$ of a tile type $t$}
  \State let $u,f$ be such that \textit{query} asked about edge $f$ of tile $u$;
  \State let $\textit{query} = \hbox{predecessor of }\textit{query}$;
  \State return $(\textit{query}, u, f)$;
\Else
  \State let $\textit{new-query} = (\textit{query}, e)$;
  \If{\textit{new-query} is not in the query list}
    \State return failure, and recommend adding \textit{new-query} to the list;
  \Else
    \State let $(u, f) = \textit{new-answer}=\hbox{answer to \textit{new-query} in \textit{answers}}$;
    \State return $(\textit{new-query}, u, f)$;
  \EndIf
\EndIf
\EndProcedure
\vskip0.5em
\Procedure{SubtileNeighbour}{$i_1,\textit{answers},\textit{query},t_1$}
\State $\textit{lookup-result} \leftarrow \hbox{look up }\Int{i_1}{e}\hbox{ in adjacency map for }t_1$;
\While {\textit{lookup-result} is of the form $\Ext{u}{v}$}
    \State $\textit{new-query},t'_1,u'\leftarrow\textsc{SupertileNeighbour}(\textit{answers},\textit{query},t_1,u)$;
    \State $\textit{lookup-result} \leftarrow \hbox{look up }\Ext{u'}{-v}\hbox{ in adjacency map for }t'_1$;
    \State replace $(t_1,\textit{query})\leftarrow(t'_1,\textit{new-query})$;
\EndWhile
\State now expect that \textit{lookup-result} is of the form $\Int{i'_1}{e'}$;
\State return $(i'_1,\textit{query},t_1)$;
\EndProcedure
\end{algorithmic}
\end{algorithm}

The simple definition of `neighbourhood' for the chair system is
easily converted into this form:
\begin{itemize}
\setlength\itemsep{0em}\setlength\parskip{0em}
\item[$Q_0$.] What type is the central tile? (in this case it will always be `chair')
\item[$Q_1$.] If the answer to $Q_0$ was `chair', what borders that chair on edge 0?
\item[$Q_2$.] If the answer to $Q_0$ was `chair', what borders that chair on edge 1?
\item[\;] $\vdots$
\item[$Q_8$.] If the answer to $Q_0$ was `chair', what borders that chair on edge 7?
\end{itemize}

To refine a general substitution system, we initialise the list of
queries in this fashion:\ for each tile type, we make a set of queries
following up $Q_0$ and asking about every edge of that tile. Then we
find tile subtypes by searching the substitution system for finitely
long tile addresses which provide enough information to answer all
queries in the list. Then we attempt to determine the deflation
diagram for each of those subtypes, by answering the same list of
queries about each of its subtiles, using no information except the
answers to the same queries for the supertile. If this attempt fails
due to not being able to determine what is on the other side of a
spur, the failure comes with a description of an additional query that
would have allowed the algorithm to make more progress. We collect as
many of these extra queries as possible, and then restart the
algorithm from scratch with the extended query list.

The algorithm terminates when both the set of queries and the set of
tile subtypes are closed:\ the set of queries is such that knowing the
answers for a supertile allows the answers for all its subtiles to be
determined, and every tile subtype appearing as a subtile in a
deflation diagrams is also a subtype for which a deflation diagram is
known in turn.

The key part of this algorithm is finding the answers to the list of
subtype-determining queries. Algorithm~\ref{refine-query} gives the
general procedure for doing this, supposing that it has access to some
type of information identifying a particular tile, allowing the tile's
type to be determined, and a method for identifying the neighbour of
the tile.

To determine a starting list of tile subtypes,
Algorithm~\ref{refine-query} is used in a form where the `information
identifying a tile' is a finitely long tile address giving the tile's
position within a deflation of some particular supertile, and a tile's
neighbour is determined by passing that address to
Algorithm~\ref{algrec} to calculate the similar address of the
neighbour. If this procedure completes without Algorithm~\ref{algrec}
recursing too deep and running off the end of the input address, then
the tile's subtype is successfully determined.

\begin{algorithm}[t]
\caption{Compute a refinement of a whole tiling}\label{algrefine}
\begin{algorithmic}
\State $\textit{queries}\leftarrow\{Q_0\}\cup\{(0, t, e)\}$;
\For{$n=1,2,3,\ldots$}\Comment{length of tile address to explore}
  \State $\textit{subtypes}\leftarrow\emptyset$;
  \For{\textbf{all} finite tile addresses of length $n$}
    \State let $\textit{answers}=\textsc{AnswerQueries}(address)$;\Comment{via Algorithm~\ref{algrec}}
    \If{all queries answered successfully}
      \State $\textit{subtypes}\leftarrow\textit{subtypes}\cup\{(t,\textit{answers})\}$;
    \EndIf
  \EndFor
  \State $\textit{incomplete}\leftarrow(\textit{subtypes}\neq\emptyset)$;\Comment{need at least one subtype}
  \For{$(t,\textit{super-answers})$ in \textit{subtypes}}\Comment{make deflation diagrams}
    \For{each deflation of $t$}
      \For{each subtile index $i$ of that subtype}
        \State let $\textit{sub-answers}=\textsc{AnswerQueries}(\textit{super-answers},Q_0)$;\Comment{via Algorithm~\ref{refine-subtypes}}
        \If{Algorithm~\ref{refine-subtypes} failed and suggested an extra query}
          \State add the new query to \textit{queries};
          \State $\textit{incomplete}\leftarrow\textbf{true}$;
        \EndIf
        \If{$\textit{sub-answers}\not\in\textit{subtypes}$}
          \State $\textit{incomplete}\leftarrow\textbf{true}$;
        \EndIf
        \State enter \textit{sub-answers} as type of subtile $i$ in deflation of \textit{subtype};
      \EndFor
    \EndFor
  \EndFor
  \If{not \textit{incomplete}}
    \State return \textit{subtypes} and deflation diagrams;
  \EndIf
\EndFor
\end{algorithmic}
\end{algorithm}

However, when determining the subtype of a tile in a deflation
diagram, a full address for the subtype is not available, because it
is necessary to ensure that the subtypes of subtiles can be derived
\emph{only} from the answers to the queries for the supertile. So for
this purpose we reuse Algorithm~\ref{refine-query} using a different
representation of a tile identity, and a modified algorithm for
determining its neighbour.

At the layer of the supertile $t$, for which we know the answers to
all the subtype queries, our representation of a tile identity for
these purposes simply nominates a particular query $Q_i$, and the tile
referred to is whichever tile $Q_i$ is asking about. Given a tile
identity specified in that form, finding the neighbour of the tile
along some edge $e$ is normally done by extending the query $Q_i$ into
a longer one, incorporating the answer given in the tile identity: `If
the answer to $Q_i$ was [this edge of that tile], what's on the other
side of its edge $e$?'. The only exception is if $e$ is the same edge
that $Q_i$ asked about, in which case we instead \emph{shorten} the
previous query by a step, returning to whatever query it was itself a
followup to.

However, when determining tile subtypes in a deflation diagram, we
need to find neighbours at the \emph{subtile} layer, not the supertile
layer. So the full algorithm for finding neighbours in this context
uses the same structure as the recursive Algorithm~\ref{algrec}, using
the adjacency map for each supertile, but where Algorithm~\ref{algrec}
recursed by calling itself, this method instead calls the supertile
neighbour procedure described in the previous paragraph.
Algorithm~\ref{refine-subtypes} shows this in full.

Algorithm~\ref{algrefine} shows the top-level structure of the overall
refinement procedure, using all the previous algorithms in this
section as subroutines.

One aspect omitted from the pseudocode for brevity is handling
substitution systems involving multiple layers. Tile types at all
layers potentially need refining into multiple subtypes. So when
Algorithm~\ref{algrefine} iterates over all tile addresses of length
$n$, it must include addresses of tiles \emph{at all layers}, not just
length-$n$ addresses of tiles in the base layer. If a tile type
appears in multiple layers, it should be considered to be a separate
tile of the same shape for each layer it appears in, so that the list
of queries for each layer can vary independently.

It is logical that this refinement procedure should be able to
transform an ambiguous substitution system into an unambiguous one,
because it incorporates information into each tile type about that
tile's immediate neighbourhood -- i.e.\ about exactly what a
transducer needs to know to produce its output.

\begin{figure}[t]
\centering
\includegraphics*[alt={The rhombus and square tiles from the Ammann-Beenker tiling, shown with deflation diagrams. The square is shown bisected into two triangles, which allows those triangles to appear individually in the deflation diagrams so that each tile's deflation covers the same region as the original tile. Edges of the tiles are marked with single and double arrows to indicate type and direction; the tile vertices are marked with shapes made from isosceles right triangles which join together into a thick arrow shape at a complete vertex. The rhombus is shown in two mirror-image forms, coloured red and green, differing only in the orientation of the vertex marks at its pointy ends.}]{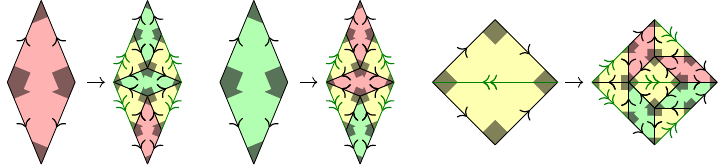}
\caption{Deflation diagrams for the Ammann-Beenker tiling} \label{ammannbeenker}
\end{figure}

\begin{figure}[t]
\centering
\includegraphics*[alt={Eight Ammann-Beenker rhombi as shown in the previous figure, with a sharp corner of each one all meeting at a single vertex. The pattern of red and green rhombi is such that the vertex marks at the centre form an arrow shape. At the next 8 vertices where two rhombi meet, the vertex marks are such that an arrow shape could still be formed by whatever fills in the remaining gap.}]{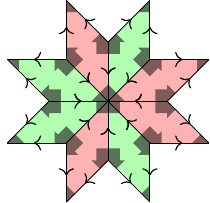}
\caption{The only way 8 Ammann-Beenker rhombi can meet at a vertex} \label{ammannbeenker8}
\end{figure}

However, a single run of the refinement algorithm is not always
sufficient. The Ammann-Beenker tiling (\cite{10.5555/19304},
Figures~10.4.14 and~10.4.15) is an example of a substitution system
requiring two passes of this refinement algorithm before a transducer
can be successfully generated. Figure~\ref{ammannbeenker} shows the
substitution system in a form suitable for use with these algorithms.
The Ammann-Beenker tiling includes a constraint about the order of
tiles that can appear around each vertex, shown in the figure as
shaded markings which are required to form into a thick arrow shape
around the vertex. The markings on the rhombus are asymmetric, and in
order to meet the constraint, both mirror-image forms of it must occur
in a specific pattern around any vertex where 8 rhombus tiles meet,
shown in Figure~\ref{ammannbeenker8}.

An unambiguous version of the system must therefore distinguish all 8
of the rhombi around such a vertex. Knowing each tile's immediate
neighbours is not enough to distinguish all 8 rhombi, because (for
example) there are two instances of a green rhombus with both
neighbours red. But a second run of the refinement algorithm includes
information about the \emph{once-refined} type of each of the
rhombus's neighbours, i.e.\ information about the original type of the
rhombus's \emph{second-order} neighbours, which includes 5 of the
rhombuses around an 8-rhombus vertex -- and that \emph{is} enough
information to distinguish all 8, because between 5 rhombuses are 4 of
the vertex's 8 incident edges, and one of those must be the
green/green or red/red edge. Full details of the necessary refinement
are given in section~\ref{ammann-beenker}.

The more times this algorithm is run on a system, the more each tile
type incorporates information from a larger and larger area of the
tiling around the tile. After some finite number of refinements, each
tile will describe a neighbourhood of itself sufficiently large that,
given any full tile address containing those refined tile types, the
regions of the plane forced by the neighbourhoods of larger and larger
supertiles expand in all directions and cover the whole plane. So it
should always be possible to make a substitution system unambiguous
after applying the refinement algorithm enough times.

\subsubsection{Remerging tile types}\label{reduce}

\begin{figure}[b]
\centering
\includegraphics*[alt={Three deflation diagrams for different types of chair, similar to the previous figure with seven of them. This time the chair types are labelled `x', `y' and `z'.}]{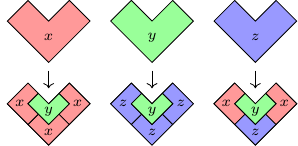}
\caption{Deflation diagrams for the 3-chair refinement} \label{chair3-deflations}
\end{figure}

In the previous section, as an example, we derived a refinement of the
chair substitution system containing 7 different tile types. That
refinement is unambiguous, where the original system is ambiguous.

In the deflation diagrams shown in Figure~\ref{chair7-deflations},
chair types 1, 2 and 3 all have identical deflation diagrams. So do
chair types 5, 6 and 7. We can reduce the number of tile types by
merging each of those sets into a single combined type, obtaining
(say) types $x=\{1,2,3\},y=\{4\},z=\{5,6,7\}$, to produce the
resulting simpler system shown in Figure~\ref{chair3-deflations}.

This system remains unambiguous:\ merging the synonymous tile types
does not cause Algorithm~\ref{determinise} to fail. That is to be
expected, because a `coarse' tile address given in terms of the three
types $x,y,z$ can be readily converted back to a `fine' one using the
seven original types $1,\ldots,7$ with only a finite amount of
lookahead:\ to recover the fine subtype of a chair at level $n$, it is
necessary only to look at its coarse supertile type at level $n+1$.
All the fine types merged into that coarse type had exactly the same
deflation diagram, so the coarse supertile type plus the subtile index
is enough to recover the fine type of the subtile.

Even if this reduction were performed more than once (if the first
round of merging tile types caused some deflation diagrams to become
identical which were previously not), this argument would still
hold:\ each reduction can be undone by a finite state machine using a
single symbol of lookahead, therefore finitely many reductions can be
undone using finite lookahead. So there is no reason why reducing a
system in this way should prevent a finite-state transducer from being
able to handle it.

\subsection{Putting it all together}

The algorithms presented in this section combine into a
general-purpose method for computing a deterministic finite-state
transducer for tiling substitution systems.

Given a substitution system, one first uses Algorithm~\ref{algadj}
(section~\ref{adjmatcher}) to make a DFA that matches the neighbour
language $\mathcal{N}$ of pairs of neighbouring tile edge addresses.
Then Algorithm~\ref{determinise} (section~\ref{transducer}) converts
that into a transducer capable of taking one tile edge address as
input, and producing the neighbouring one as output. If
Algorithm~\ref{determinise} cannot succeed, this can be detected via
the technique in Section~\ref{detectfailure}. In that situation, the
substitution system is ambiguous, but all hope is not lost:\ run
Algorithm~\ref{algrefine} (section~\ref{refinement}) to refine it into
a substitution system more likely to be unambiguous, and try the whole
procedure again.

All of these algorithms are practical in their running time, and
suitable for real-world use in software. In the author's own
implementations, compiled code (in Rust) can process a typical
substitution system in well under a second, and even unoptimised
interpreted code (in Sage) in under a minute. And a transducer only
needs to be constructed once for a given tiling, and then can be
reused for as much work as is needed.

Some applications of these finite-state automata are given in
section~\ref{applications}.

\section{Infinite supertile boundaries}\label{infsup}

So far, we have discussed the question of finding the address of a
tile's neighbour on the assumption that the two tiles share an
$n$th-order supertile, for some $n$. The recursive
Algorithm~\ref{algrec} must ascend the levels of supertiles until it
reaches such an $n$ before it can begin generating any output at all;
the transducer constructed by Algorithm~\ref{determinise} signals that
a shared supertile has been reached by entering an accepting state,
and if further input symbols are provided after that, emitting them
unchanged as output (checking them for validity as a side effect).

But $n$ can be arbitrarily large. For any $n$, it is easy to find
examples of adjacent tiles which have no supertile in common before
the $n$th layer, simply by making a single $n$th-order tile $T$,
deflating it $n$ times, and choosing two 0th-order tiles that came
from different immediate subtiles of $T$. It follows that arbitrarily
long paths exist through the transducer state machine which never
enter an accepting state.

The transducer has finitely many states. So if \emph{arbitrarily} long
paths exist in it which never reach an accepting state, then there
must also exist \emph{infinitely} long paths with the same property.
That is, there exists an infinite sequence of input symbols for which
the transducer never enters an accepting state -- but it still
produces an infinite sequence of output, since no transducer state has
more than a finite amount of pending output, and the invariant is
always maintained that total length of output (including pending
output in the current state) equals the total length of input consumed
so far.

An input sequence like this describes a tile $t$, and its output a
neighbouring tile $t'$, such that the two tiles have no supertile in
common at all. The union of the patches of tiling deflated from $t$'s
supertiles does not cover the whole plane, and the same is true for
$t'$. \cite{Goodman-Strauss1999} describes these as `infinitely large
supertiles'; one might also describe them as supertiles of infinite
\emph{order}. They are both, so we shall simply say `infinite
supertile'.

The recursive Algorithm~\ref{algrec} cannot handle tile addresses of
this type. It cannot begin generating any output until it has found an
$n$th-order supertile in common between its input and output tiles,
and in this case, no such supertile exists. So it will keep recursing
until it runs out of stack, and never generate any useful output. If
one imagines each $n$th-order supertile's patch of tiling to be
surrounded by a wall of height $n$, then Algorithm~\ref{algrec} must
climb over such a wall step by step in order to get over it to the
next tile -- but an \emph{infinite} supertile boundary, being the
union of a wall of every integer height, prevents it from passing at
all.

However, our deterministic transducer has no trouble crossing those
infinitely high walls, because it generates its output address from
the low order upwards, before even \emph{knowing} that the wall is
infinitely high. Indeed, it never receives enough information to be
certain of that. For all the transducer knows, at any point, the next
few symbols might lead to an accepting state, and there might be a
shared finite-order supertile between its input and output after all.

Proper infinite supertiles (that is, not covering the entire plane)
often occur in the most interesting instances of a given aperiodic
tiling, such as instances with full reflective or rotational symmetry,
or `singular' tilings with near-symmetry. So an algorithm for
calculating neighbour addresses which can cross their boundaries is
desirable, because it allows those particularly interesting tiling
instances to be generated and investigated conveniently.

An example of an infinite supertile boundary occurs in the P2 Robinson
triangle substitution system shown in Figure~\ref{robinson1}. From
that figure we can see that the upper part of edge 1 of the $A$
triangle, under deflation, becomes edge 2 of a $B$ triangle (which we
have indexed as subtile number 1 in figure~\ref{robinson2}).
Conversely, the upper part of edge 2 of a $B$ becomes edge 1 of an $A$
(again subtile number 1). So we can construct the following address,
describing an $A$ tile whose sequence of supertile types alternates
between $B$ and $A$ forever, and in particular describing edge 1 of
the lowest-order $A$:
$$\mathcal{A}=\Edge{A}{1}, \Sup{1}{B}, \Sup{1}{A}, \Sup{1}{B},
\Sup{1}{A}, \Sup{1}{B}, \Sup{1}{A}, \ldots$$
This is a concrete example of an address that will cause
Algorithm~\ref{algrec} to attempt to recurse to infinite depth, and
fail to terminate. Every time it is called to traverse edge 1 of a
type-$A$ tile, it will recurse to attempt to traverse edge 2 of its
type-$B$ supertile, and \emph{vice versa}.

The infinite supertile containing this address is a half-plane, whose
boundary is a straight line. From first principles we can construct an
address for a supertile that can occupy the other half-plane and fit
next to this one, by symmetry. The Robinson-triangle system is
invariant under the combined transformation of replacing each tile
with its mirror image, and interchanging the tile types
$A\leftrightarrow B,U\leftrightarrow V$. (This also interchanges edge
indices 1 and 2, since the edges of each triangle are indexed
anticlockwise.) The image of address $\mathcal{A}$ under this symmetry
gives rise to a mirror-image half-plane, with the same pattern of long
and short tile edges along its boundary, fitting exactly to the first
one:
$$\mathcal{B}=\Edge{B}{2}, \Sup{1}{A}, \Sup{1}{B}, \Sup{1}{A},
\Sup{1}{B}, \Sup{1}{A}, \Sup{1}{B}, \ldots$$

\begin{figure}
\centering
\includegraphics*[alt={A patch of Penrose P2 tiling, namely the famous `infinite cartwheel pattern' with 10-way near-symmetry. The symmetric regions of the diagram are in light colours, and the parts which vary under at least one rotation are shaded dark. A thick line down the centre shows the single axis of true reflection symmetry.}]{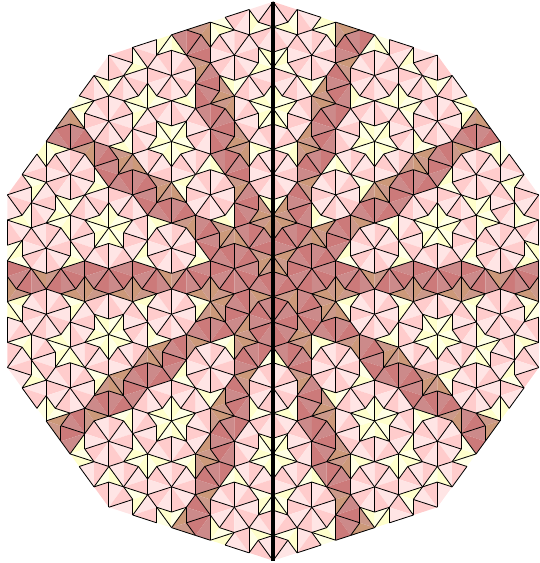}
\caption{The Penrose P2 infinite cartwheel pattern} \label{cartwheel}
\end{figure}

This is perhaps the simplest example of an infinite-order supertile
boundary. It is also a good example of the principle that tile
addresses of this kind tend to coincide with tilings of particular
interest, because these tile addresses represent the two tiles at the
centre of the famous `infinite cartwheel pattern', an instance of P2
which not only has reflection symmetry but \emph{almost} has 10-way
rotation symmetry too\footnote{This phenomenon of near-symmetry
derives from the `singular pentagrids' described in
\cite{DEBRUIJN198153}:\ the tiling is ultimately generated from a
pentagrid which \emph{actually} has the symmetry group of a decagon,
but at the points where more than two pentagrid edges meet, gaps are
created which can be filled with Penrose tiles in multiple ways, and
no way to fill the gaps can preserve all the symmetry. We shall use
the word `singular' for other nearly-symmetric instances of
substitution tilings, by analogy.}, in that \emph{most} of the plane
is invariant under a rotation of $36^\circ$, with only the centre of
the image and a few thin spokes varying. Figure~\ref{cartwheel} shows
this pattern, with the non-rotationally-symmetric parts shaded, and
the line of true reflection symmetry running vertically down the
centre.

The example transducer in Figure~\ref{p2transducer} demonstrates that
these two addresses $\mathcal{A},\mathcal{B}$ are indeed mapped to one
another, without any difficulty arising from the infinite supertile
boundary: the symbols of address $\mathcal{B}$ shown above are
generated as output from the transducer, if one follows transitions
from \textsc{start} on the symbols of $\mathcal{A}$, and \emph{vice
versa}.

\subsection{Eventually periodic tile addresses}\label{periodic}

The address $\mathcal{A}$ shown above has the useful property of being
\emph{eventually periodic}. Apart from the initial special symbol
containing an edge type instead of a supertile index, the rest of the
address repeats the two symbols $\Sup{1}{B},\Sup{1}{A}$ for ever.

Eventually periodic tile addresses are particularly well suited to a
finite-state transducer, because given a description of an eventually
periodic address in the form (initial segment, repeating segment), a
finitely long computation can deliver the neighbouring tile's address
in the same form, simply by waiting until a second visit is made to
the same combination of a state of the transducer and a position in
the repeating part of the input. Then all the output generated since
the first visit to that combination is guaranteed to repeat forever,
since in every subsequent repetition, the transducer will receive the
same sequence of input symbols starting from the same
state.\footnote{As a by-product of this algorithm, it is also possible
to detect the boundary between two infinite supertiles. An edge
between two tiles is part of an infinite supertile boundary precisely
when a transducer generates infinite output without ever entering an
accepting state. In the case of eventually periodic tile addresses,
the transducer's sequence of state transitions is also eventually
periodic, so this can be determined in finite time.}

It follows that the class of eventually periodic tile addresses is
closed under the neighbour relation:\ every tiling of the plane based
on a substitution system is either composed entirely of eventually
periodic addresses, or has none of them at all.

Since the class of eventually periodic addresses appears to include
many (if not all) of the tilings of special interest, like symmetric
and singular tilings, it's useful to have a way to talk about them
concisely. A useful notation is to describe a whole \emph{infinite
supertile} at a time, by mentioning only the eventually repeating
segment of the tile address, normalised to be (a) of minimum length
(never talk about repetitions of \textit{abab} when you can just say
\textit{ab}), and (b) starting at a position in the address which is a
multiple of its own length.

For example, the address $\mathcal{A}$ in the previous section is part
of the infinite supertile $\left[\Sup{1}{A},\Sup{1}{B}\right]^*$,
whereas its neighbour $\mathcal{B}$ is part of the infinite supertile
$\left[\Sup{1}{B},\Sup{1}{A}\right]^*$. The repeating sequence in each
of these infinite supertiles is a cyclic rotation of the other,
implying that the two supertiles can be obtained from one another by
deflation. In this case they also appear in the same tiling of the
whole plane with a boundary between them, but in other cases, they are
entirely separate tilings.

\subsection{Uniqueness}

When a deterministic transducer crosses an infinite supertile boundary
and returns a unique address for the tile on the far side, it does not
necessarily mean that it has found the only \emph{layout of tiles}
that physically fits to the far side of the boundary (even counting
matching rules implicit in the substitution system, such as types and
directions of edges).

A transducer constructed by this method only permits tile adjacencies
that would arise in finitely large patches generated by the input
substitution system. Any other way to fit tiles together is not
considered. So a transducer will only be able to say that there is a
unique layout of tiles, with their types, such that every finitely
large patch crossing the boundary is consistent with the substitution
system.

If an infinite supertile does extend to the rest of the plane in
multiple ways that physically fit, then it will have multiple
representations in the tile addressing system, and which one you give
as input to the transducer controls which of the possible extensions
you receive as output. Sections~\ref{spectre} and~\ref{hat} exhibit
example cases in which this occurs.

\section{Applications}\label{applications}

Now that we have a suite of algorithms for dealing with substitution
tilings using finite state automata, what are they useful for?

\subsection{Tiling generation}\label{tiling-generation}

Most obviously, a finite-state transducer for a substitution system
provides a highly practical method of generating patches of the tiling
described by the system, for whatever purpose you might need one --
mathematical research, puzzle games, or simply pretty pictures.

Perhaps the most obvious method of generating a patch of a Penrose
tiling, or hats, or Spectres, is to start from a single order-$n$ tile
or metatile, for some $n$ large enough for your purposes; apply the
substitution system $n$ times to deflate that single $n$th-order tile
into a large patch of 0th-order tiles; then pick out from that patch a
smaller region of whatever size and shape you really wanted, and throw
the rest of the output tiles away.

In place of that approach, we present a family of approaches based on
knowing an address for each tile plotted. In general, one chooses a
starting tile, equipped with both a position in the plane and an
address within the substitution system. Then one calculates the
address of a neighbouring tile via the combinatorial algorithms
presented in Section~\ref{algorithms}, obtaining in particular its
lowest-order tile type, and finding out which pair of edges of the two
tiles coincide, which is enough to place that tile in the plane too --
so now that second tile also has known plane coordinates and a known
address, enabling the same procedure to be repeated. One can search
outwards from the starting tile in any pattern desired, discovering
new tiles on demand, and stop as soon as all the tiles overlapping the
desired region have been found.

The simplest method of doing this is to run the neighbour algorithm
\emph{in full} to compute a complete address for each new tile. One
could give the starting tile's address to a deterministic transducer,
and run it until the transducer reaches an accepting state.
(Equivalently, and requiring less complicated programming, one could
not bother constructing any transducers at all, and do this same job
directly by using the recursive Algorithm~\ref{algrec}.)

This approach only works if there is no possibility -- or, failing
that, at least zero probability -- of encountering an infinite
supertile boundary, since in that situation disaster occurs:\ the
transducer will \emph{never} reach an accepting state, or
Algorithm~\ref{algrec} will recurse forever and overflow its stack. If
the input coordinates are generated at random, then this disaster has
probability 0 (although in principle arbitrarily high-order boundaries
are still possible, so however large your stack is, it \emph{could}
overflow). If the aim is to plot an entire patch of tiling deflated
from some single $n$th-order supertile, then the disaster is \emph{by
definition} impossible, since either algorithm will report failure by
running out of input rather than running for ever.

To plot non-random instances of a tiling which involve infinite
supertile boundaries with \emph{eventually periodic} addresses, this
approach can still be used, but instead of running the transducer
until it reaches an accepting state, one runs it until it has
calculated the full eventually periodic address of the destination
tile, by the procedure described in section~\ref{periodic}. This also
allows infinite supertile boundaries to be detected, and displayed on
the diagram if desired.

What about plotting a \emph{fully general} patch of tiling, given only
a stream of symbols from an input tile address about which you know
nothing? In this case, neither of these approaches works:\ there is no
guarantee that a transition algorithm will reach an accepting state or
terminate its recursion, and even if it does begin repeating a
periodic cycle of output, the algorithm cannot know that, without
knowing the same thing about the input. In this situation, the simple
approach of computing each tile's address before moving on to the next
tile cannot work reliably.

An alternative algorithm in this situation is to propagate symbols of
the address one by one across the set of tiles generated so far. Keep
track of a tree of links between the tiles, with each tile linked to a
unique `parent', being the neighbour tile from which it was first
discovered. For each tile other than the starting tile, store the
currently known prefix of its address, and remember the current state
of a transducer that is generating its address from that of its
parent. Then, for each new symbol received from the input stream,
append it to the address of the starting tile, and propagate that
symbol outwards to all the tile's children by feeding it to each
child's transducer. If a transducer emits one or more output symbols
as a result, propagate those in turn to that tile's children, and so
on. A transducer is initialised for each exposed edge of a known tile,
and a new tile is added to the collection as soon as that transducer
generates the first output symbol for that tile, allowing its type and
position to be known.

At every stage of this algorithm, a finite number of symbols propagate
along each edge of this tree, so no intermediate step can fail to
terminate. And every desired tile must eventually be discovered,
because the lag between the deterministic transducer's input and
output is bounded by a constant $k$, so that if $n$ symbols of a
parent tile's address are known, then so are at least $n-k$ symbols of
its child. By induction, if a tile is separated from the starting tile
by $m$ edges of the tree, then at least $n-mk$ of its symbols can be
calculated from $n$ symbols of the starting tile. So for any finite
output region, a finite number of input symbols suffices to discover
at least the lowest-order type and position of every tile in the
region.

This family of techniques in general has many advantages over the
approach of deflating a single starting tile:

\begin{description}
  \item[Only generate the tiles you really need.] Generating a much
    larger patch than necessary, only to throw most of the tiles away,
    is a waste of computation. The transducer technique does not
    calculate the address of any tile until it is already known that
    the tile overlaps the target region of the plane.

    In some cases the deflation technique can reduce its computational
    costs by detecting early when a high-order intermediate tile will
    not overlap the output region at all, and discarding it from the
    list of tiles to deflate to the next level. This approach works
    well in systems such as Robinson triangles, where a supertile and
    its subtiles occupy exactly the same region of the plane, and that
    region is geometrically simple. But it works badly in the hat and
    Spectre tilings, where the $n$-times deflation of an order-$n$
    supertile becomes a more and more complicated shape and the
    geometry of the tiling is distorted during each deflation, so it
    is hard work to determine \emph{whether} a high-order supertile
    will contribute to the output region.

  \item[Ignores higher-order tile geometry completely.] The geometric
    distortion at every layer of the hat or Spectre tiling does not
    complicate the transducer approach at all. This entire suite of
    algorithms is completely unconcerned with the geometry of tiles,
    only with their combinatorics. It matters not at all to them
    whether the geometric effects of deflation are simple or complex.
    Even the reflecting nature of Spectre substitution systems, in
    which the whole plane reverses handedness in every deflation, is
    ignored as irrelevant.

  \item[Avoids dealing with very large or precise coordinates.] One
    way or another, the deflation technique needs to handle plane
    coordinates at high precision, because even if it does not
    \emph{generate} every single tile in the large patch corresponding
    to its starting tile, it must at least consider locations
    everywhere in that patch to decide which of them to process
    further. Storing and processing large numbers is inconvenient and
    error-prone:\ if the algorithm runs in floating point, then
    rounding errors become more significant as the coordinates become
    more precise, and if plane coordinates are instead represented
    exactly as an integer linear combination of some set of algebraic
    numbers, then the integer coefficients risk overflowing.

    The transducer technique never needs to consider any plane
    coordinate that is not a vertex of a tile in the final output. All
    the information about large supertiles is encapsulated in the
    combinatorial tile addresses, without dealing with their locations
    in the plane. Those tile addresses are essentially discrete, so it
    is easy for algorithms to handle them without any risk of
    imprecision. Overflow is still possible, but it manifests as
    running off the end of a finitely long address string and needing
    to extend it.

  \item[No need to commit to $n$ in advance.] In the deflation
    technique, you begin by choosing how large a patch to generate, by
    deciding on how high-order a supertile to begin from. In the
    transducer technique, you \emph{can} work the same way, but more
    flexible approaches are also possible.

    In particular, if you want to generate a \emph{random} patch of a
    substitution tiling, then there is no need to precommit to the
    length of the tile addresses at all. Instead, you can literally
    `make it up as you go along', by starting with a very short
    initial address, and extending it lazily as necessary:\ if the
    symbols of the starting tile address turn out to be insufficient
    to fill the whole output region, choose additional symbols at
    random to extend it, and proceed with tiling generation as if
    those additional symbols had been known all along.

  \item[Random tilings from the exact limiting distribution.]
    Continuing on the theme of generating a random patch of
    tiling:\ in a typical substitution tiling the relative density of
    tile types in the plane will converge to a limiting distribution
    the larger a patch you consider (namely, the eigenvector with the
    largest eigenvalue of the matrix that describes how many of each
    subtile type are produced by deflating each supertile type).

    Using the deflation technique for making a patch of tiling from a
    single $n$th-order supertile, you can approach this limiting
    distribution as closely as you like by choosing $n$ sufficiently
    large, but the computational cost increases as you do so, so you
    must trade off how closely you want to approximate the limiting
    distribution with the amount of work you are prepared to do.

    The transducer technique requires no tradeoff. When you invent
    each supertile symbol to extend the input address, you can choose
    that supertile symbol \emph{directly} from the exact limiting
    distribution.

  \item[Can handle infinite supertile boundaries.] As discussed
    already, transducer-based techniques can describe and successfully
    generate a region of tiling intersecting an infinite supertile
    boundary, because they do not depend on having any finite-order
    supertile in common between a pair of neighbouring tiles.
\end{description}

\subsection{Address-finding}

In the previous section we discussed the problem of turning a tile
address into a patch of tiling. Another application of tiling automata
is to go in the opposite direction. Suppose you are given a patch of
tiles covering a region of the plane, and you wish to find an address
that describes the same patch in terms of a particular substitution
system.

(Perhaps the patch of tiles was generated by a different substitution
system for the same tiling; or perhaps from a different method, such
as a cut-and-project scheme. Or perhaps you do not even know that, and
have merely received a patch in the form of a drawing, from some
unknown source.)

The input to this procedure is precisely the output of the previous
one:\ a collection of tiles, each one given its type, and the knowledge
of which pairs of tiles are adjacent, along which edge of each one. We
aim to recover an address describing one of the tiles in the patch --
any one will do -- such that regenerating a tiling from that address
as described in the previous section will produce the same layout of
tiles as the input.

To begin with, there is no hope of finding a \emph{unique} tile
address. The argument in the previous section is sufficient:\ if the
maximum input/output lag of the transducer for a substitution system
is $d$, and two tiles have a path of length $j$ between them, then
$jd$ symbols of the address of one tile are sufficient to determine
the type of the other. So a finite patch of tiling can only constrain
finitely many symbols of the address of any of its tiles.

Therefore, an address-finding algorithm must deal in a representation
of many possible tile addresses at once. A convenient representation
is an NFA.

As a first step, we use Algorithm~\ref{algadj} from
section~\ref{adjmatcher} to generate a recogniser $R$ for pairs of
neighbouring edge addresses. As mentioned in section~\ref{transducer},
this can also be regarded as a nondeterministic transducer taking one
edge address as input and producing another as output. (For these
purposes, a nondeterministic transducer is good enough, and it is
convenient that this form of automaton has no lag between the input
and the output.)

To analyse a patch of tiles, the essential plan is to focus on one
tile at a time, and have an NFA $A$ that matches all possible
addresses for that tile that we have not yet ruled out. We begin by
choosing a starting tile within the input patch, and initialising our
NFA $A$ to one that recognises any legal address at all for a tile of
that type. (This is easy to construct from first principles as a DFA,
and of course that is also a valid NFA.)

Next, we step across edge $e$ of the starting tile $t$, arriving in
its neighbour tile of type $u$ via edge $f$. We must now construct a
new NFA $A'$ for the possible addresses of that destination tile.

To do this, we construct the Cartesian product $A\times R$ of the
previous NFA with the neighbour-recogniser automaton. A state of
$A\times R$ is an ordered pair consisting of a state $a\in A$ and a
state $r\in R$; the product automaton has a transition
$(a,r)\xrightarrow{\alpha}(a',r')$ on a symbol $\alpha$ if there
exists $\beta$ such that the transition $a\xrightarrow{\beta}a'$
exists in $A$, and $r\xrightarrow{(\beta,\alpha)}r'$ exists in $R$.
That is, the NFA for $t$'s address accepts $\beta$ as the next symbol
from state $a$, and the neighbour recogniser accepts $\alpha$ as the
next symbol of $u$'s address if $\beta$ is the next symbol of $t$'s.
Finally, we delete every transition from the start state of the
product machine which does not take $\Edge{t}{e}$ as the first input
symbol and deliver $\Edge{u}{f}$ as the output symbol.

The resulting NFA describes a set of infinitely long tile addresses
which are consistent with being the address of a tile of type $u$,
whose neighbour via the edges $e,f$ is a tile of type $t$ described by
the input NFA $A$. We reduce this NFA by deleting unreachable states,
and translate to a canonical form (converting complex ordered-pair
representations of its states back to something computationally simple
like integers), to produce $A'$.

By repeating this procedure, we can walk around the entire input patch
of tiling, via a path that visits every tile. (A Hamilton cycle is not
necessary; visiting a tile more than once is harmless.) As we move
around, our NFA evolves in two ways:\ it changes to reflect the tile
we are currently standing on, and it also becomes more and more
refined as the algorithm incorporates more knowledge about nearby
tiles in the tiling.

Once we have visited all tiles, the NFA should be such that any valid
path through it describes an address consistent with the tiles
visited. For a reasonably large patch of tiling, it is likely that the
NFA will give the first few symbols of the output tile address
unambiguously.

If desired, at this point the output NFA can be determinised, via the
usual technique of making a DFA whose states correspond to sets of NFA
states. A DFA is easier to read, because it will begin with a unique
transition for each address symbol that is uniquely determined. (In
the NFA, one might have to trace many branches with the same symbols
on them to check if one of them had a different idea.)

The remaining challenge is to choose just one address. For some
purposes (depending on what the output is going to be used for), it is
enough to choose randomly, or arbitrarily. However, if the input patch
of tiling was a particularly interesting instance, such as a singular
or symmetric one, then probably one of the possible addresses is the
one most desired, and it is probably a simple, eventually periodic
address with a short period. The author has found that a good approach
is to end the walk around the tiling on a tile as close as possible to
the pattern's centre of symmetry, and then search the NFA to find the
\emph{simplest} eventually periodic address that it accepts, measured
by the combined length of its initial segment and repeating segment.
This in turn can be done by iterating over candidate pairs of lengths
in increasing order, and for each one, searching the automaton to find
whether any address of that shape is accepted at all.

\subsection{Analyses of tilings}

One more possibility, given a finite-state automaton describing a
tiling, is to transform the automaton in various ways to derive
automata that describe classes of addresses, or address pairs, of
particular interest.

For example, suppose one is interested in infinite supertile
boundaries in general. As mentioned in section~\ref{periodic}, an edge
between two tiles is part of an infinite supertile boundary if and
only if the neighbour address recogniser, given both tiles' addresses
as input, never enters an accepting state. Therefore, by deleting all
accepting states from the neighbour recogniser, a DFA is obtained
which matches precisely those pairs of addresses that lie on such a
boundary.

If two address pairs lie on the same infinite supertile boundary, then
they have some finite-order supertile in common, so their paths
through this DFA will at some point reach the same state. Therefore,
one can imagine an equivalence relation on address pairs in which two
pairs are considered equivalent if they differ in only finitely many
places; an infinite supertile boundary as a whole corresponds to an
equivalence class of this relation.

In particular, suppose we restrict our attention once again to the
useful class of \emph{eventually periodic} addresses. As mentioned in
section~\ref{periodic}, if an individual tile address is described as
an initial segment and a repeating segment, then its containing
infinite supertile is defined by the repeating segment alone. Pairs of
infinite supertiles bordering on each other are defined by a cyclic
path within the winnowed DFA. By searching for those cycles in
increasing order of length, it is possible to automate the process of
enumerating the infinite supertile boundaries in the
eventually-periodic class. An analysis similar to this is begun for
the Spectre tiling in \cite{Paolini}; this DFA technique permits it to
be automated for any tiling.

\subsubsection{Finding out why a substitution system was ambiguous}\label{ambiguityfinding}

If the substitution system was originally ambiguous, and
Algorithm~\ref{algrefine} was used to turn it into an unambiguous one
admitting a transducer, then a different analysis of this kind can be
used to find addresses in the unrefined system which correspond to
more than one address in the refined one -- in other words, to find
the particular hard-to-distinguish tilings which \emph{gave rise} to
the ambiguity in the substitution system.

Some of this information can be extracted from an attempt to
determinise a transducer for the original (coarse) tiling. If the
section~\ref{detectfailure} technique is used to detect that
determinisation has failed, then a by-product is the ability to
identify a particular input to the transducer with two possible
outputs, via the cycle found in the graph $G$ and a path to that cycle
from the start state. However, this only describes \emph{one}
ambiguous input -- whichever was the first one to be detected during
the failed transducer construction. The following technique will find
them all systematically.

To do this, we first make a DFA which matches the `diagonal' subset
$\mathcal{D}\subset\mathcal{T}\times\mathcal{T}$ of pairs of tile
addresses in the refined system: the language of strings made by
zipping the \emph{same} tile address with itself. This is matched by a
DFA isomorphic to the one for $\mathcal{T}$ itself, consisting of a
\textsc{start} state, and a state for each $(\textit{layer},
\textit{tile type})$ pair, with transitions derived from the
substitution system in an obvious manner. We aim to find a pair of
distinct tile addresses $A,A'$ accepted by this DFA which have the
same image under the `forgetful' map $\phi$ which discards tile
subtypes to recover the corresponding address in the original
unrefined system.

If $A,A'$ have the same subtype for their $n$th-order supertiles, for
some $n$, then they must exactly match in all lower-order tile
subtypes as well, because the subtype of any tile together with a
subtile index uniquely determines the type of the next lower-order
tile, and so on inductively. So if $A\neq A'$, then there must be some
$k\geq0$ such that the two addresses exactly match in the first $k$
symbols, and thereafter, have \emph{no} symbol in common.

Now consider what happens if we zip together the addresses $A,A'$ into
a combined string containing ordered pairs of a symbol from each
address, and use that as input to the product DFA
$\mathcal{D}\times\mathcal{D}$. While processing the initial $k$
symbols, the product DFA will remain in a `diagonal' state, i.e.\ one
of the form $(s,s)$. After that, every transition will be on a pair of
symbols which match in their coarse tile type, and subtile index, but
\emph{never} match in their tile subtype from refinement.

We can construct a subset of $\mathcal{D}\times\mathcal{D}$ matching
all such strings by starting from every diagonal state and using only
transitions on symbols of the form $\Bigl(\Tile{t},\Tile{t'}\Bigr)$ or
$\Bigl(\Sup{i}{t},\Sup{i}{t'}\Bigr)$ for tiles $t,t'$ which are
distinct refinements of the same coarse tile type. Any infinite path
in this subset DFA corresponds to the tail of some pair of ambiguous
tile addresses, and a full ambiguous tile address can be constructed
by prefixing to it any path from the DFA start state to the diagonal
state it began from.

\section{Results}

In this section we shall present some results of applying these
algorithms to well-known tilings. Full-page figures at the end of the
article give the complete details of each substitution system suitable
for input to these algorithms: the edge types and edge deflation rules
are shown (constituting the `geometric source code' described in
section~\ref{source}), and each deflation diagram is shown in full,
with both sides of each edge labelled with an appropriate $\Ext{u}{v}$
or $\Int{i}{j}$ symbol, so that the full combinatorial adjacency maps
can be read off these diagrams, or cross-checked against them if
regenerating them from the geometry.

Where a refinement of the system is needed to make it unambiguous,
that refinement is also shown in a figure. The refinement diagrams do
not need to repeat the detailed adjacency maps from their source
tiling, since those are unchanged; they only need to show how many
subtypes of each original tile type exist, and which subtypes appear
in each deflation.

\subsection{Penrose P2 and P3 tilings}\label{penrose}

The Penrose P2 tiling of kites and darts, and the P3 tiling of thin
and thick rhombs, can each be generated by either a substitution
system of Robinson triangles, or one made from the whole tiles. Both
systems have strengths. The Robinson-triangle systems are
geometrically simple in that each tile's deflation occupies the same
region of the plane as the original tile (if scaled appropriately),
and preserve all the symmetry in the original tilings. On the other
hand, the whole-tile system is naturally more practical if your aim is
to generate actual Penrose tilings for display, since it does not
require a post-processing step of recombining the two triangular
halves of each tile.

(As a practical matter, this recombination could be done by a separate
software layer consuming the output of a transducer, but another
approach would be to augment the substitution system with a one-off
base layer whose parent is the main layer of Robinson triangles,
deflating one half of each tile into the tile itself and the other
half into an outline doubling back on itself and enclosing no area.)

All four of these substitution systems -- \{ P2 and P3 \} $\times$ \{
triangles, whole tiles \} -- are unambiguous.
Algorithm~\ref{determinise} successfully constructs a deterministic
transducer for each one.

Figures~\ref{p2-triangles-prototiles}, \ref{p2-triangles-edge}
and~\ref{p2-triangles-adjmap} show the P2 triangle system in full;
figures~\ref{p2-whole-prototiles}, \ref{p2-whole-edge}
and~\ref{p2-whole-adjmap} show the P2 whole-tile system;
figures~\ref{p3-triangles-prototiles}, \ref{p3-triangles-edge}
and~\ref{p3-triangles-adjmap} show the P3 triangle system;
figures~\ref{p3-whole-prototiles}, \ref{p3-whole-edge}
and~\ref{p3-whole-adjmap} show the P3 whole-tile system.

\subsection{Ammann-Beenker tiling}\label{ammann-beenker}

The Ammann-Beenker tiling is usually described as a two-tile system,
consisting of a square and a rhombus. However, the rhombus comes in
two mirror-image forms (although these are indistinguishable if you do
not show the vertex markings that force the tiling to be aperiodic).
For the purposes of these algorithms, mirror-image tiles must be
counted as different.

When choosing a substitution system for the Ammann-Beenker tiling, one
also has a similar choice to the Penrose tilings:\ either bisect the
square tile into a pair of right triangles so that each tile deflates
to a patch of subtiles occupying the same area of the plane, or leave
the square whole and instead complicate the shapes of the deflation
maps.

As mentioned in section~\ref{refinement}, the Ammann-Beenker tiling is
ambiguous, because when eight rhombus tiles appear around a single
vertex, they can only do so in one way, which requires a particular
arrangement of the two mirror-image forms (see
Figure~\ref{ammannbeenker8}). This is equally true of either
substitution system. Moreover, both systems require two passes of
Algorithm~\ref{algrefine} to refine them from the original three- or
four-tile system into one that admits a deterministic transducer.

The unambiguous refinements of both systems involve eight variants of
the rhombus (four of each of the two mirror-images). When the square
is left whole, there are eight variants of that too; when it is
bisected, each half-square splits into \emph{five} variants, which can
combine in a subset of the possible combinations.

Figures~\ref{ammann-beenker-triangles-prototiles},
\ref{ammann-beenker-triangles-edge}
and~\ref{ammann-beenker-triangles-adjmap} show the Ammann-Beenker
substitution system involving two rhombus tiles and the two triangular
halves of the square, and
Figure~\ref{ammann-beenker-triangles-refined} shows its unambiguous
refinement. Figures~\ref{ammann-beenker-whole-prototiles},
\ref{ammann-beenker-whole-edge} and~\ref{ammann-beenker-whole-adjmap}
show the system with a whole square, and
Figure~\ref{ammann-beenker-whole-refined} shows the unambiguous
refinement of that system.

\subsection{Spectre tiling}\label{spectre}

\cite{Smith_Sep2024} presents two substitution systems for the Spectre
tiling. One, shown in figure~2.1 of the paper, is a minimalist
system containing two tile types:\ a single Spectre, and the `Mystic'
double tile consisting of two Spectres at $30^\circ$ to each other.
The other, spread across figures~4.2, 5.1, and~5.3, involves nine
hexagonal metatiles and eight different types of edge, with one set of
deflation rules turning hexagons into more hexagons, and a second set
turning them into Spectre tiles.

Both of these systems can be represented in the combinatorial form
required as input to these algorithms. The 9-hexagon system is
essentially in that form already, but the simpler Spectre/Mystic
system is more difficult to represent, because the paper does not
state how the edges of a single Spectre or Mystic correspond to the
edges of the patch of tiles deflated from each one. However, such an
edge mapping does exist; information equivalent to it is available in
various sources such as \cite{Cheritat}, and the author is indebted to
Pieter Mostert in particular for discussions involving the related
`comet-rhomb tiling' from which it can be derived.

To specify the edge mapping for the Spectre/Mystic system, we first
classify the edges into two distinct types, and assign a direction to
each one. The edge types correspond to the direction of the edge in
the plane: edges whose direction is a multiple of $60^\circ$ have one
type, and those at an odd multiple of $30^\circ$ the other type.
Directions are assigned in alternation around the 14 edges of the tile
(counting the long edge as two collinear edges, as usual).

This assignment of edge types separates the previously identical
Spectre tiles into two distinct types, again depending on their
orientation mod $60^\circ$, with most Spectres being of one type, and
a small number being the `skew' Spectres with all the edge types
swapped. Each skew Spectre is then fused with one of its neighbours to
make the Mystic double-tile. Figure~\ref{spectreh7h8tiles} shows the
resulting two tile types, with their edge markings and directions.

With the edges distinguished in this way, it is now possible to write
down a pair of deflation rules that turn each edge consistently into a
pattern of other edges. Figure~\ref{spectreh7h8rules} shows these
deflation rules.

If the Spectre and Mystic outlines are deflated by applying those
rules, each one generates an output region which can be tiled with a
number of Spectres, including one skew one forming half of a Mystic.
However, since the Spectre deflation reflects the plane, these two
regions must either be tiled with reflected Spectres, or be reflected
themselves before tiling with Spectres in the original handedness. (We
have chosen the latter approach in our diagrams.)

The deflated outline also includes some zero-thickness spurs. In both
outlines a long spur of six edges protrudes outwards from the patch
outline, and then retraces its steps to its starting point. The Mystic
deflation also includes two further spurs, one intruding \emph{inward}
into the region and crossing over interior edges of its tiling with
Spectres, and one running back and forth along part of the outline.

Figures~\ref{spectreh8} and~\ref{spectreh7} show the deflation
diagrams for the two tiles in full\footnote{As a practical matter,
when generating a Spectre tiling from this system, one could add a
base layer to the substitution system which separates the Mystic tile
into its two Spectres. This is a simple transformation, and left to
the reader.}. The difficult spurs in the Mystic deflation are shown
detached from the main tiling, with an indication of where they should
properly start from.

The algorithms in this article take all of this complication in their
stride. The adjacency maps for the two deflations describe the
complicated spurs in the same way as described in section~\ref{spurs},
mapping some $\Ext{u}{v}$ to each other, and whenever two pairs of
edges coincide, making an arbitrary choice about which one to map to
which\footnote{The adjacency map in Figure~\ref{spectreh7} shows a
particular choice of which edges to map to which, but this is not
critical to the algorithm. For example, the diagram shows the mappings
$\Ext{15}{-4}\leftrightarrow\Ext{12}{-2}$ and
$\Int{2}{9}\leftrightarrow\Int{3}{4}$ on a pair of edges which
coincide if the separate spur is replaced at its proper starting point
(shown by a red diamond). The alternative choice
$\Ext{15}{-4}\leftrightarrow\Int{3}{4}$ and
$\Int{2}{9}\leftrightarrow\Ext{12}{-2}$, mapping the top side of each
of those edges to the bottom side of the \emph{other} one, would work
equally well and make no difference to the success of the
algorithms.}.

The alternative 9-hexagon system is not shown in full here, partly
because \cite{Smith_Sep2024} gives all the details required, and also
because it is redundant, for the following reason.

The 9-hexagon system is unambiguous. If it is transcribed into a
combinatorial description in the style of this article, with a
`Spectres' layer and a `hexagonal metatiles' layer and two sets of
deflation rules, then that description admits a deterministic
transducer.

On the other hand, the Spectre/Mystic system is ambiguous: it does not
admit a deterministic transducer. However, one pass of
Algorithm~\ref{algrefine} refines it into an unambiguous system
consisting of nine tile types: eight different variants of the single
Spectre tile, and still only one Mystic.

\begin{figure}
\centering
\includegraphics*[alt={Nine deflation diagrams, eight for Spectre tiles marked with different Greek letters, and one for the Mystic tile. Each tile is deflated to a pattern of other Spectres and one Mystic, again marked with the same Greek letters to indicate which tile types are generated from which.}]{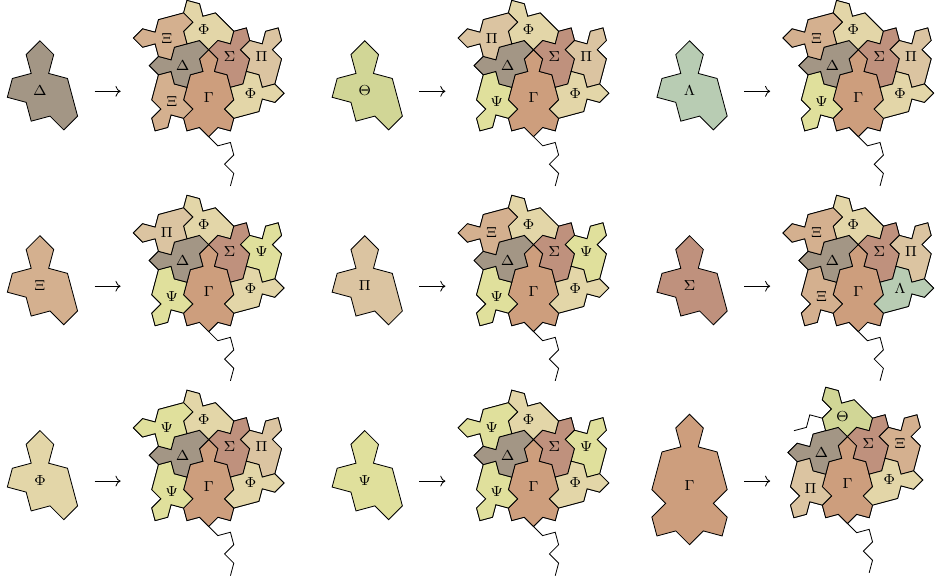}
\caption{The 9-tile refinement of the Spectre/Mystic substitution system}
\label{spectre9}
\end{figure}

The unambiguous system generated by Algorithm~\ref{algrefine}
\emph{exactly} corresponds to the 9-hexagon system: each of its 9 tile
types corresponds to one of the hexagons, and every deflation diagram
in the refined system matches the corresponding diagram in Figure~5.1
of \cite{Smith_Sep2024}, \emph{except} that the metatile shapes are
Spectres and Mystics instead of hexagons. Figure~\ref{spectre9} shows
these deflation diagrams in full, with the same tile names (and
colours) as \cite{Smith_Sep2024} uses for the corresponding hexagons.

Moreover, a deterministic transducer constructed from the 9-hexagon
system is \emph{equivalent} to one constructed from the refined
Spectre/Mystic system shown here. Having labelled the tile types
compatibly, tile addresses in the two systems look essentially the
same (up to differences of detail such as whether the lowest-order
Spectre is labelled directly with one of the nine types, or whether it
is deflated from a lowest-order hexagon with the same type), and the
two transducers agree on the neighbour address for any input tile.
This occurs because \emph{after} the transducer is built, there are no
remaining references to the edges of any supertile, so the fact that
the supertiles are hexagonal in one system and Spectre-shaped in the
other is no longer visible.

So applying the refinement algorithm to the Spectre tiling has not
generated a new and interesting substitution system. Instead, it has
re-derived one of the systems in the paper from the other one (with
the small difference of more complicated tile shapes but fewer edge
types). However, even without any new discoveries, the fact that the
two analyses agree is a useful check that the algorithms presented
here are behaving as expected.

Also, it is useful to have re-derived the 9-tile system as a
refinement of the Spectre/Mystic system, because we now have an
explicit mapping between the two systems. Therefore it becomes
possible to apply the procedure from Section~\ref{ambiguityfinding} to
find out details of \emph{why} the Spectre/Mystic system did not
already admit a deterministic transducer: what infinite-order
supertile(s) have multiple representations?

Running that algorithm reports that there is just \emph{one} ambiguous
infinite supertile, with three distinct addresses in the refined
system which all look the same in the unrefined Spectre/Mystic system.
In the unrefined system, the `simplest' address of a tile in the
ambiguous supertile is a Spectre tile which is subtile 7 of a
larger Spectre (in the arbitrary indexing of subtiles used in
Figure~\ref{spectreh8}), which is subtile 7 of its still larger
parent, and so on for ever:
$$\hbox{spectre}\xleftarrow{7}\hbox{spectre}\xleftarrow{7}\hbox{spectre}\xleftarrow{7}\hbox{spectre}\xleftarrow{7}\cdots$$
In the refined system, this address separates into three. One of them
has the same tile type $\Psi$ at every level, so that it is invariant
under deflation. The other two alternate between the $\Xi$ and $\Pi$
tile types, so that each one deflates to the other:
$$\begin{aligned}
\Psi\xleftarrow{7}\Psi\xleftarrow{7}\Psi\xleftarrow{7}\Psi\xleftarrow{7}\cdots \\
\Xi\xleftarrow{7}\Pi\xleftarrow{7}\Xi\xleftarrow{7}\Pi\xleftarrow{7}\cdots \\
\Pi\xleftarrow{7}\Xi\xleftarrow{7}\Pi\xleftarrow{7}\Xi\xleftarrow{7}\cdots \\
\end{aligned}$$

In fact, the two infinite supertiles with $\Pi,\Xi$ tiles are
neighbours. Specifically, a transducer constructed from this
substitution system reports that each of the following two edge
addresses maps to the other one:
$$\begin{aligned}
\Edge{\Xi}{0}\xleftarrow{7}\Pi\xleftarrow{7}\Xi\xleftarrow{7}\Pi\xleftarrow{7}\cdots \\
\Edge{\Pi}{9}\xleftarrow{7}\Xi\xleftarrow{7}\Pi\xleftarrow{7}\Xi\xleftarrow{7}\cdots \\
\end{aligned}$$

So these two supertiles occur in the same Spectre tiling of the plane,
side by side. This is interesting because it means that some rigid
motion of the plane -- in fact, a $\frac16$ rotation -- maps one of
these identical supertiles to the other, and therefore, leaves an
infinite region of the plane invariant. Under the same rotation, what
happens to the rest of the plane?

A plausible answer would be that the whole plane is symmetric under
that rotation. There are Spectre tilings with a rotational symmetry.
However, there are none with \emph{6-way} rotational symmetry, only
2-way and 3-way. The answer in this case is more interesting: the
Spectre tiling containing these two infinite supertiles is singular in
the same sense as the Penrose infinite cartwheel pattern, in that
\emph{almost} all of the plane is invariant under that $60^\circ$
rotation, with only a small number of narrow paths of tiles varying.
Figure~\ref{spectre-px} shows the central section of this singular
tiling (which can be found by other methods too).

The third version of the same infinite supertile, involving a $\Psi$
tile at every level, occurs in a separate but extremely similar
Spectre tiling. It also has $60^\circ$ near-symmetry, but moreover, if
the two tilings are overlaid so that the $\Psi$ supertile corresponds
to one of the $\Pi,\Xi$ supertiles, then they match each other almost
completely, with only one of the asymmetric paths varying between the
two tilings. (Which path varies depends on which of $\Pi,\Xi$ and
$\Xi,\Pi$ you overlaid the $\Psi$ supertile on.)

\subsection{Hat tiling}\label{hat}

\cite{Smith_Jul2024} presents two substitution systems for the hat
tiling. One, shown in figure~2.11 of the paper, is a minimalist system
very similar to the Spectre/Mystic system, consisting of two tiles, of
which one is a single hat, and the other is a fusion of a reflected
hat with one of its neighbours, in exactly the same way that the
Mystic is a fusion of the rare skew Spectre with a normally oriented
neighbour. These deflate to patches of hats labelled $H_8$ and $H_7$
respectively. The other system consists of four metatiles $H,T,P,F$ in
simple polygonal shapes, with one deflation rule to turn them into
more metatiles, and another to turn them into hats (shown in
figures~2.8 and~4.1 of the paper respectively).

As with the Spectre/Mystic system, in order to handle either of these
hat substitution systems using the techniques in this paper, we first
need to know exactly how the edges match up between each tile type and
its deflation.

The rules for the $\{H_7,H_8\}$ system are very similar to those for
the Spectre/Mystic system: as in that system, we classify the edges
into two types based on orientation\footnote{In the hat tiling, this
classification of edge types also matches the difference in length:
the hat tiling inherits from its underlying kite tiling the property
that an edge of length 1 and one of length $\sqrt{3}$ always differ in
orientation by an odd multiple of $30^\circ$.}, and assign directions
in alternation around the tile. Figure~\ref{hath7h8tiles} shows the
two tiles marked in this fashion, and Figure~\ref{hath7h8rules} shows
the deflation rules for the two edge types.

Transforming each edge of the two tile types via these rules expands
each one into a patch which can be tiled by one double-hat fusion and
either 6 or 7 ordinary hats. Figures~\ref{hath8} and~\ref{hath7} show
the results of this process\footnote{As with the Spectre system, one
might wish to add a base layer to this system separating the
double-hat tile back into its two individual hats. Another useful
optional base layer, in this or the {\HTPF} system, is one that
divides each hat into the eight kites of the underlying kite tiling,
which allows the output of these algorithms to be handled in a
discrete manner as subsets of the kites in that grid.}.

The {\HTPF} system is presented in two forms. In one form, the
deflation diagrams of the metatiles overlap; in the other, selected
metatiles are removed from each deflation diagram so that each tile
has a unique supertile. For these algorithms, only the non-overlapping
form is suitable.

Edge annotations for the non-overlapping {\HTPF} system are derived in
\cite{Baake_2025} (figures~4 and~5), and reproduced here in
Figure~\ref{htpf-prototiles} (showing the edges subdivided as
mentioned in section~\ref{edgesubdivision}, with edge types and
directions), and Figures~\ref{htpf-edge-meta-meta}
and~\ref{htpf-edge-meta-hats} showing the rules for transforming each
edge type into a sequence of sub-edges during deflation of metatiles
to other metatiles or hats respectively. The resulting adjacency maps
are shown in Figures~\ref{htpf-adjmap-meta-meta}
and~\ref{htpf-adjmap-meta-hats}.

\emph{Both} of these substitution systems are ambiguous. In each case,
refinement via Algorithm~\ref{algrefine} resolves the ambiguity by
creating a system with more tile types.

\begin{figure}
\centering
\includegraphics*[alt={Deflation diagrams for a refinement of the HTPF substitution system, containing two different types of H tile and three F.}]{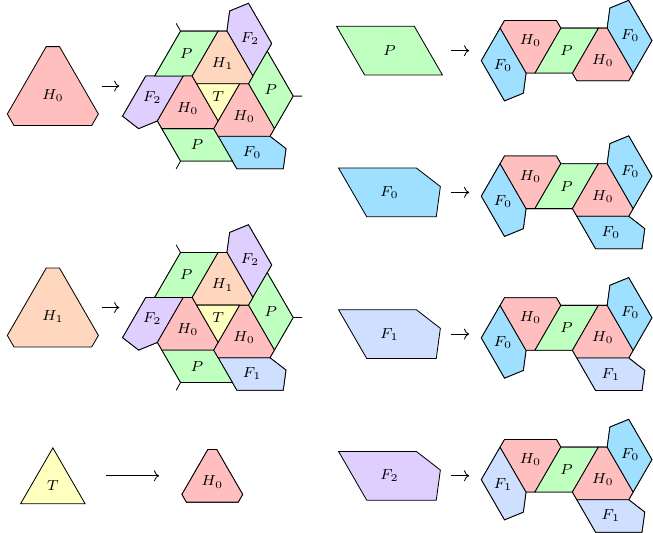}
\caption{The {\HHTPFFF} substitution system}\label{hhtpfff}
\end{figure}

The unambiguous version of {\HTPF} is a substitution system one might
call {\HHTPFFF}: it consists of two variants of the original $H$ tile,
three variants of $F$, and just one variant of each of $T$ and $P$.
These are shown in Figure~\ref{hhtpfff}.

\begin{figure}
\centering
\includegraphics*[alt={Deflation diagrams for seven types of single hat, distinguished by numbers 1 to 7 and colours, and one double-hat.}]{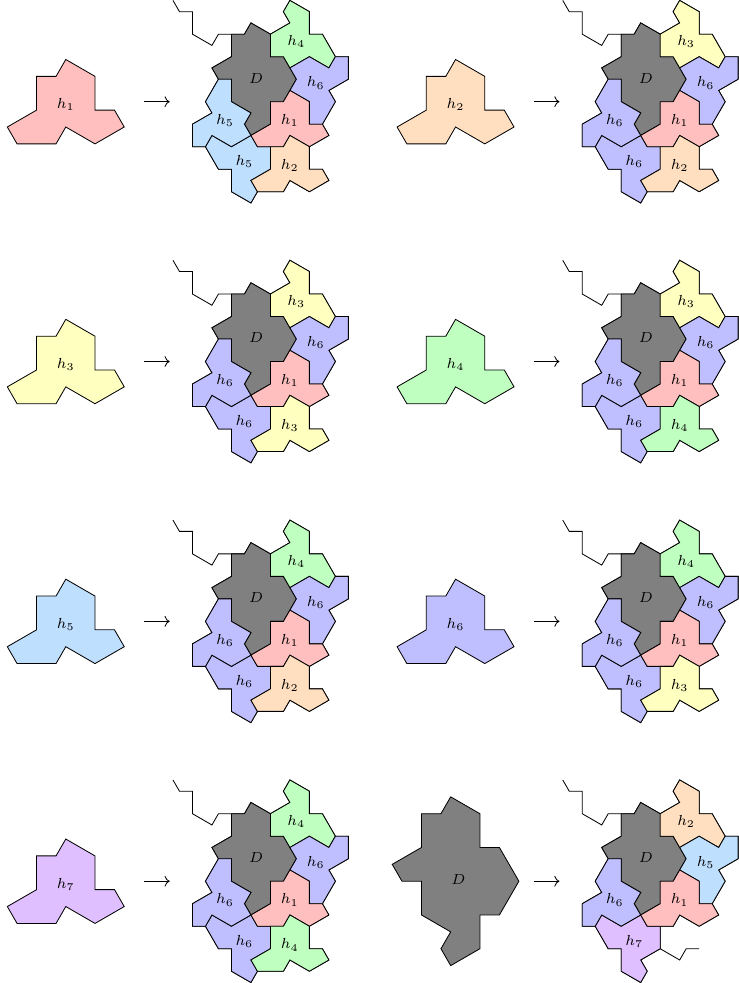}
\caption{The 8-tile refinement of the $H_7/H_8$ system for hats}\label{hat8}
\end{figure}

\begin{figure}
\centering
\includegraphics*[alt={Diagrams of the four triangles obtained by bisecting the Penrose P2 kites and darts. The edges are marked with four different types, via colours and arrows. The three- and four-arrow edges force the triangles to rejoin into whole kites and darts, and the one- and two-arrow edges match the ones in previous figures.}]{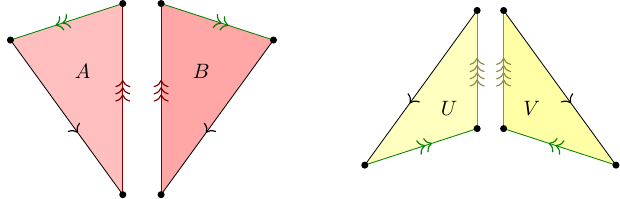}
\caption{Diagrams of the P2 Robinson-triangle tiles, with labelled edges}
\label{p2-triangles-prototiles}
\end{figure}

\begin{figure}
\centering
\includegraphics*[alt={Mapping rules showing how each of the edge types from the previous figure is transformed into a pattern of other edges during deflation.}]{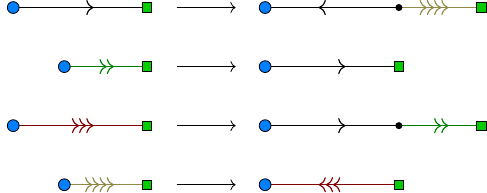}
\caption{Edge deflation rules for P2 Robinson triangles}
\label{p2-triangles-edge}
\end{figure}

\begin{figure}
\centering
\includegraphics*[alt={A detailed map of each P2 half-tile triangle and how it deflates into smaller triangles, marked up with symbols on both sides of each edge.}]{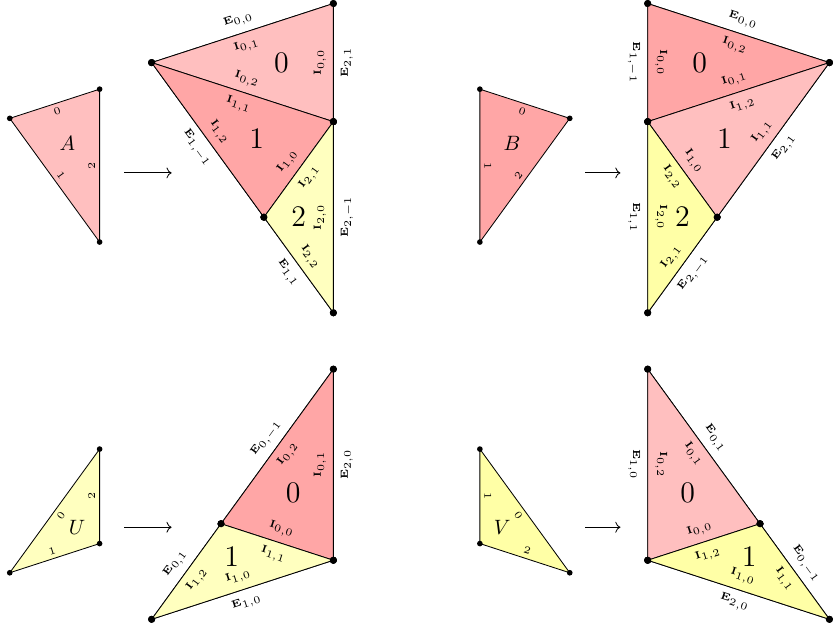}
\caption{Adjacency maps for P2 Robinson triangles}
\label{p2-triangles-adjmap}
\end{figure}

\begin{figure}
\centering
\includegraphics*[alt={Diagrams of the four triangles obtained by bisecting the Penrose P3 rhombs. The edges are marked with four different types, via colours and arrows. The three- and four-arrow edges force the triangles to rejoin into whole rhombs, and the one- and two-arrow edges match the ones in previous figures.}]{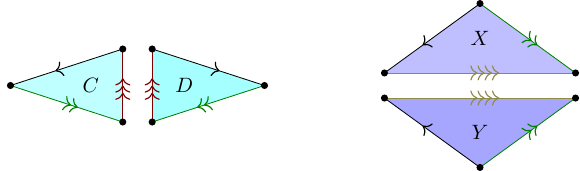}
\caption{Diagrams of the P3 Robinson-triangle tiles, with labelled edges}
\label{p3-triangles-prototiles}
\end{figure}

\begin{figure}
\centering
\includegraphics*[alt={Mapping rules showing how each of the edge types from the previous figure is transformed into a pattern of other edges during deflation.}]{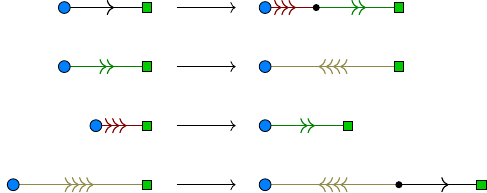}
\caption{Edge deflation rules for P3 Robinson triangles}
\label{p3-triangles-edge}
\end{figure}

\begin{figure}
\centering
\includegraphics*[alt={A detailed map of each P3 half-tile triangle and how it deflates into smaller triangles, marked up with symbols on both sides of each edge.}]{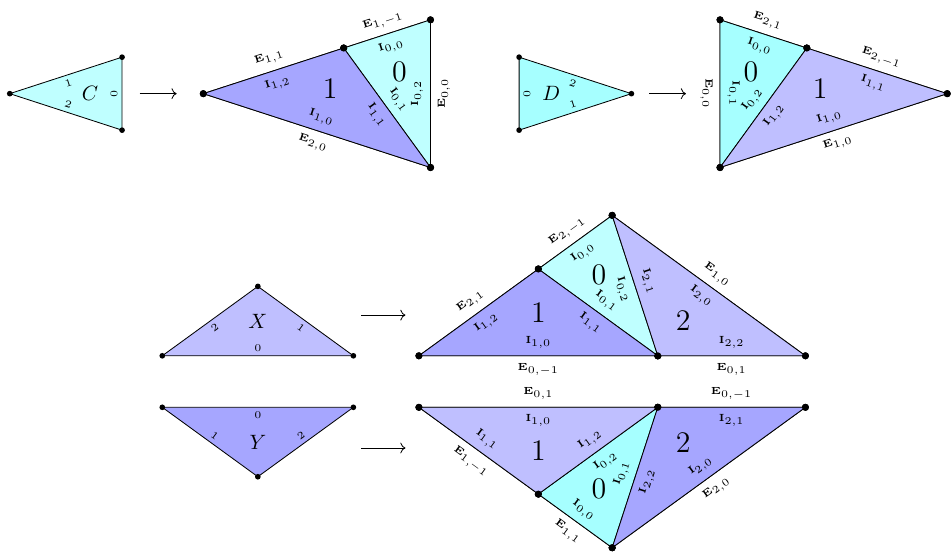}
\caption{Adjacency maps for P3 Robinson triangles}
\label{p3-triangles-adjmap}
\end{figure}

\begin{figure}
\centering
\includegraphics*[alt={Diagrams of the Penrose P2 kite and dart, whole rather than bisected. The edges are marked with two different types, via colours and arrows.}]{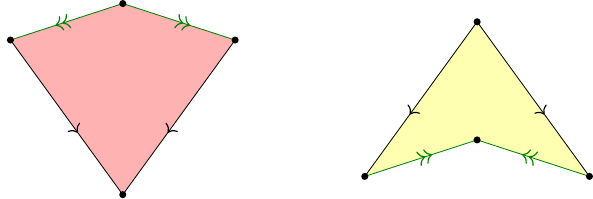}
\caption{Diagrams of the P2 whole tiles, with labelled edges}
\label{p2-whole-prototiles}
\end{figure}

\begin{figure}
\centering
\includegraphics*[alt={Mapping rules showing how each of the edge types from the previous figure is transformed into a pattern of other edges during deflation. The pattern obtained from the single-arrow edge has angles between the output edges.}]{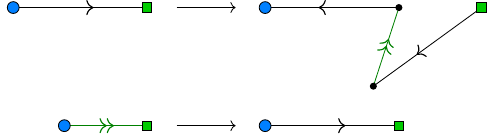}
\caption{Edge deflation rules for P2 whole tiles}
\label{p2-whole-edge}
\end{figure}

\begin{figure}
\centering
\includegraphics*[alt={A detailed map of each P2 tile and how it deflates into smaller tiles, marked up with symbols on both sides of each edge.}]{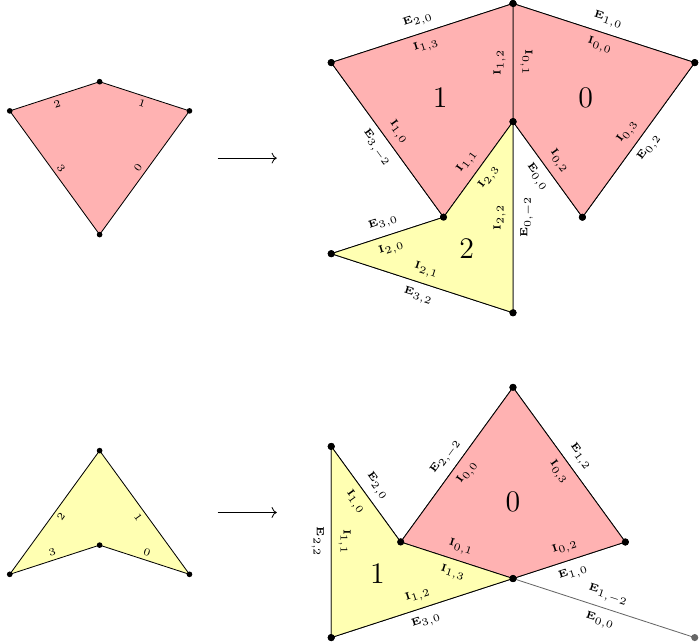}
\caption{Adjacency maps for P2 whole tiles}
\label{p2-whole-adjmap}
\end{figure}

\begin{figure}
\centering
\includegraphics*[alt={Diagrams of the Penrose P3 rhombs, whole rather than bisected. The edges are marked with two different types, via colours and arrows.}]{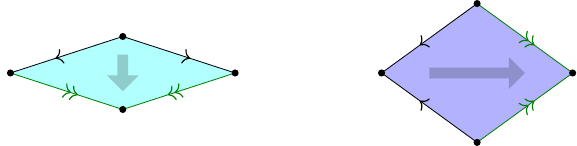}
\caption{Diagrams of the P3 whole tiles, with labelled edges}
\label{p3-whole-prototiles}
\end{figure}

\begin{figure}
\centering
\includegraphics*[alt={Mapping rules showing how each of the edge types from the previous figure is transformed into a pattern of other edges during deflation.  Both of the output patterns have angles between the edges.}]{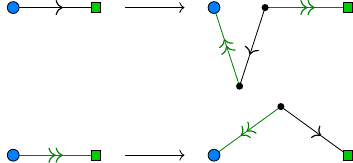}
\caption{Edge deflation rules for P3 whole tiles}
\label{p3-whole-edge}
\end{figure}

\begin{figure}
\centering
\includegraphics*[alt={A detailed map of each P3 tile and how it deflates into smaller tiles, marked up with symbols on both sides of each edge.}]{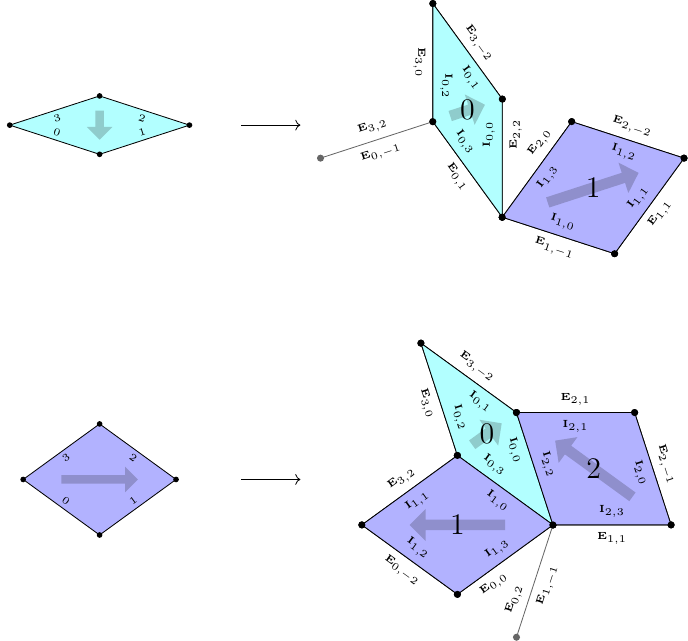}
\caption{Adjacency maps for P3 whole tiles}
\label{p3-whole-adjmap}
\end{figure}

\begin{figure}
\centering
\includegraphics*[alt={Diagrams of the Ammann-Beenker tiles, in a form containing two mirror-image rhombi and two right triangles obtained by bisecting the square tile. The edges are marked with two different types via colours and arrows, as in an earlier diagram. The two-arrow edges force the two types of triangle to join back into the whole square. All remaining edges have the same type.}]{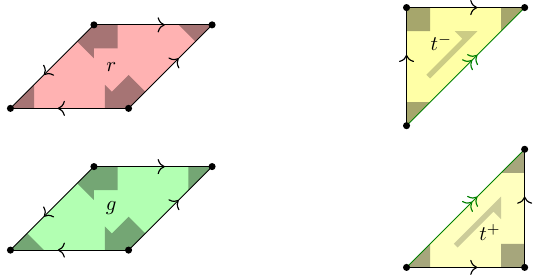}
\caption{Diagrams of the Ammann-Beenker tiles (with the square
  bisected), with labelled edges}
\label{ammann-beenker-triangles-prototiles}
\end{figure}

\begin{figure}
\centering
\includegraphics*[alt={Mapping rules showing how each of the edge types from the previous figure is transformed into a pattern of other edges during deflation.}]{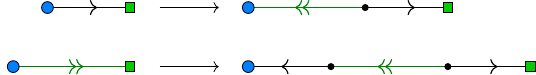}
\caption{Edge deflation rules for the Ammann-Beenker tiling}
\label{ammann-beenker-triangles-edge}
\end{figure}

\begin{figure}
\centering
\includegraphics*[alt={A detailed map of each Ammann-Beenker tile or half-tile and how it deflates into smaller tiles, marked up with symbols on both sides of each edge.}]{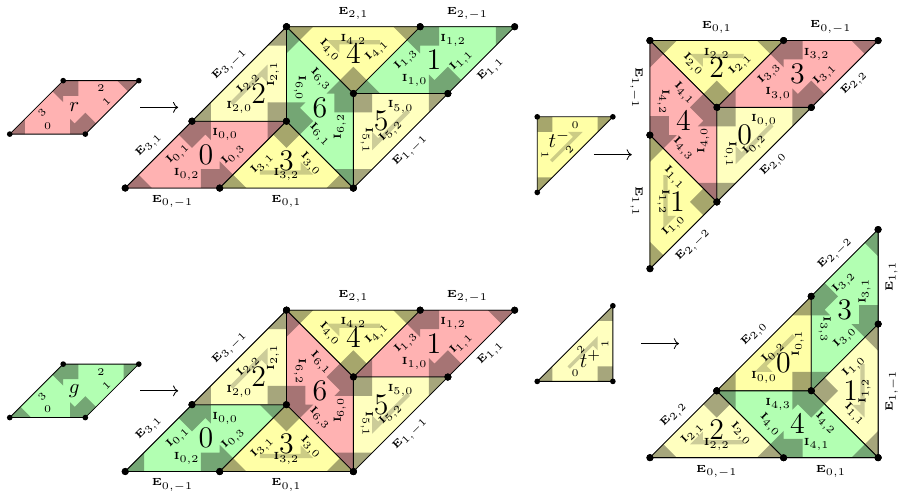}
\caption{Adjacency maps for the Ammann-Beenker tiling}
\label{ammann-beenker-triangles-adjmap}
\end{figure}

\begin{figure}
\centering
\includegraphics*[alt={Diagrams similar to the deflation maps in the previous figure, with the detailed edge labels omitted since they are unchanged, but with each tile cloned into multiple subtypes, distinguished by subscripts and slightly different colours, showing which subtypes are deflated from which other ones.}]{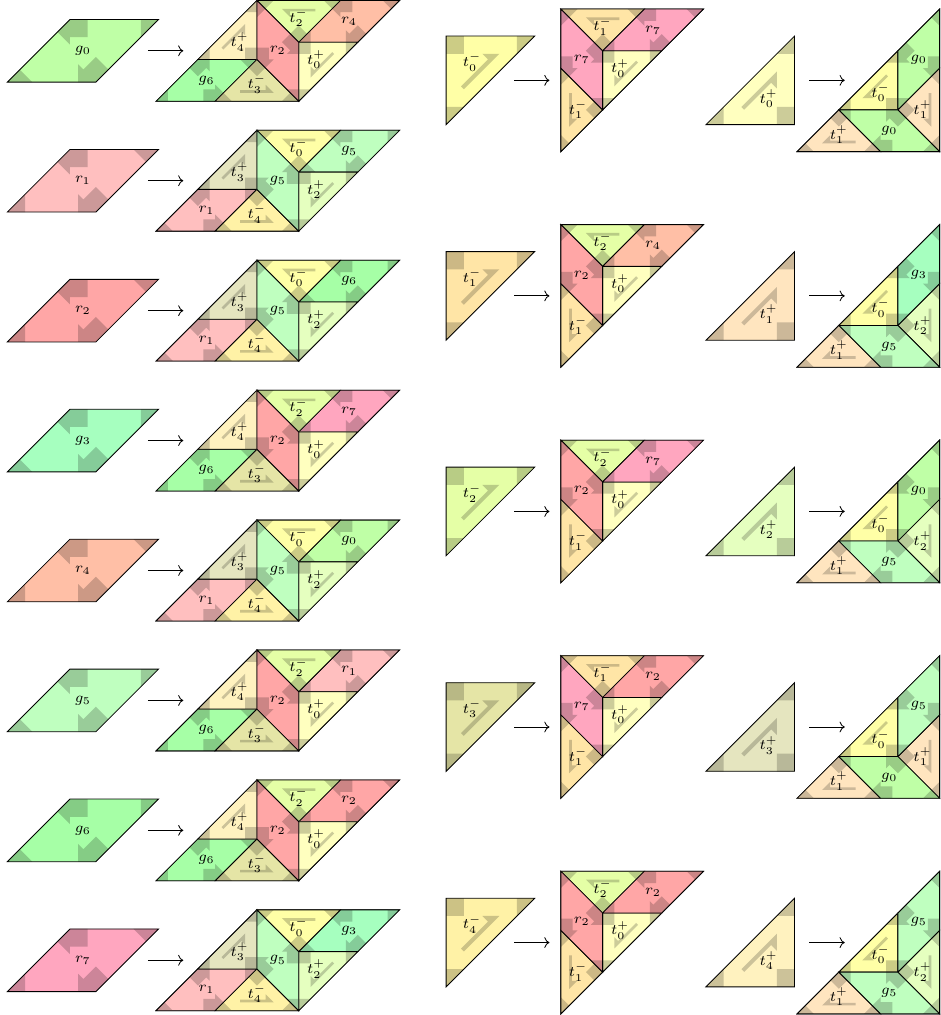}
\caption{Unambiguous refinement of the Ammann-Beenker system, with the square bisected}\label{ammann-beenker-triangles-refined}
\end{figure}

\begin{figure}
\centering
\includegraphics*[alt={Diagrams of the Ammann-Beenker tiles, in a form containing two mirror-image rhombi and a whole square tile. The edges all have the same type, but are marked with arrows indicating direction.}]{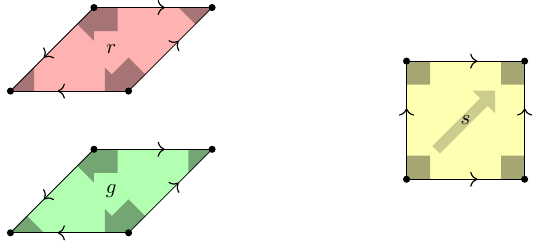}
\caption{Diagrams of the Ammann-Beenker tiles with the square whole,
  with labelled edges}
\label{ammann-beenker-whole-prototiles}
\end{figure}

\begin{figure}
\centering
\includegraphics*[alt={A single mapping rule showing how the directed edges from the previous figure is transformed into a pattern of other edges, with angles between them, during deflation.}]{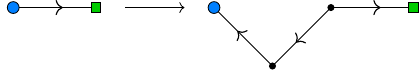}
\caption{The single edge deflation rule for the Ammann-Beenker tiling
  with the square whole}
\label{ammann-beenker-whole-edge}
\end{figure}

\begin{figure}
\centering
\includegraphics*[alt={A detailed map of each Ammann-Beenker tile and how it deflates into smaller tiles, marked up with symbols on both sides of each edge.}]{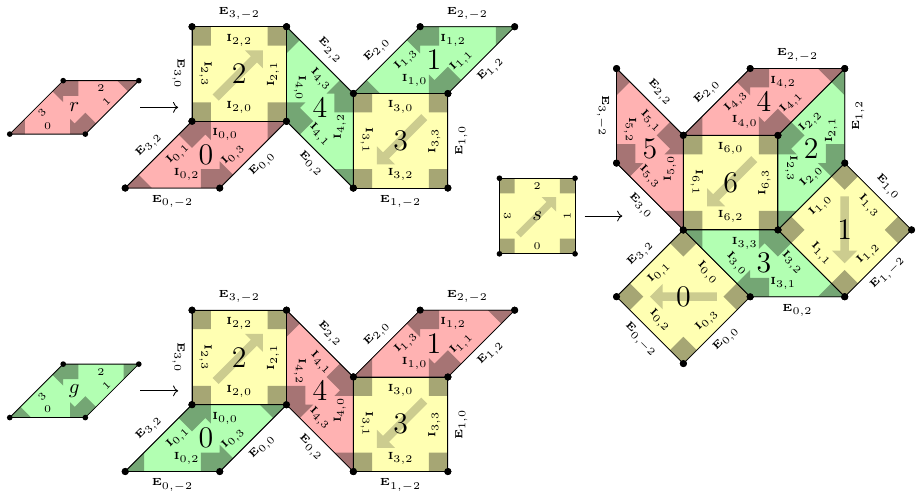}
\caption{Adjacency maps for the Ammann-Beenker tiling, with the square whole}
\label{ammann-beenker-whole-adjmap}
\end{figure}

\begin{figure}
\centering
\includegraphics*[alt={Diagrams similar to the deflation maps in the previous figure, with the detailed edge labels omitted since they are unchanged, but with each tile cloned into multiple subtypes, distinguished by subscripts and slightly different colours, showing which subtypes are deflated from which other ones.}]{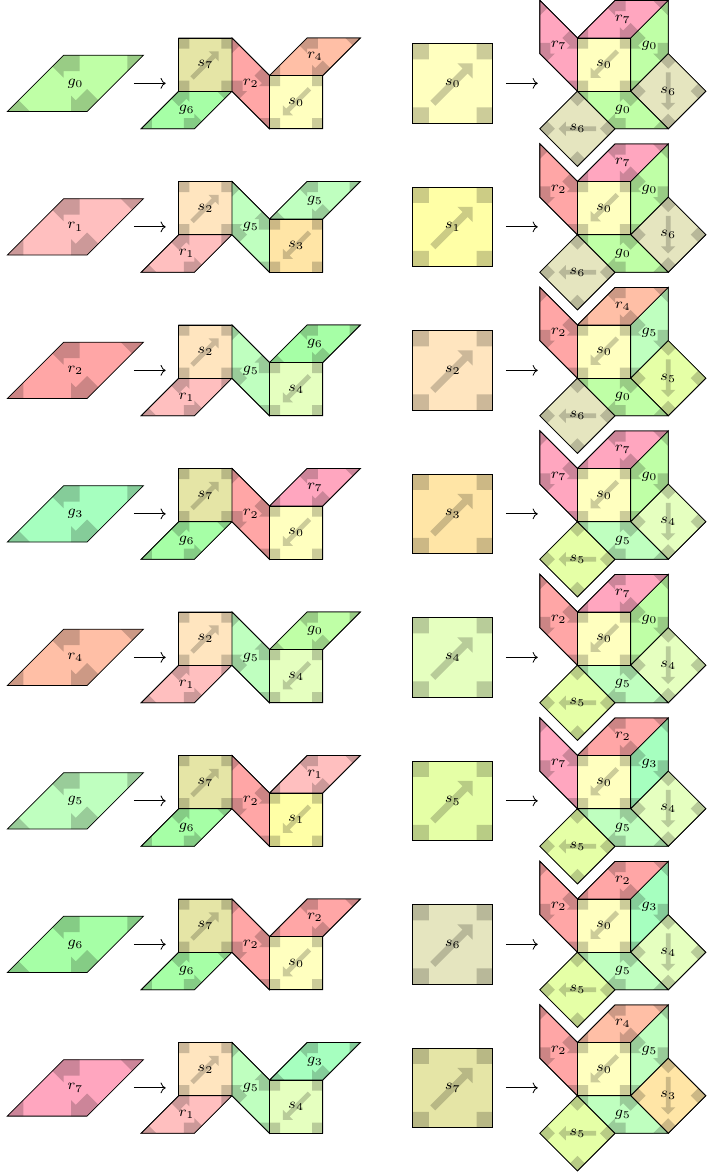}
\caption{Unambiguous refinement of the Ammann-Beenker system, with the square whole}\label{ammann-beenker-whole-refined}
\end{figure}

\begin{figure}
\centering
\includegraphics*[alt={Diagrams of the Spectre and Mystic tiles. The edges are marked with two different types, via colours and arrows. A faint line within the Mystic shows how it ultimately decomposes into two Spectres.}]{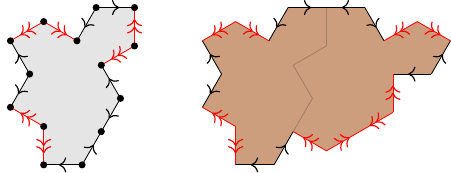}
\caption{Tile diagrams for the Spectre/Mystic system, with labelled
  edges} \label{spectreh7h8tiles}
\end{figure}

\begin{figure}
\centering
\includegraphics*[alt={A detailed map of the patch of Spectres (and one Mystic) obtained by deflating a single Spectre tile, marked up with symbols on both sides of each edge.}]{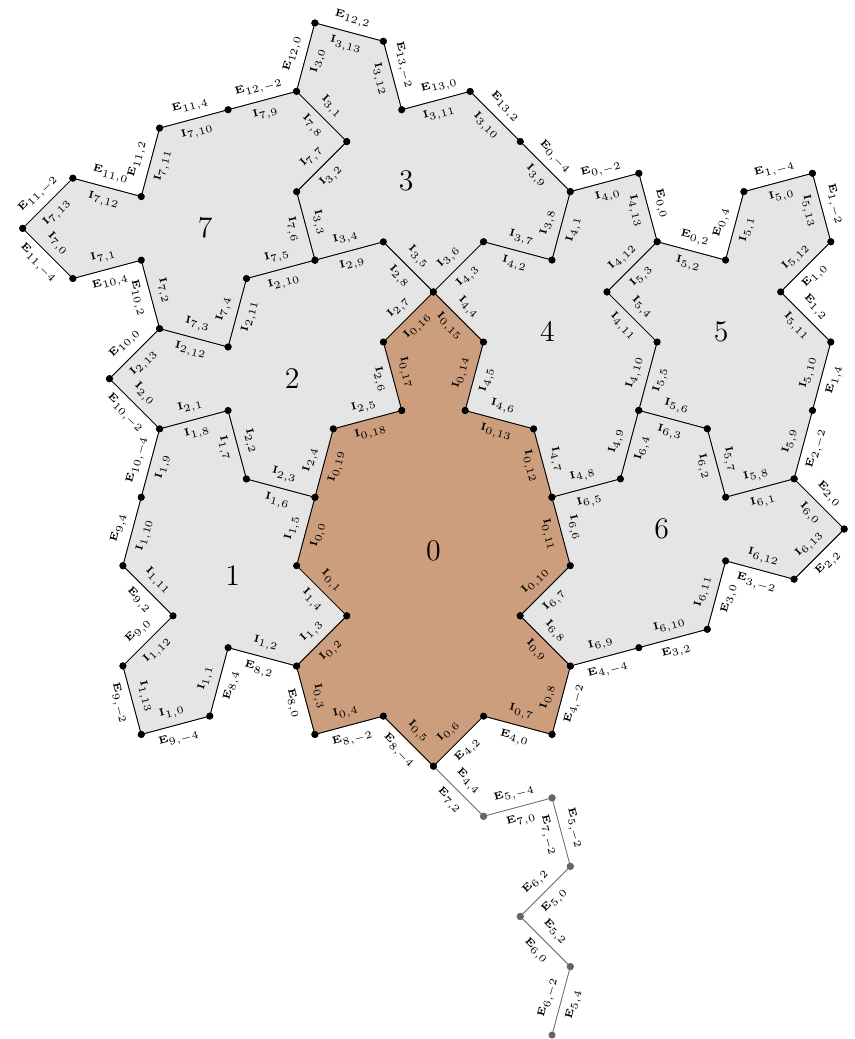}
\caption{Adjacency map for deflating a Spectre to a patch of
  Spectres}
\label{spectreh8}
\end{figure}

\begin{figure}
\centering
\includegraphics*[alt={Mapping rules showing how each of the edge types from the previously shown diagram of Spectre and Mystic tiles is transformed into a non-straight path of other edges during deflation.}]{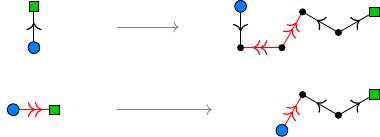}
\caption{Edge deflation rules for the Spectre/Mystic substitution system}
\label{spectreh7h8rules}
\end{figure}

\begin{figure}
\centering
\includegraphics*[alt={A detailed map of the patch of Spectres (and one Mystic) obtained by deflating a single Mystic tile, marked up with symbols on both sides of each edge. In this diagram two zero-thickness spurs are shown separately, which would otherwise cause overlapping edges and make the diagram too hard to read. The vertex where each spur should start is indicated by a matching symbol on the spur and in the main diagram: a blue square and a red diamond.}]{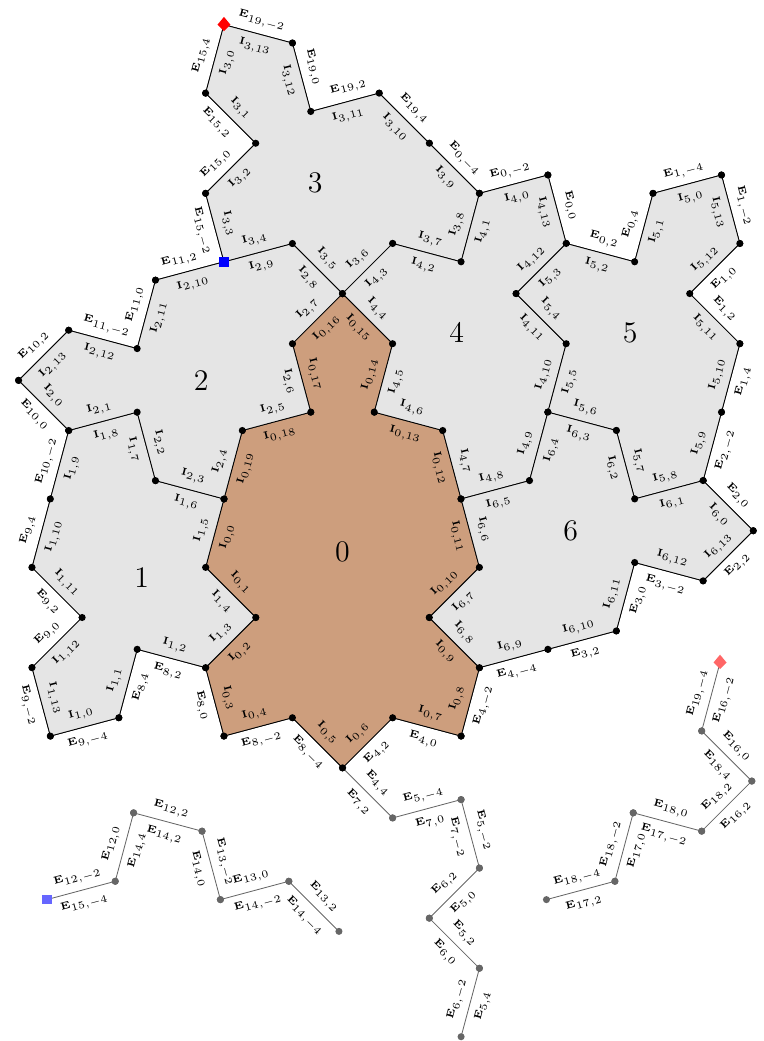}
\caption{Adjacency map for deflating a Mystic to a patch of Spectres
  (two spurs are shown displaced from their true locations to avoid
  overlapping text)} \label{spectreh7}
\end{figure}

\begin{figure}
\centering
\includegraphics*[alt={A large patch of Spectre tiling, most of which has 6-way rotational symmetry, except for a few paths of Spectres which are asymmetric, shown in a darker shade. Two regions at the top of the diagram are coloured red and yellow and outlined by thick boundaries, pointing out that those two regions are exactly identical, with even the shaded paths of Spectres within them agreeing.}]{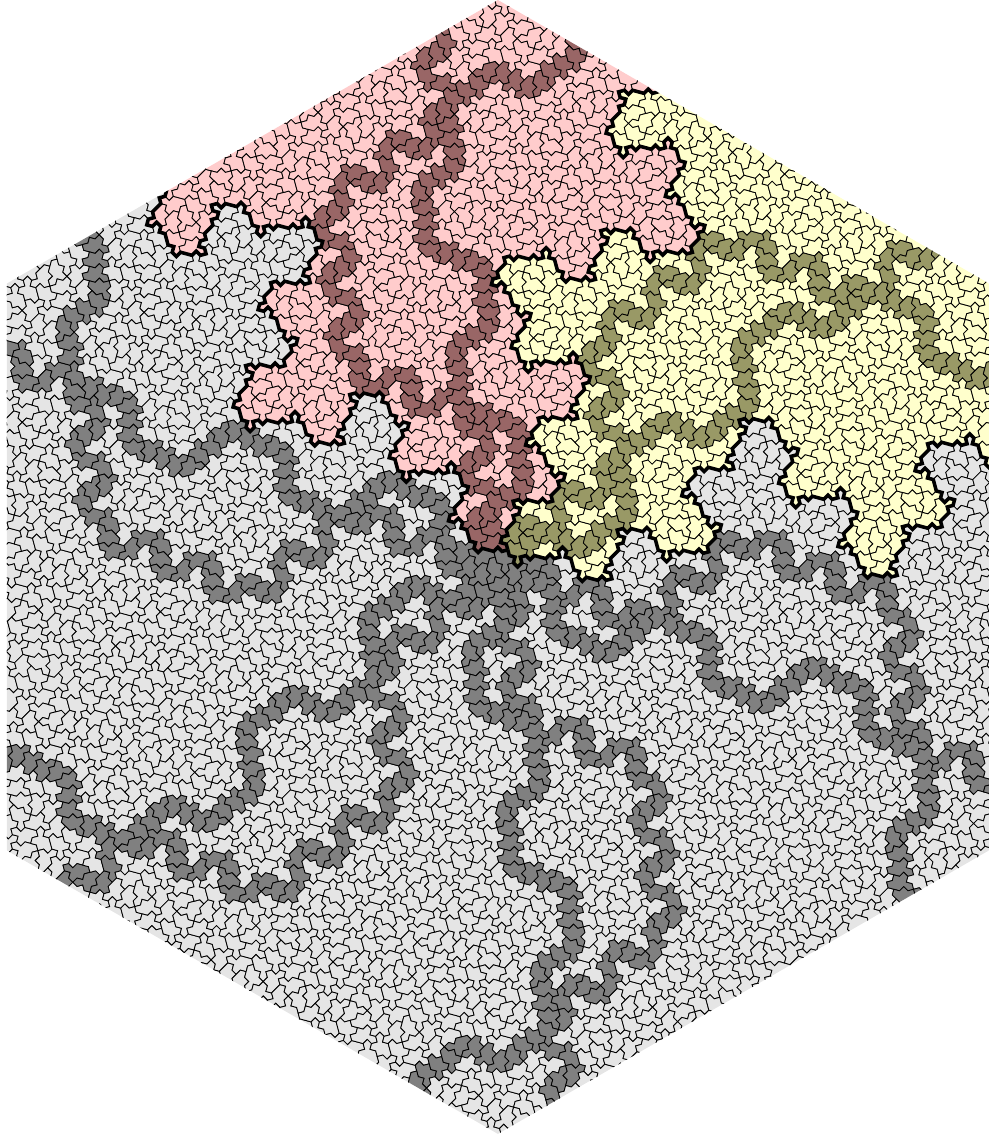}
\caption{The singular Spectre tiling containing two congruent infinite
  supertiles $(\Pi,\Xi)$ and $(\Xi,\Pi)$ (red and yellow, top). The
  light tiles are invariant under any rotation by a multiple of
  $60^\circ$.} \label{spectre-px}
\end{figure}

\begin{figure}
\centering
\includegraphics*[alt={Diagrams of the hat and double-hat tiles. The edges are marked with two different types, via colours and arrows. A faint line within the double-hat shows how it ultimately decomposes into two hats, one of which is reflected.}]{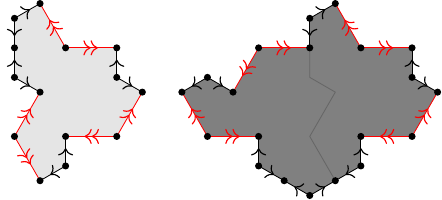}
\caption{Marked hat tiles in the $\{H_7,H_8\}$ substitution
  system} \label{hath7h8tiles}
\end{figure}

\begin{figure}
\centering
\includegraphics*[alt={A detailed map of the patch of hats (and one double-hat) obtained by deflating a single hat tile, marked up with symbols on both sides of each edge.}]{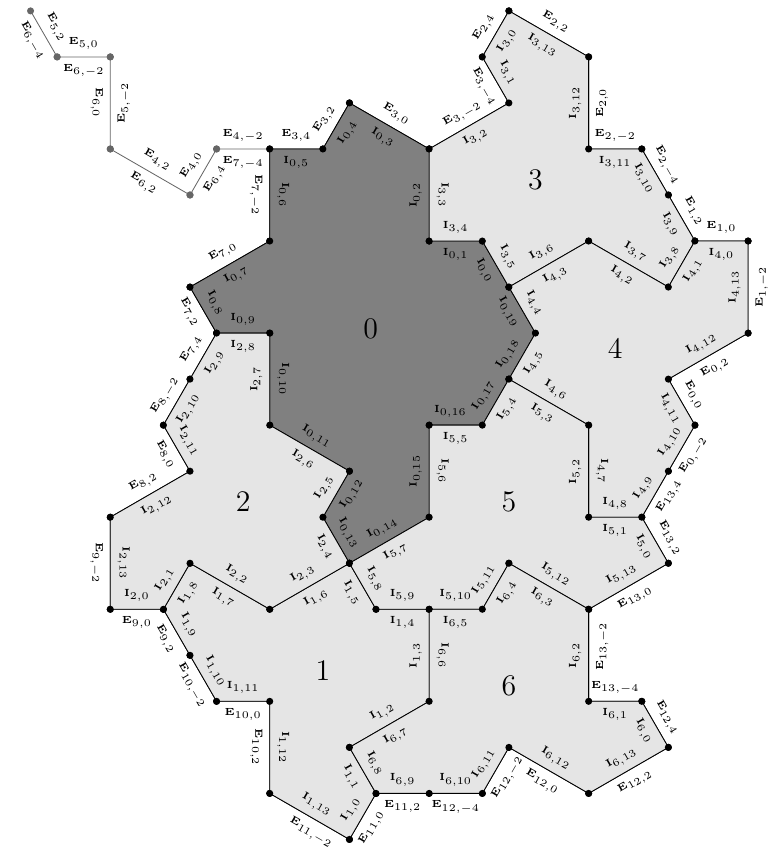}
\caption{Adjacency map for deflating a single hat to an $H_8$
  patch} \label{hath8}
\end{figure}

\begin{figure}
\centering
\includegraphics*[alt={Mapping rules showing how each of the edge types from the previously shown diagram of hat and double-hat tiles is transformed into a non-straight path of other edges during deflation.}]{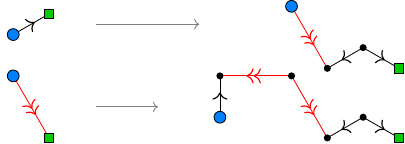}
\caption{Edge deflation rules for the $\{H_7,H_8\}$ substitution
  system} \label{hath7h8rules}
\end{figure}

\begin{figure}
\centering
\includegraphics*[alt={A detailed map of the patch of hats (and one double-hat) obtained by deflating a single double-hat tile, marked up with symbols on both sides of each edge. In this diagram two zero-thickness spurs are shown separately, which would otherwise cause overlapping edges and make the diagram too hard to read. The vertex where each spur should start is indicated by a matching symbol on the spur and in the main diagram: a blue square and a red diamond.}]{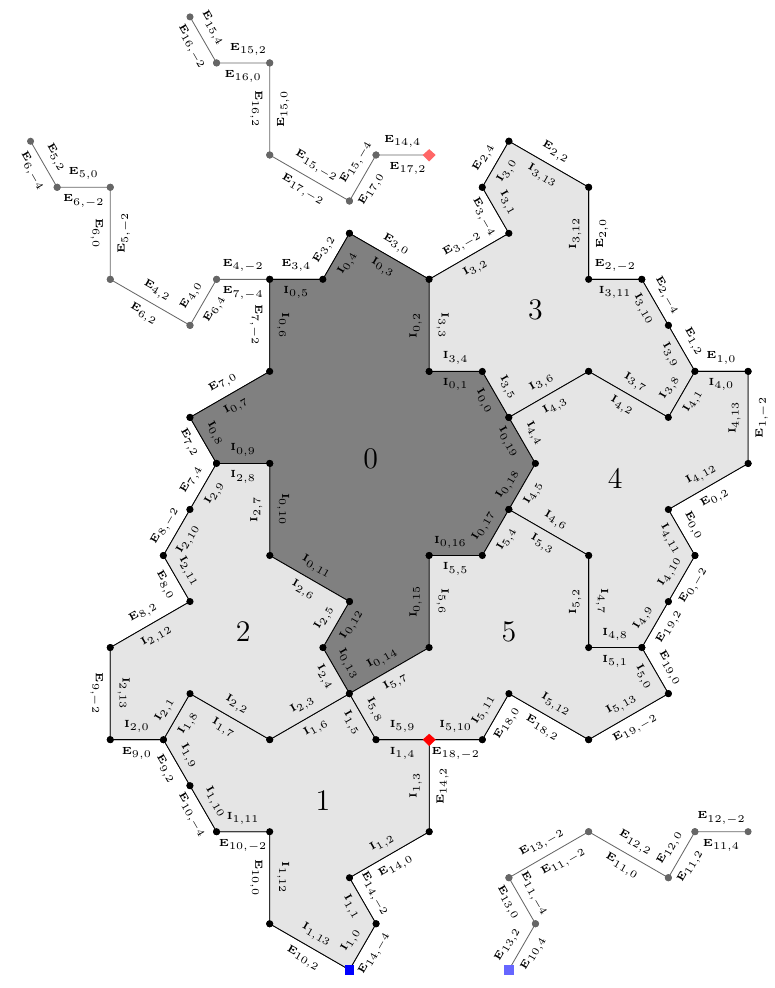}
\caption{Adjacency map for deflating a double hat to an $H_7$ patch
  (two spurs are shown displaced from their true locations to avoid
  overlapping text)} \label{hath7}
\end{figure}

The unambiguous version of the $\{H_7,H_8\}$ system derived by the
same procedure involves seven variants of the single hat tile, and
still just one variant of the double-hat. These are shown in
Figure~\ref{hat8}.

As in the previous section, we can apply the algorithm in
section~\ref{ambiguityfinding} to each of these refined systems, to
discover the reason why each system was originally ambiguous.

For the 8-tile refinement of the $\{H_7,H_8\}$ system, the answer is
extremely similar to the answer for the Spectre/Mystic system,
described in the previous section. Again, there is just one infinite
supertile in the unrefined system which has three distinct
representations in the refined system. Moreover, again, two of those
are neighbours within the same tiling of the plane, rotated $60^\circ$
relative to each other -- and, also just like the Spectre case, that
tiling has 6-way rotational near-symmetry. The third identical
infinite supertile occurs in a distinct tiling of the plane, but a
very similar one.

The unrefined address of a tile in this infinite supertile consists of
a hat which is subtile number 6 of a larger hat, which is subtile 6 in
turn of the next larger hat, and so on for ever:
$$\hbox{hat}\xleftarrow{6}\hbox{hat}\xleftarrow{6}\hbox{hat}\xleftarrow{6}\hbox{hat}\xleftarrow{6}\cdots$$
In the refined system, the three corresponding addresses each have the
same type of hat at every layer (unlike the Spectre case, where two of
the addresses alternated tile types):
$$\begin{aligned}
h_2\xleftarrow{6}h_2\xleftarrow{6}h_2\xleftarrow{6}h_2\xleftarrow{6}\cdots \\
h_3\xleftarrow{6}h_3\xleftarrow{6}h_3\xleftarrow{6}h_3\xleftarrow{6}\cdots \\
h_4\xleftarrow{6}h_4\xleftarrow{6}h_4\xleftarrow{6}h_4\xleftarrow{6}\cdots \\
\end{aligned}$$

\begin{figure}
\centering
\includegraphics*[alt={A large patch of hat tiling, most of which has 6-way rotational symmetry, except for a few paths of hats which are asymmetric, shown in a darker shade. Two regions at the top of the diagram are coloured red and yellow and outlined by thick boundaries, pointing out that those two regions are exactly identical, with even the shaded paths of hats within them agreeing.}]{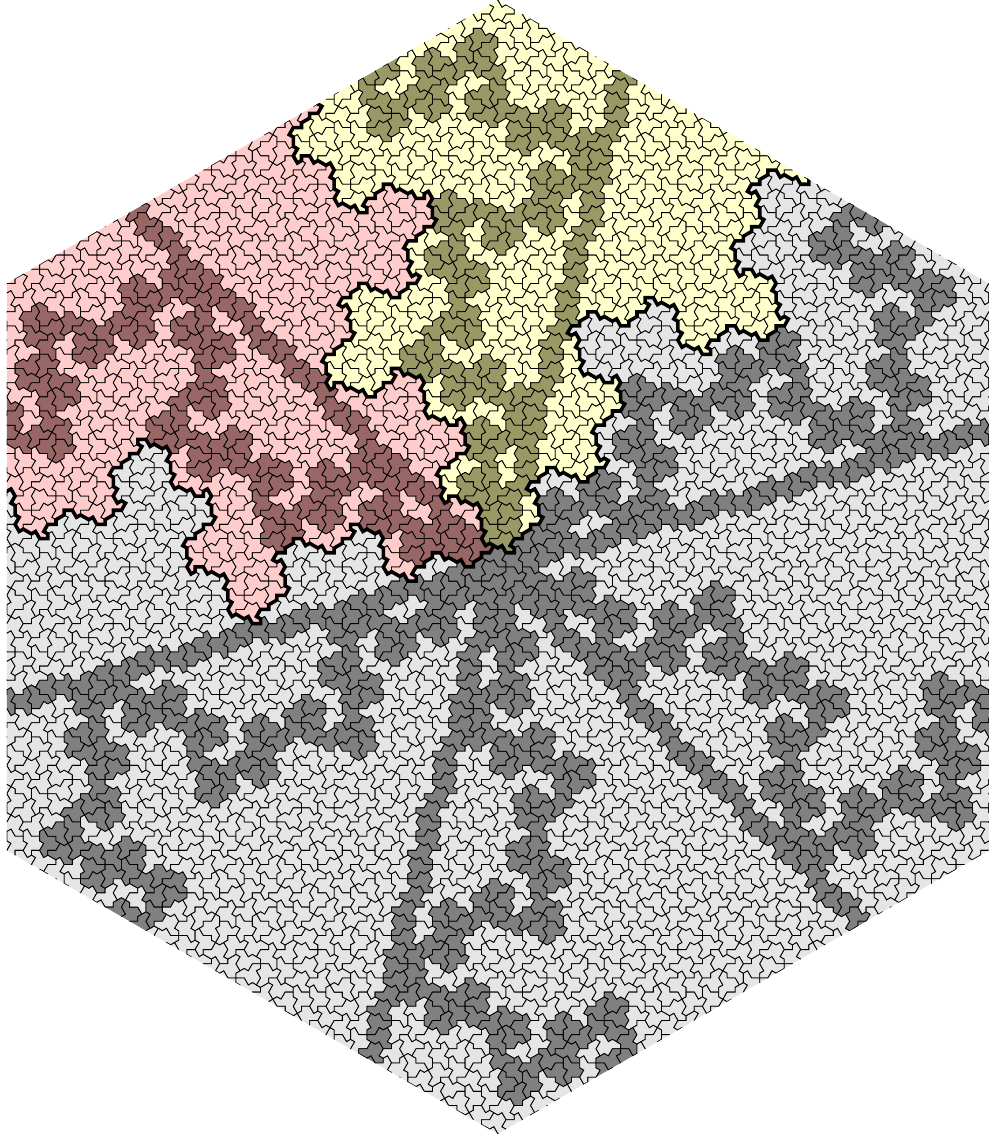}
\caption{The singular Hat tiling containing two congruent infinite
  supertiles $(h_2)$ and $(h_4)$ (red and yellow, top). The light
  tiles are invariant under any rotation by a multiple of
  $60^\circ$.} \label{hat-singular}
\end{figure}

\begin{figure}
\centering
\includegraphics*[alt={Diagrams of the HTPF hat metatiles and the two handednesses of hat. The edges of the HTPF metatiles are marked with five different types, via colours and arrows. The edges of the hat tiles themselves are coloured to distinguish the short and long edges, but have no arrows.}]{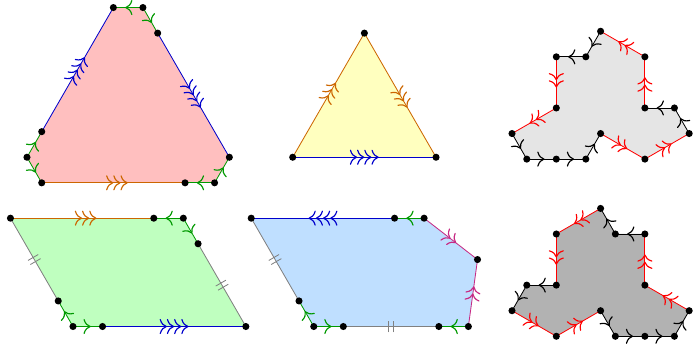}
\caption{Diagrams of the {\HTPF} metatiles and the hat, with labelled edges}
\label{htpf-prototiles}
\end{figure}

\begin{figure}
\centering
\includegraphics*[alt={Mapping rules showing how each of the HTPF edge types from the previous figure is transformed into a path of other edges during deflation to other HTPF metatiles. The purple two-arrow edge from the pointy end of the F metatile is shown at an awkward angle, because it is naturally skew to all other edge types, and this way its deflation appears at a sensible angle.}]{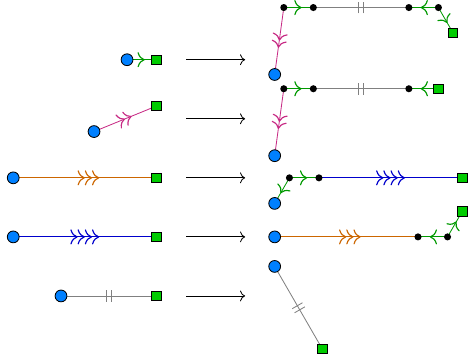}
\caption{Edge deflation rules for {\HTPF} metatiles to other metatiles}
\label{htpf-edge-meta-meta}
\end{figure}

\begin{figure}
\centering
\includegraphics*[alt={Mapping rules showing how each of the edge types from the diagram of HTPF metatiles is transformed into a path of other edges during final deflation to hats.}]{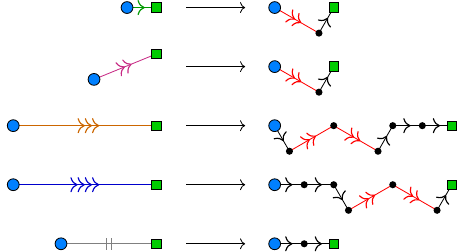}
\caption{Edge deflation rules for {\HTPF} metatiles to hats}
\label{htpf-edge-meta-hats}
\end{figure}

\begin{figure}
\centering
\includegraphics*[alt={A detailed map of each HTPF metatile and how it deflates into smaller metatiles, marked up with symbols on both sides of each edge.}]{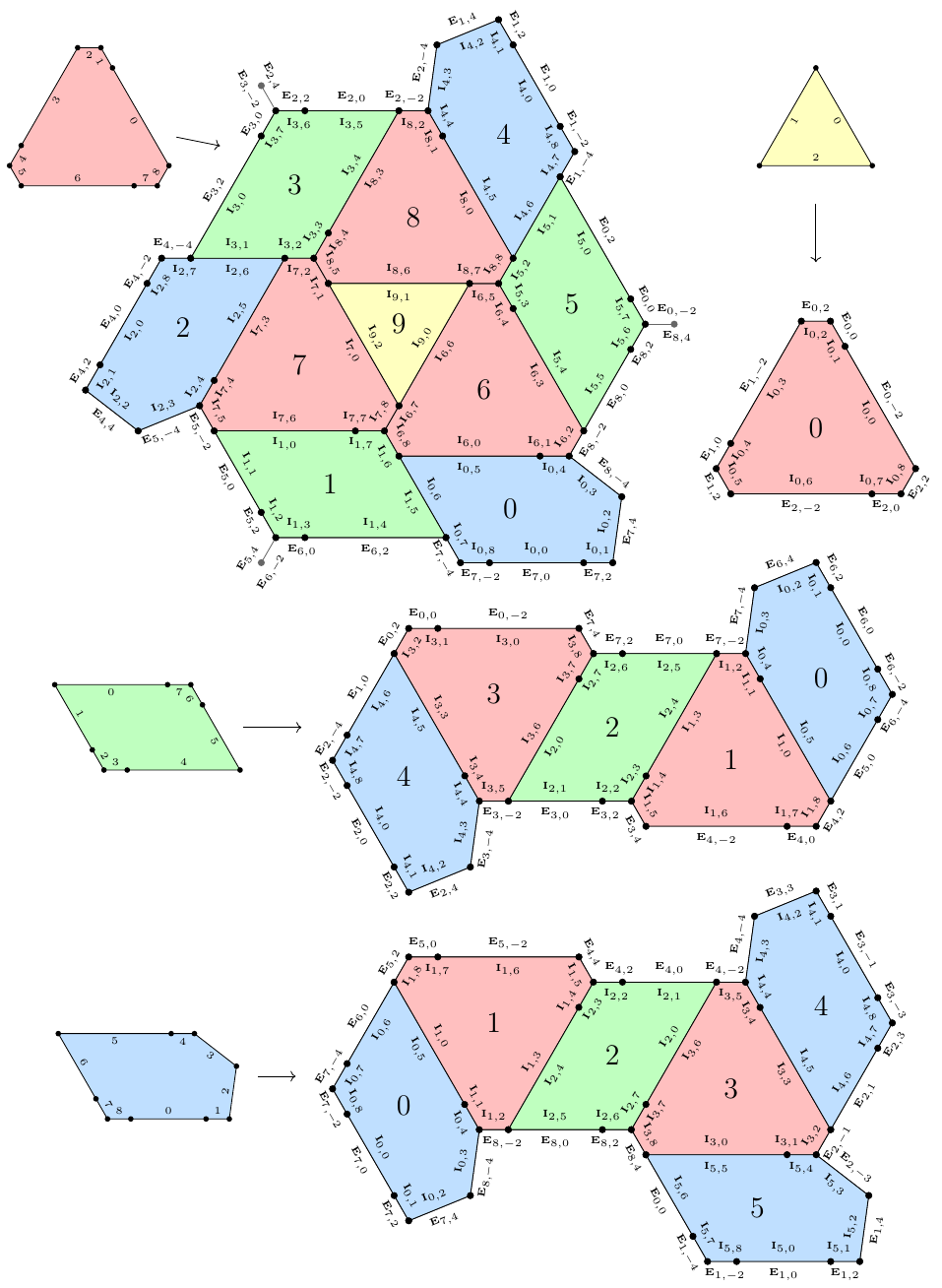}
\caption{Adjacency maps for deflating {\HTPF} metatiles to more metatiles}
\label{htpf-adjmap-meta-meta}
\end{figure}

\begin{figure}
\centering
\includegraphics*[alt={A detailed map of each HTPF metatile and how it deflates into hats, marked up with symbols on both sides of each edge.}]{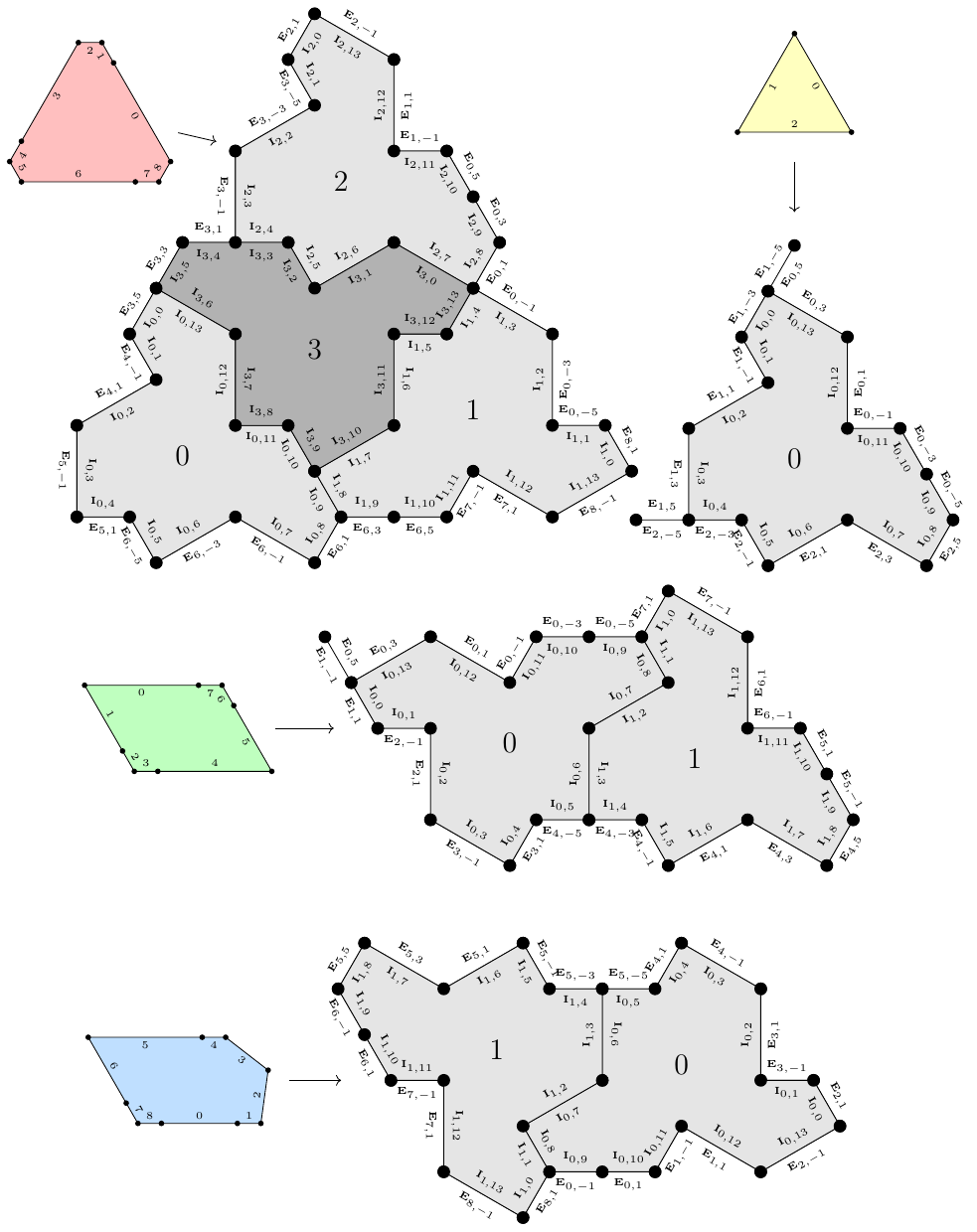}
\caption{Adjacency maps for deflating {\HTPF} metatiles to hats}
\label{htpf-adjmap-meta-hats}
\end{figure}

The two of these supertiles that appear together in the same tiling
are the ones with $h_2$ and $h_4$ tile types. The tiling itself is
shown in Figure~\ref{hat-singular}. Like the singular Spectre tiling,
the asymmetric paths come in two forms; unlike the Spectre tiling, one
of those forms is a straight line rather than a wiggly path. (If the
hats are instead drawn as turtles, the straight lines are aligned
naturally to the underlying kite tiling.)

The third copy of this infinite supertile, involving the $h_3$
address, is part of a second 6-way singular hat tiling. If the two
tilings are overlaid so that the $h_3$ supertile corresponds to the
$h_4$ supertile, then the two tilings differ only in a single straight
line of hats; if instead the two tilings are overlaid so that the
$h_3$ supertile corresponds to $h_2$, then they instead differ in a
wiggly path. \cite{Socolar_2023} exhibits these two tilings and their
two minimal differing paths.

Applying the section~\ref{ambiguityfinding} algorithm to the
{\HHTPFFF} substitution system instead of the 8-tile refinement of
$\{H_7,H_8\}$ reveals no new interesting cases. In this system, there
is again a single infinite supertile with multiple representations,
which in the unrefined {\HTPF} system has an address consisting
entirely of $F$ metatiles, with each one being subtile number 5 of its
supertile (in the numbering of Figure~\ref{htpf-adjmap-meta-meta}):
$$\hbox{hat}\xleftarrow{0}F\xleftarrow{5}F\xleftarrow{5}F\xleftarrow{5}F\xleftarrow{5}\cdots$$
This infinite supertile has just two possible representations in the
refined {\HHTPFFF} system, because both $F_0$ and $F_1$ can appear as
subtile 5 of some version of $F$, but $F_2$ cannot. For each of $F_0$
and $F_1$, subtile 5 is the same type as itself, so the two possible
addresses are
$$\begin{aligned}
\hbox{hat}\xleftarrow{0}F_0\xleftarrow{5}F_0\xleftarrow{5}F_0\xleftarrow{5}F_0\xleftarrow{5}\cdots \\
\hbox{hat}\xleftarrow{0}F_1\xleftarrow{5}F_1\xleftarrow{5}F_1\xleftarrow{5}F_1\xleftarrow{5}\cdots \\
\end{aligned}$$
These correspond exactly to the $h_3$ and $h_2$ supertiles from the
other substitution system. The $h_4$ variant of that supertile does
not appear separately in this system.

\section{Future possibilities}

One useful operation not implemented by this suite of algorithms is
the ability to translate the same tiling from one substitution system
to another, by creating a transducer that would consume an address in
the first system and emit an address in the second.

Such a transducer would not necessarily be able to keep the invariant
seen in the transducers here, of emitting exactly one output symbol
per input symbol, because not all substitution systems for the same
tiling deflate by the same factor. For example, there is a
`half-deflation' system for the P2 and P3 Penrose tilings
\cite{DEBRUIJN198139} which deflates each one to the other, such that
two deflations return to the same tiling you started with, having
performed what would normally be a full deflation.
\cite{mathBlock2025} also presents systems for the hat tiling which
take two deflations to repeat the effect of a single deflation in the
$\{H_7,H_8\}$ or {\HTPF} systems (although in this case the
intermediate tiling is still a hat tiling).

If an algorithm of this form were possible, then it might also be
possible to use it to match a substitution system against
\emph{itself}, to determine if there exist any pairs of distinct
addresses which give rise to an identical tiling of the plane. For
example, the system in \cite{mathBlock2025} has such pairs, because in
some cases, half-inflating a perfectly symmetric instance of the hat
tiling leads to a `singular' (nearly symmetric) instance, and so the
tile addresses in the symmetric tiling must break its symmetry in
order to specify which way round the not-quite-symmetric singular
deflation will be.

\section{Implementation}

An implementation of the algorithms in this article, written in Rust,
is in the \href{https://crates.io/}{\texttt{crates.io}} package repository under the name
\href{https://crates.io/crates/substitution-tiling-transducers}{\texttt{substitution-tiling-transducers}}.

\section{Acknowledgments}

I am grateful to Robin Houston, Craig Kaplan, Pieter Mostert and
mathBlock for useful discussions and encouragement.

\clearpage

\clearpage

\clearpage

\clearpage

\clearpage

\clearpage

\clearpage

\clearpage

\clearpage

\clearpage

\clearpage

\clearpage

\clearpage

\clearpage

\clearpage

\clearpage

\clearpage

\clearpage

\bibliography{writeup.bib}

\end{document}